\documentclass[10pt]{amsart}
\pdfoutput=1
\usepackage{ mathrsfs }
\UseRawInputEncoding
\usepackage{amsmath, amsthm, amssymb,slashed,stmaryrd}
\usepackage{mathtools}
\usepackage{ifpdf}
\usepackage[pdftex]{graphicx}
\usepackage{tikz}
\usetikzlibrary{}
\usetikzlibrary{decorations.markings,matrix,arrows,decorations.pathmorphing,backgrounds,patterns,shapes.geometric,calc}
\newcommand*{\halfway}{0.5*\pgfdecoratedpathlength+.5*3pt}
\usetikzlibrary{matrix,arrows,calc}
\usepackage[all]{xy}
\UseRawInputEncoding
\usepackage{tikzit}
\tikzstyle{new edge style 0}=[-, dashed=]

\usepackage[pdftex,plainpages=false,hypertexnames=false,pdfpagelabels]{hyperref}
\long\def\todo#1{{\color{red} {#1}}}

\usepackage{amsthm}

\theoremstyle{definition}
 \newtheorem{thm}{Theorem}[section]
    
 \newtheorem{cor}[thm]{Corollary}
  \newtheorem{hyp}[thm]{Hypothesis}
 
 \newtheorem{prop}[thm]{Proposition}
  
 \newtheorem{conj}[thm]{Conjecture}
 \newtheorem{defn}[thm]{Definition}
 
 \newtheorem{ex}[thm]{Example}
  
 \newtheorem*{thm*}{Theorem}
 \theoremstyle{remark}
 \newtheorem{rmk}[thm]{Remark}

\def\beq{\begin{eqnarray}}
\def\eeq{\end{eqnarray}}
 \newcommand{\bp}{\begin{proof}[Proof]}
 \newcommand{\ep}{\end{proof}}

\DeclareSymbolFont{bbold}{U}{bbold}{m}{n}
\DeclareSymbolFontAlphabet{\mathbbold}{bbold}
\def\one{\mathbbold{1}}

\def\QFT{{\sf QFT}}
\def\qft{{\sf qft}}
\def\sqft{{\sf sqft}}

\def\MF{{\rm MF}}
\def\Fer{{\rm Fer}}
\def\sdim{{\rm sdim}}
\def\cH{{\mathcal{H}}}
\def\Ham{{\sf Ham}}
\def\sHam{{\sf sHam}}
\def\qm{{\sf QM}}
\def\sqm{{\sf SQM}}
\def\PU{{\rm PU}}

\def\U{{\rm U}}
\def\dom{{\rm dom}}

\def\Tr{{\rm Tr}}

\def\Pf{{\sf{Pf}}}

\def\s{{\sf s}}

\def\rC{{\rm C}}
\def\rT{{\rm T}}
\def\sL{{\boldsymbol\cap}}
\def\sR{{\boldsymbol\cup}}
\def\sI{{\sf I}}
\def\sS{{\sf S}}

\def\Sym{{\rm Sym}}
\def\Vect{{\sf Vect}}
\def\spt{{\rm spt}}

\def\Ra{{\sf R}}
\def\NS{{\sf NS}}

\def\Bord{\hbox{{\sf Bord}}}

\def\sTr{{\sf sTr}}

\def\O{{\rm O}}
\def\Pin{{\rm Pin}}
\def\KO{{\rm KO}}

\def\K{{\rm K}}

\def\Th{{\rm Th}}

\def\Cl{{\rm Cl}}
\def\cCl{\mathbb{C}{\rm l}}

\def\cl{{\rm cl}}

\def\H{{\rm H}}
\def\HH{{\mathbb H}}
\def\Spin{{\rm Spin}}

\def\U{{\rm U}}

\def\SO{{\rm SO}}

\def\TA{{\sf Alg}}

\def\Map{{\sf Map}}

\def\TMF{{\rm TMF}}

\def\cG{\mathcal{G}}

\def\X{\mathfrak{X}}

\def\EBord{{\sf EBord}}

\def\Pf{{\sf Pf}}

\def\Pic{{\rm Pic}}

\newcommand{\nsq}{{\sq^{\!\nabla}}\!}

\def\pt{{\rm pt}}

\def\bS{{\mathbb{S}}}
\def\A{{\mathbb{A}}}

\def\R{{\mathbb{R}}}

\def\id{{{\rm id}}}
\def\C{{\mathbb{C}}}
\def\Q{{\mathbb{Q}}}
\def\Z{{\mathbb{Z}}}

\def\End{{\rm End}}
\def\SL{{\rm SL}}
\def\MP{{\rm MP}}

\def\Sym{{\sf Sym}}

\newcommand{\op}{{\sf{op}}}   
\newcommand{\sq}{/\!\!/}

\def\leftsquigarrow{\ensuremath{\rotatebox[origin=c]{-180}{$\rightsquigarrow$}}}

\def\torus{\!\!\hbox{\begin{tikzpicture}[baseline=(basepoint)];
\node (A) at (0,0) {};
\node (AA) at (-.25,0) {};
\node (B) at (0,.1) {};
\node (BB) at (-.25,.1) {};
\node (C) at (0,-.3) {};
\node (CC) at (-.25,-.3) {};
\node (D) at (0,-.4) {};
\node (DD) at (-.25,-.4) {};
\draw[-,bend right=90] (A) to (C);
\draw[-,bend left=90] (AA) to (CC);
\draw[-,bend right=90] (B) to  (D);
\draw[-,bend left=90] (BB) to  (DD);
\draw (-.12,.046) ellipse (.02 and .05);
\draw (-.12,-.346) ellipse (.02 and .05);
\path (-.4,-.2) coordinate (basepoint);
\end{tikzpicture}}
}

\def\upmac{\!\!\rotatebox[origin=c]{-90}{\hbox{\begin{tikzpicture}[baseline=(basepoint)];
\node (A) at (0,0) {};
\node (B) at (0,.1) {};
\node (C) at (0,-.3) {};
\node (D) at (0,-.4) {};
\draw[-,bend right=80] (A) to (C);
\draw[-,bend right=80] (B) to  (D);
\draw (-.12,.032) ellipse (.02 and .05);
\draw (-.12,-.332) ellipse (.02 and .05);
\path (-.4,-.15) coordinate (basepoint);
\end{tikzpicture}}
}}

\def\downmac{\!\!\rotatebox[origin=c]{90}{\hbox{\begin{tikzpicture}[baseline=(basepoint)];
\node (A) at (0,0) {};
\node (B) at (0,.1) {};
\node (C) at (0,-.3) {};
\node (D) at (0,-.4) {};
\draw[-,bend right=80] (A) to (C);
\draw[-,bend right=80] (B) to  (D);
\draw (-.12,.032) ellipse (.02 and .05);
\draw (-.12,-.332) ellipse (.02 and .05);
\path (-.4,-.3) coordinate (basepoint);
\end{tikzpicture}}
}}

\newcommand\nc{\newcommand}
\nc\mf\mathfrak
\nc\mc\mathcal
\nc\mb\mathbb

\begin{document}

\title[Elliptic cohomology and quantum field theory]{Elliptic cohomology and quantum field theory}
\author{Daniel Berwick-Evans}
\thanks{Department of Mathematics, University of Illinois Urbana-Champaign, \url{danbe@illinois.edu}}
\thanks{The author gratefully acknowledges support from NSF grants DMS-2205835 and DMS-2340239.}


\maketitle 

\setcounter{tocdepth}{1} 
\tableofcontents

\section{Introduction}

In the late 1980s, Witten used methods from theoretical physics to construct an elliptic genus valued in integral modular forms~\cite{Witten_Dirac}. His arguments suggested a much deeper connection between elliptic cohomology and quantum field theory. Pursuing these ideas, Segal formulated a precise definition of 2-dimensional quantum field theory as a projective representation of a bordism category~\cite{SegalCFT}. During the same period,  Hopkins and collaborators grew the subject of elliptic cohomology, culminating in the construction of topological modular forms (TMF) and a string orientation lifting Witten's genus to a map of ring spectra~\cite{HopkinsICM2002,Lurie_Elliptic}.
Through a variety of insights and efforts spanning several  decades, a rich picture is emerging that connects TMF with the geometry of field theories. In Stolz and Teichner's framework, the basic conjecture can be summarized as follows~\cite{ST11}.

\begin{conj}\label{conj}
For each $n\in \Z$ and smooth manifold $X$ there exist dashed arrows 
\beq
\begin{tikzpicture}[baseline=(basepoint)]
    \node (A) at (0,0) {$  \left\{\begin{array}{c}  {\rm families\ of\ } \\ {n\hbox{-}{\rm dimensional\ string }}\\ {\rm  \ manifolds } \ M\to X \end{array}\right\}$};  
  \node (B) at (8,-1) {$\TMF^{-n}(X)$};
  \node (C) at (0,-2.5) {$\left\{\begin{array}{c}  {\rm fully\ extended,\ degree} \ -\!n \\ 2\hbox{-}{\rm dimensional} \\  {\rm supersymmetric\ field }\\ {\rm  theories\ over\ } X\end{array}\right\}$};
  \draw[->] (A) to node [above=4pt] {${\rm topological\ index}$} (B);
  \draw[->,dashed] (A) to node [left] {${\rm quantize}$} (C);
  \draw[->>,dashed] (C) to node [below] {${\rm cocycle}_{-n}$} (B);
    \path (0,-1) coordinate (basepoint);
    \end{tikzpicture}\label{eq:conj}
\eeq
making the triangle commute, where ``topological index" is the string orientation, ``quantize" is quantization of the supersymmetric $\sigma$-model on the fibers of $M\to X$, and ``${\rm cocycle}_{n}$" sends a degree~$n$ field theory over $X$ to a class in~$\TMF^{n}(X)$. For negative $n$, the empty manifold is the only string manifold and the content of~\eqref{eq:conj} is the existence of a cocycle map.
\end{conj}

In brief, Conjecture~\ref{conj} proposes a geometric construction of TMF whose cocycles are 2-dimensional quantum field theories. It also proposes a $\TMF$-index theorem: the composition of dashed arrows for $n\ge 0$ defines an analytic index which the triangle~\eqref{eq:conj} equates with the topological index constructed by Ando--Hopkins--Rezk~\cite{AHR}. Stolz and Teichner have explained how this proposed $\TMF$-index theorem generalizes Atiyah and Singer's~\cite{ST04}; see Conjecture~\ref{thm11} below.

\subsection{Why the conjecture is interesting}
The central role of topological K-theory in 20th century mathematics comes from the way it intermingles with many disciplines. In geometry, vector bundles determine K-theory classes. In analysis, the space of Fredholm operators represents K-theory. In algebra, representations of compact Lie groups determine (equivariant) K-theory classes. Many important results rely on these interconnections, with the Atiyah--Singer index theorem being the crowning example that interfaces analysis, geometry, and topology. 
Within algebraic topology, TMF is considerably deeper than K-theory, e.g., $\pi_*\TMF$ detects many more classes in the stable homotopy groups of spheres. 
However, the absence of a geometric description has largely confined the reach of TMF to the realm of homotopy theory. Resolving Conjecture~\ref{conj} would change this landscape. 


\vspace{.025in}

Complementing this motivation from algebraic topology, a current theme in theoretical physics is to understand various spaces of quantum field theories, including the landscape problem in string theory~\cite{Susskind,Douglas_space} and topological phases of matter~\cite{Kitaev,FreedMoore}. Conjecture~\ref{conj} asserts that one such space of quantum field theories has a very interesting topology. The quantum theories in question have the minimal supersymmetry in dimension~2; this imposes relatively meager constraints on observables, making the standard computations difficult. Physical reasoning assuming the validity of Conjecture~\ref{conj} has led to a slew of surprising predictions, many of which connect back with mathematics, e.g., see~\cite{GJF2,GPPV,Tachikawa,TachikawaYamashita,TachikawaYamashita2,TYY,LinPei}. To mention two specific examples, the abundant torsion in TMF indicates new deformation invariants of quantum field theories that refine elliptic genera~\cite{GJFW}, and the $576$-periodicity of TMF predicts a surprising periodicity in the chiral free fermions~\cite{DouglasHenriques}. All told, there is an expectation that resolving Conjecture~\ref{conj} will provide important mathematical tools for studying spaces of quantum field theories in general, with deep insight into the structure of 2-dimensional field theories and string theory in particular~\cite{Witten_Dirac}. 

\vspace{.025in}

There are also anticipated applications of~\eqref{eq:conj} outside of homotopy theory and theoretical physics. We highlight three areas. 

\vspace{.025in}

In differential geometry, Stolz conjectured that the Witten genus obstructs the existence of a positive Ricci curvature metric~\cite{stolz_conj}. The analogous $\hat{A}$-genus obstruction to positive scalar curvature uses the geometry of Dirac operators~\cite{Lichnerowicz}. In a related vein, Hopkins sketched the detection of exotic smooth structures on free loop spaces via a $\TMF$-generalization of the signature operator \cite[Ch.~10]{DFHH} suggested by Witten~\cite{Witten_Dirac}. This flavor of loop space index theory remains out of reach without a solid handle on the geometry of field theories in~\eqref{eq:conj}. 

\vspace{.025in}


In algebra, lattice vertex operator algebras can encode sporadic groups, most famously through Borcherd's proof of monstrous moonshine~\cite{Borcherds}. The $\theta$-function of a lattice also determines a topological modular form~\cite[\S5]{HopkinsICM2002}. The cocycle map in~\eqref{eq:conj} is expected to clarify these and other coincidences surrounding the monster group, 2-dimensional conformal field theories, and elliptic cohomology~\cite{Moravamoon,TheoTMM,TMFmoon,TMFmoon2}. 

\vspace{.025in}

In representation theory, Grojnowski described elliptic generalizations of standard examples built from equivariant cohomology and equivariant K-theory~\cite{Grojnowski}. Pursuing these ideas using complex analytic equivariant elliptic cohomology (roughly, equivariant $\TMF\otimes \C$) has lead to some incredibly rich objects \cite{GKV,ZhaoZhong,YangZhao,FRV,RTV,AganagicOkounkov,RSVZ,GanterRam,KRWell,RWschubert,RWell}.
In these constructions---in contrast to $\K$-theory and ordinary cohomology---the relevant equivariant elliptic cohomology classes do not have a well-understood geometric description. Correspondingly, the true origin of these elliptic-themed algebras remains a bit opaque. It seems inevitable a better understanding of~\eqref{eq:conj} would provide further insights
 into these deep objects, e.g., via integral refinements.

\begin{rmk}
There is an important enhancement of Conjecture~\ref{conj}, stating that \emph{concordance classes} of degree~$n$ field theories over $X$ are in bijection with $\TMF^n(X)$. In other words, $\TMF$ is a \emph{complete} deformation invariant of 2-dimensional supersymmetric field theories. This is strictly more difficult than Conjecture~\ref{conj} and introduces technicalities we wish to avoid in this survey. For a brief discussion we refer to Remark~\ref{rmk:weirdsheaf2} below. 
\end{rmk}




%



\subsection{Frameworks for supersymmetric quantum field theory}
A reoccurring theme in mathematical quantum field theory is the interplay between the analysis of operator algebras and the geometry of spacetime. For example, Jones's contributions---as further amplified by Witten---used quantum field theory to forge a deep connection between the geometry of knots and the analysis of von Neumann algebras \cite{Jones,WittenJones,AtiyahJonesWitten}. It is expected that Conjecture~\ref{conj} will continue on this path. An influential language designed to probe this interface of geometry and analysis is due to Atiyah~\cite{AtiyahTQFT} and Segal~\cite{SegalCFT}: a $d$-dimensional quantum field theory is a symmetric monoidal functor from a category of $d$-dimensional bordisms to topological vector spaces,
\beq\label{eq:SegalQFT}
(\Bord_d,\coprod)\to (\Vect,\otimes).
\eeq
Geometry is in the source of~\eqref{eq:SegalQFT} and functional analysis is in the target. 
For example, in dimension~2 the source can be chosen to capture moduli spaces of Riemann surfaces (as in Segal's original work on conformal field theory), and the target can be chosen to encode the analysis of nuclear Fr\'echet spaces under the projective tensor product. 

\vspace{.05in}

In addition to~\eqref{eq:SegalQFT}, there are now many schools of thought on mathematical quantum field theory, each bringing their own insights to~\eqref{eq:conj}, e.g., vertex algebras~\cite{CDO1,HuKriz,GanterLaures}, factorization algebras~\cite{costello_WG1,gradygwilliam}, and conformal nets~\cite{Wassermann,CN1}. Ultimately, development of each of these languages as well as dictionaries that compare them will quicken progress on Conjecture~\ref{conj} and other outstanding problems in mathematical quantum field theory. For concreteness, below we will examine Conjecture~\ref{conj} through the lens of~\eqref{eq:SegalQFT} following Stolz and Teichner~\cite{ST11}. 

\vspace{.05in}

To begin, for a \emph{rigid geometry}~$\cG$ and a smooth stack~$\X$, there is a symmetric monoidal category $\Bord_\cG(\X)$ of \emph{bordisms with $\cG$-structure over $\X$}~\cite[Definitions~2.21, 2.46, 2.48, and 4.4]{ST11}. The bordisms are the \emph{spacetime} or \emph{worldsheet} on which a quantum field theory is defined, and the stack $\X$ contributes \emph{background fields}. In the setting of Conjecture~\ref{conj}, background fields are simply maps to the smooth manifold $X$. However, it will be important to allow for more general background fields such as $G$-bundles with connection, i.e., background gauge fields. 

\vspace{.05in}

The following two classes of examples give a flavor for the categories $\Bord_\cG(\X)$.

\begin{ex}[Euclidean bordism categories]\label{ex:Eucrigid}
A prototypical rigid geometry is flat spin geometry on $\R^d$, characterized by the isometry group $\R^d\rtimes \Spin(d)$ acting on $\R^n$ by translations and orientation-preserving rotations (through the double cover $\Spin(d)\to \SO(d)$). In this case we write $\EBord_d(\X)$ for the corresponding \emph{$d$-dimensional Euclidean bordism category over~$\X$}. The morphisms in $\EBord_d(\X)$, denoted $\Phi\colon \Sigma\to \X$, are compact spin $d$-manifolds $\Sigma$ with flat metric and geodesic boundary together with a map $\Phi\colon \Sigma\to \X$. The source and target of a bordism comes from restricting $\Phi$ to a partition of the boundary $\partial\Sigma=\partial_{\rm in}\Sigma\coprod \partial_{\rm out} \Sigma$ into \emph{incoming} (source) and \emph{outgoing} (target) components.\footnote{For brevity we have omitted several important technical details in the definition. For example, a well-defined composition in $\EBord_d$ requires that $\partial_{\rm in}\Sigma$ and $\partial_{\rm out} \Sigma$ be collared. We refer to~\cite{ST11} for details.} 
\end{ex}

\begin{ex}[Super Euclidean bordism categories]
Generalizing Example~\ref{ex:Eucrigid}, \emph{super Euclidean} rigid geometries determine bordism categories $\EBord_{d|\delta}(\X)$ depending on a pair of natural numbers~$d|\delta$ (i.e., a  superdimension) as well as additional data~\cite[\S4.2]{ST11}. 
Roughly, a morphism in $\EBord_{d|\delta}(\X)$ is again a compact, spin $d$-manifold $\Sigma$ with flat metric and geodesic boundary together with a map $\Phi\colon \Sigma\to \X$. The new feature is that the $\delta$-dimensional complex spinor bundle $\bS\to \Sigma$ is regarded as ``odd," and \emph{supersymmetries} are automorphisms of $\bS\to \Sigma$ that mix even directions in $\Sigma$ with odd directions in $\bS$. To first approximation, $\EBord_{d|\delta}(\X)$ is just $\EBord_d(\X)$ with the addition of these odd automorphisms. Making this precise requires the language of supermanifolds~\cite{Freed5,DM}. 
\end{ex}

\vspace{.05in}

This leads to our first notion of supersymmetric field theory (elaborations will follow). 

\vspace{.05in}

\begin{defn}[{Sketch of \cite[Definition~5.2]{ST11}}] \label{defn:twistedQFT}
Let $\TA$ denote the 2-category of $\Z/2$-graded topological algebras, bimodules and bimodule maps with symmetric monoidal structure from the $\Z/2$-graded tensor product. A \emph{twisted $d|\delta$-dimensional supersymmetric Euclidean field theory over a smooth stack $\X$} is a natural transformation $E$
\beq
&&\begin{tikzpicture}[baseline=(basepoint)];
\node (A) at (0,0) {$\EBord_{d|\delta}(\X)$};
\node (B) at (5,0) {$\TA$};
\node (C) at (2.5,0) {$E \Downarrow$};
\draw[->,bend left=15] (A) to node [above] {$\one$} (B);
\draw[->,bend right=15] (A) to node [below] {$\alpha$} (B);
\path (0,0) coordinate (basepoint);
\end{tikzpicture}\label{eq:twistedEFT}
\eeq
between the constant functor $\one$ (to the monoidal unit of~$\TA$) and a functor $\alpha$ called the \emph{twist}. Twisted field theories over $\X$ are the objects of a groupoid. 
\end{defn}

We comment on some structural aspects of this definition. First, field theories are natural in $\X$: for a map $\X\to \X'$ the evident functor $\EBord(\X)\to \EBord(\X')$ induces a pullback from twisted field theories over $\X'$ to twisted field theories over~$\X$. Second, twists and twisted field theories have a multiplication: the monoidal structure on $\TA$ provides
$$
(E\colon \one \Rightarrow \alpha,E'\colon \one \Rightarrow \alpha')\xmapsto{\otimes} (E\otimes E'\colon \one\simeq \one\otimes \one\Rightarrow \alpha\otimes \alpha').
$$
A twist $\alpha$ is \emph{invertible} if there exists a twist $\alpha'$ and an equivalence $\alpha\otimes \alpha'\simeq \one$ with the trivial twist. All of the twists considered below will be invertible. 

Invertible twists are one way to capture \emph{anomalous} quantum field theories~\cite{Freedanomaly}.  Loosely speaking, projective symmetries determine anomalies, and twisted field theories are projective representations of bordism categories. 
Variants of Definition~\ref{defn:twistedQFT} have appeared (often independently) in the literature, e.g., the \emph{relative field theories} of Freed and Teleman~\cite{FreedTelemanRelative}. Stolz and Teichner's chosen terminology comes from a developing connection between anomalous quantum field theories and twisted cohomology classes~\cite{Stoffeltwists}. 

\begin{rmk}
We caution that the twists in Definition~\ref{defn:twistedQFT} are not related to topological twists from~\cite{Wittentoptwist}, e.g., the $A$- and $B$-twists from mirror symmetry~\cite{WittenTQFT}. 
\end{rmk}

\begin{ex}[Untwisted theories]
Let $\Vect$ be the symmetric monoidal category of $\Z/2$-graded (topological) vector spaces over $\C$ equipped with the tensor product. An \emph{untwisted field theory} is a symmetric monoidal functor $\EBord_{d|\delta}(\X)\to \Vect$, which determines a natural transformation~\eqref{eq:twistedEFT} for the trivial twist $\alpha=\one$~\cite[Lemma 5.7]{ST11}. 
\end{ex}

Part of the data of a twist is an algebra $\alpha(Y)$ for each object $Y\in \EBord_{d|\delta}(\X)$, and part of the data of a twisted field theory is a left $\alpha(Y)$-module, $E(Y)$. In this way, twisted field theories can be thought of as untwisted theories valued in $\alpha(Y)$-modules rather than vector spaces; see~\cite[page~52]{ST11}. 

\begin{ex}[Degree twists]
Specializing to $d|\delta=1|1$ and $2|1$ in Definition~\ref{defn:twistedQFT}, for each $n\in \Z$ there is a \emph{degree twist} \cite[\S5.3]{ST11} (see also Definitions~\ref{defn:11degreetw} and~\ref{defn:21degreetw} below)
\beq\label{eq:degtwist}
\EBord_{d|1}(\pt)\xrightarrow {\alpha_n}\TA,\qquad \alpha_n\otimes \alpha_m\simeq \alpha_{n+m} \quad n,m\in \Z, \ d\in \{1,2\}.
\eeq
For $d=1$, the $n$th degree twist is determined by the $n$th Clifford algebra $\cCl_n$, see~\eqref{eq:degreetwist11}. For $d=2$ it is determined by the \emph{$n$ (chiral) free fermions}, which includes the data of the infinite dimensional Clifford algebra $\cCl(L\R^n)$ where $L\R^n$ is the free loop space of $\R^n$, see~\eqref{eq:degreetwist21}. 
\end{ex}

To accurately capture the desired examples, we enhance twisted field theories with additional structure and property called \emph{reflection positivity}. In brief, this is $\Z/2$-equivariance data for twists and twisted field theories relative to the involution on bordisms that reverses orientations and the involution on $\TA$ from complex conjugation on~$\C$. Relative to this equivariance data, hermitian pairings in the theory are required to be inner products. For overviews, see~\cite[\S3]{FreedHopkins} or~\cite[page~18]{KontsevichSegal}. The full definition of reflection positivity in the context of Definition~\ref{defn:twistedQFT} is under development; a working definition appears in~\cite[\S{C.3}]{DBEtorsion} and we sketch the relevant features in Definition~\ref{defn:RP1} below. 

\begin{defn}\label{def:QFT}
Let $\QFT_{d|\delta}^\alpha(\X)$ denote the groupoid whose objects are reflection positive $\alpha$-twisted $d|\delta$-dimensional Euclidean field theories over $\X$, and whose morphisms are isomorphisms of twisted field theories compatible with the reflection structure. 
\end{defn}

\begin{defn}
For $n\in \Z$ and $\X$ a smooth stack, let $\QFT_{d|1}^n(\X):=\QFT_{d|1}^{\alpha_n}(\X)$ denote the category of field theories twisted by the pullback of~\eqref{eq:degtwist} along the map~$\X\to \pt$. 
\end{defn}

Specializing to a smooth manifold $\X=X$, the categories $\QFT_{d|1}^n(X)$ are the best-known approximation to the desired objects in Conjecture~\ref{conj}. But for there to be any hope of proving the conjecture, a further modification is required: excision demands fully-extended field theories, see~\cite[Conjecture 1.17]{ST11} and~\cite[\S1]{ST04}. This process of \emph{extending down} necessitates a 2-category of super Euclidean bordisms and a 3-category that deloops~$\TA$. Neither the source nor target category are uniquely determined by Definition~\ref{def:QFT}. Instead, there are a multitude of generalizations and (at present) it unclear which one to take.

\subsection{Enhancing the conjecture with Lurie's higher-categorical equivariance}\label{sec:higherequiv}

One way to narrow the zoo of possible definitions is to demand more structure in Conjecture~\ref{conj}. This makes it harder to guess the ``wrong" definition of field theory. A promising option is Lurie's higher-categorical equivariant elliptic cohomology, called \emph{2-equivariant elliptic cohomology} in~\cite[\S5]{Lurie_Elliptic}. This associates to each compact Lie group $G$ a derived scheme~$\mathcal{M}_G$, and to each level $[\ell]\in \H^4(BG;\Z)$ a line bundle $\mathcal{L}_\ell$ over~$\mathcal{M}_G$. The ordinary scheme underlying $\mathcal{M}_G$ is the moduli of flat $G$-bundles on elliptic curves and the ordinary line bundle underlying $\mathcal{L}_\ell$ has as sections nonabelian $\theta$-functions. Just as global sections of the sheaf of functions on the moduli stack of derived elliptic curves yields $\TMF$, global sections of $\mathcal{L}_\ell$ over $\mathcal{M}_G$ determines 2-equivariant $\TMF$, denoted~$\TMF_G^\ell$~\cite{LurieIII,LenartDavid}.

A compact Lie group $G$ acting on a manifold~$X$ determines the stack $\X=X\nsq G$,\footnote{The stack $X\nsq G$ classifies $G$-bundles with connection $(P,\nabla)$ and a $G$-equivariant map $P\to X$. Discarding the connection datum gives a forgetful functor $X\nsq G\to X\sq G$ to the usual quotient stack.} leading to $G$-equivariant refinements of Conjecture~\ref{conj}, compare \cite[\S1.7]{ST11}.

\begin{conj}\label{conjG}
For a compact Lie group $G$ acting on a manifold $X$ and a level $[\ell]\in \H^4(BG;\Z)$, there is a cocycle map generalizing~\eqref{eq:conj}
\beq
&&\begin{tikzpicture}[baseline=(basepoint)];
\node (A) at (0,0) {$\QFT^\ell_{2|1}(X\nsq G)$};
\node (B) at (5,0) {$\TMF_G^\ell(X)$};
\draw[->,dashed] (A) to node [above] {cocycle} (B);
\path (0,0) coordinate (basepoint);
\end{tikzpicture}\label{eq:bijection2}
\eeq
valued in $\ell$-twisted $G$-equivariant $\TMF$ of~$X$. Furthermore, for families of $G$-equivariant string manifolds there is an analog of the commuting triangle~\eqref{eq:conj}. 
\end{conj}

In physical language, the category $\QFT^\ell_{2|1}(\pt\nsq G)$ contains field theories with \emph{background gauge fields} and a specified \emph{'t Hooft anomaly}. By general theory~\cite{Freed_Det}, such anomalies determine line bundles over the moduli of $G$-bundles on genus~1 surfaces (i.e., elliptic curves). These line bundles are central to Grojnowski's twisted complex analytic equivariant elliptic cohomology \cite{Grojnowski,GKV}. Initial investigations suggest that twisted Euclidean field theories over $\pt\nsq G$ determine sections of such line bundles~\cite{BET0,BET1}. Hence when comparing to Lurie's framework, the expectation is that $\QFT^\one(-\nsq G)$ encodes the derived scheme $\mathcal{M}_G$ and $\QFT^\ell(-\nsq G)$ encodes sections of the $\theta$-line, $\Gamma(\mathcal{M}_G,\mathcal{L}_\ell)$. 

In addition to restricting the possible definitions of fully-extended field theory, Conjecture~\ref{conjG} illuminates a path towards constructing a comparison map between field theories and $\TMF$:  Lurie's sketches a universal property of 2-equivariant elliptic cohomology~\cite[\S5.5]{Lurie_Elliptic}. Conjecture~\ref{conjG} proposes a 2-equivariant structure for field theories, and so Lurie's universal property would lead to the desired comparison map. In this way, higher equivariance can be seen as an essential ingredient to unlocking Conjecture~\ref{conj}.

\begin{rmk}
The link between equivariant elliptic cohomology and 2-dimensional gauge theory is not new. The earliest ideas are perhaps due to Grojnowski~\cite{Grojnowski} and Segal \cite[Theorem 5.3]{Segal_Elliptic} \cite[Definition 2.4]{SegalElliptic}, though others have independently discovered similar parallels. For example, Witten's rigidity theorem for the Dirac operator on loop space can be understood in terms of $S^1$-equivariant elliptic cohomology \cite{BottTaubes,RosuEquivariant,AndoBasterra,RosuDelocalized}, and for finite groups it has long been clear that $G$-equivariant elliptic cohomology is related to \emph{twisted sectors} in string theory \cite{DHVW,Segal_Elliptic,DevotoI,DevotoII,HKR,MoravaHKR}. 
Loop groups and their positive energy representations also appear in both 2-dimensional field theory (e.g., WZW models) and equivariant elliptic cohomology \cite{Liu2,Kitchloo,FHT3,KitchlooII,GanterEllipticWCF,Kiran}. Equivariance also leads to natural elliptic power operations with evident connections to the geometry of string theory \cite{DMVV,Ando,Ganterstringy,GanterOrb,GanterHecke,Huan,ZhenMatt,Powerops}. 
 \end{rmk}

\subsection{Disclaimer and overview}
The chapter has two primary goals: (i) provide motivation and intuition for Conjecture~\ref{conj}, and (ii) describe the central challenges that remain in its verification. Both of these goals require an understanding of examples from theoretical physics, e.g., the sigma models and free fermions. Attempting to fit these (partially-constructed) examples within a mathematical framework for quantum field theory unfortunately introduces many technicalities that tend to obscure the larger picture. For this reason, each section below starts with somewhat ad hoc definitions of quantum mechanical systems and quantum field theories that explicitly capture analytical data found in the physics literature. For supersymmetric quantum mechanics in particular, the relevant data of a theory is totally unambiguous and the examples are well-understood. 
After connecting these ad hoc definitions with objects in algebraic topology, we contextualize the analytic data in the paradigm~\eqref{eq:twistedEFT} as part of a representation of a bordism category. This bordism definition points towards new structures in 2-dimensional quantum field theory that are difficult to glean from naive generalizations of supersymmetric quantum mechanics alone. To emphasize the ongoing development of this mathematical framework, we label a statement with sketched proof as a \emph{Hypothesis}; the references given indicate that a closely related statement has been proved, though with technical caveats we do not address here. We expect that any reasonable definition ought to make the stated Hypothesis true, and so this presentation is also intended to provide specific signposts to verify that proposed (fully-extended) definitions are on the right track.

We provide a brief overview of the chapter. Section~\ref{sec:QM} introduces quantum mechanics first from the analytic perspective before translating the basic structure into the bordism framework. Section~\ref{sec:SQM} repeats this for supersymmetric quantum mechanics, and sketches a version of the K-theory variant of Conjecture~\ref{conj}. Section~\ref{sec:4} starts with the analytic perspective on quantum field theory as a generalization of quantum mechanics, before translating these structure into the bordism framework. This framework leads to a variety of insights into the ultimate content of Conjecture~\ref{conj}. 

Throughout, $\cH$ will denote a (fixed) infinite-dimensional separable, complex Hilbert space. All tensor products of $\Z/2$-graded objects are the $\Z/2$-graded tensor product whose symmetry isomorphism encodes the Koszul sign.


\subsection{Acknowledgements}
I thank Emily Cliff, Andr\'e Henriques, Theo Johnson-Freyd, Matthias Ludewig, Laura Murray, Nat Stapleton, and Arnav Tripathy for many engaging conversations that shaped these ideas. I also thank Theo Johnson-Freyd and Laura Murray for feedback on a previous draft. The influence of Stephan Stolz and Peter Teichner would be difficult to overstate---I thank them for their many insights and generosity of ideas.

\section{The homotopy type of quantum mechanics}\label{sec:QM}

Let $M$ be a closed, oriented Riemannian manifold. \emph{Classical mechanics on $M$} studies paths $x\colon I\to M$ with the action functional $S(x)=\frac{1}{2}\int_I \| \dot x \| dt$. Critical points of $S$ are parameterized geodesics in~$M$. Using the exponential map and metric, such geodesics are identified with the cotangent bundle of $M$. Quantization of the symplectic manifold $T^*M$ produces \emph{quantum mechanics on $M$}, e.g., see~\cite{Weinsteinquantize}.
Below we shall only require two basic outputs of quantization: (i) the \emph{Hilbert space of states}~$\cH=L^2(M)$, and (ii) the \emph{Hamiltonian}~$H=\Delta$ given by the Laplacian on $M$. In \emph{Wick-rotated} quantum mechanics,~$(\cH,H)$ determine a representation of the semigroup~$(\R_{\ge 0},+)$,
\beq
&&\R_{\ge 0}\to \End(\cH),\qquad t\mapsto e^{-tH}\label{eq:semigroup}
\eeq
valued in self-adjoint operators on $\cH$ that are trace class for $t>0$ (using compactness of~$M$). This motivates the following, which will be made precise in~\S\ref{sec:QMdef} below.

\begin{defn}\label{defn:QMsystem}
A \emph{quantum mechanical system} is a state space $\cH$ and Hamiltonian~$H$ determining a self-adjoint, trace-class representation~\eqref{eq:semigroup}. Quantum mechanical systems~$(\cH,H)$ and~$(\cH',H')$ are \emph{equivalent} if there exists a unitary map $u\colon \cH\to \cH'$ with $H'=uHu^{-1}$. 
\end{defn}

The goal of this section is assemble quantum mechanical systems and their equivalences into a topological groupoid, i.e., a stack. Basic features of this stack generalize to the more sophisticated physical theories relevant to Conjecture~\ref{conj}. Indeed, much of the intuition for Conjecture~\ref{conj} originates in elaborations on structures in quantum mechanics. 



\begin{rmk}
There are variations to Definition~\ref{defn:QMsystem}. Dropping the trace class condition leads to \emph{noncompact} theories; examples include quantum mechanics on noncompact manifolds. Dropping the self-adjoint condition leads to Wick rotations of non-unitary theories. Although the complete picture is not fully-understood, the examples suggest that (Wick rotations of) compact, unitary quantum theories are the ones relevant to Conjecture~\ref{conj}. 
\end{rmk}


\subsection{The space of quantum mechanical theories is a $K(\Z,3)$}\label{sec:QMdef}

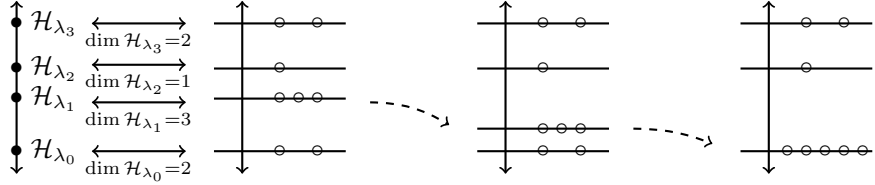
\begin{figure}
\beq\nonumber
\begin{tikzpicture}[baseline=(basepoint)];
\node (A) at (1,1.7) {$\bullet$};
\node (AA) at (1.5,1.7) {$\cH_{\lambda_3}$};
\node (B) at (1,1.1) {$\bullet$};
\node (BB) at (1.5,1.1) {$\cH_{\lambda_2}$};
\node (C) at (1,0) {$\bullet$};
\node (CC) at (1.5,0) {$\cH_{\lambda_0}$};
\node (D) at (1,.7) {$\bullet$};
\node (DD) at (1.5,.7) {$\cH_{\lambda_1}$};
\draw[thick,<->] (1,-.3) to (1,2);
\draw[<->,thick] (2,1.7) to  node [below] {$_{\dim\cH_{\lambda_3}=2}$} (3.25,1.7);
\draw[<->,thick] (2,1.15) to  node [below] {$_{\dim\cH_{\lambda_2}=1}$} (3.25,1.15);
\draw[<->,thick] (2,0) to node [below] {$_{\dim\cH_{\lambda_0}=2}$} (3.25,0);
\draw[<->,thick] (2.,.65) to  node [below] {$_{\dim\cH_{\lambda_1}=3}$} (3.25,.65);
\draw[thick,<->] (4,-.3) to (4,2);
\node (E) at (3.5,0) {};
\node (1) at (4.5,0) {$\circ$};
\node (2) at (5,0) {$\circ$};
\node (EE) at (5.5,0) {};
\node (F) at (3.5,1.7) {};
\node (1) at (4.5,1.7) {$\circ$};
\node (1) at (5,1.7) {$\circ$};
\node (FF) at (5.5,1.7) {};
\node (G) at (3.5,1.1) {};
\node (1) at (4.5,1.1) {$\circ$};
\node (GG) at (5.5,1.1) {};
\node (H) at (3.5,.7) {};
\node (1) at (4.5,.7) {$\circ$};
\node (1) at (4.75,.7) {$\circ$};
\node (1) at (5,.7) {$\circ$};
\node (HH) at (5.5,.7) {};
\draw[-,thick] (E) to (EE);
\draw[-,thick] (F) to (FF);
\draw[-,thick] (G) to (GG);
\draw[-,thick] (H) to (HH);

\draw[->,thick,dashed,bend left=15] (5.7,.65) to (6.75,.35);
\draw[thick,<->] (7.5,-.3) to  (7.5,2);
\node (J) at (7,0) {};
\node (1) at (8,0) {$\circ$};
\node (2) at (8.5,0) {$\circ$};
\node (JJ) at (9,0) {};
\node (K) at (7,1.7) {};
\node (1) at (8,1.7) {$\circ$};
\node (1) at (8.5,1.7) {$\circ$};
\node (KK) at (9,1.7) {};
\node (L) at (7,1.1) {};
\node (1) at (8,1.1) {$\circ$};
\node (LL) at (9,1.1) {};
\node (M) at (7,.3) {};
\node (1) at (8,.3) {$\circ$};
\node (1) at (8.25,.3) {$\circ$};
\node (1) at (8.5,.3) {$\circ$};
\node (MM) at (9,.3) {};
\draw[-,thick] (J) to (JJ);
\draw[-,thick] (K) to (KK);
\draw[-,thick] (L) to (LL);
\draw[-,thick] (M) to (MM);

\draw[->,thick,dashed,bend left=10] (9.2,.3) to (10.25,.1);
\draw[thick,<->] (11,-.3) to  (11,2);
\node (J) at (10.5,0) {};
\node (1) at (11.5,0) {$\circ$};
\node (2) at (12,0) {$\circ$};
\node (JJ) at (12.5,0) {};
\node (K) at (10.5,1.7) {};
\node (1) at (11.5,1.7) {$\circ$};
\node (1) at (12,1.7) {$\circ$};
\node (KK) at (12.5,1.7) {};
\node (L) at (10.5,1.1) {};
\node (1) at (11.5,1.1) {$\circ$};
\node (LL) at (12.5,1.1) {};
\node (1) at (11.25,0) {$\circ$};
\node (1) at (11.75,0) {$\circ$};
\node (1) at (12.25,0) {$\circ$};
\draw[-,thick] (J) to (JJ);
\draw[-,thick] (K) to (KK);
\draw[-,thick] (L) to (LL);
\path (0,.5) coordinate (basepoint);
\end{tikzpicture}\qquad 
\eeq

\caption{On the left is a Hamiltonian in the configuration space model, where each label $\cH_\lambda\subset \cH$ for $\lambda\in \R$ is an eigenspace of~$H$ with eigenvalue~$\lambda$. When labels collide, the subspaces add. Energies can also ``escape" off to infinity. The pictures on the right follow notation in~\cite{Wittenbreaking}: the dimension of each eigenspace is recorded by $\circ$'s (the number of \emph{states}). These dimensions determine a configuration up to isomorphism. As indicated by the dashed arrows, colliding labels adds dimensions.}
\label{fig1}
\end{figure}

Hamiltonians in Definition~\ref{defn:QMsystem} are a class of self-adjoint unbounded operators\footnote{An \emph{unbounded operator on $\cH$} is a domain $\dom(H)\subset \cH$ and a linear map $H\colon \dom(H)\to \cH$. An unbounded operator is \emph{self-adjoint} if it is self-adjoint on the closure of its domain. } 
that admit a geometric description in terms of configuration spaces \cite[\S4]{HST}. In brief, given a self-adjoint unbounded operator~$H$, the collection of mutually orthogonal eigenspaces $\{\cH_\lambda\subset \cH\}_{\lambda\in \R}$ determines a configuration on $\R$ where $\lambda\in \R$ is labeled by the subspace $\cH_\lambda\subset \cH$; see Figure~\ref{fig1}. Conversely, given a configuration on~$\R$, define an unbounded operator with domain $\bigoplus_{\lambda\in \R} \cH_\lambda$ where $H\colon \bigoplus_{\lambda\in \R} \cH_\lambda\to \cH$ is the linear map that scalar multiplies by~$\lambda$ on~$\cH_\lambda$.

The topology on unbounded operators determines a topology on configurations: a pair of configurations is close if there is an (energy cutoff) $\Lambda\in \R$ such that the eigenspaces with eigenvalues $\lambda<\Lambda$ are close; we refer to \cite[Definition 4.3, Proposition 4.4]{HST} for the precise definition of this topology. 
An unbounded operator $H$ gives a semigroup representation~\eqref{eq:semigroup} on the closure of the domain of $H$, which we extend by the zero map on the orthogonal complement of the domain~\cite[Proposition 5.9]{HST}. For a semigroup representation to be trace class, we must demand additional properties on~$H$: the nonzero labels $\lambda\in \R$ are bounded below, each subspace $\cH_\lambda\subset \cH$ has finite dimension, the spectrum ${\rm Spec}(H)\subset \R$ is discrete, and the $\lambda$ satisfy certain growth conditions.
We refer to \cite[\S3]{HST} for details.

\begin{defn} \label{defn:Ham}
The \emph{space of Hamiltonians}, denoted $\Ham$, is the space of configurations $\{\cH_\lambda\subset \cH\}_{\lambda\in \R}$ determining unbounded operators that generate trace class self-adjoint semigroup representations~\eqref{eq:semigroup}.
\end{defn}

The unitary group $\U(\cH)$ acts on unbounded operators by conjugation, and this affords a continuous action on $\Ham$ by $\U(\cH)$ equipped with the compact-open topology.\footnote{Continuity follows from the description of open subsets in $\Ham$~\cite[Definition 4.3]{HST} and continuity of the action of $\U(\cH)$ on compact operators, see \cite[Proposition A1.1]{AtiyahSegaltwistedK} and \cite{EspinozaUribe}.} The conjugation action factors through the quotient $\PU(\cH)=\U(\cH)/\U(1)$ by scalars.  


\begin{defn}\label{eq:QM}
Define the topological category $\qm:=\Ham\sq \PU(\cH)$ as the quotient groupoid for the action induced by the conjugation action of the unitary group. 
\end{defn}

We emphasize that we do not require the domain of an unbounded operator be dense in $\cH$. Hence, Definition~\ref{defn:Ham} does not require the inclusion $\bigoplus_{\lambda\in \R} \cH_\lambda \subset \cH$ to be dense. 

\begin{ex}[Topological field theories]\label{ex:TQFT}
Any subspace $V\subset \cH$ defines an unbounded operator $T_V$ with $\dom(T_V)=V$ and linear map the zero map $T_V\equiv 0\colon V\to \cH$. This unbounded operator determines a point of $\Ham$ if and only if $V$ is finite-dimensional. The associated quantum system is \emph{topological}, with configuration $\mathcal{H}_0=V$ supported at $0\in \R$. Note that the time evolution~\eqref{eq:semigroup} is constant in~$t$: there are no dynamics.
Hypothesis~\ref{hyp:QMandBord} compares this with Atiyah's definition of topological field theory. 
\end{ex}


\begin{rmk}
The objects in the groupoid $\qm$ and their depictions in Figure~\ref{fig1} have the following interpretation in quantum mechanics. For an object in $\qm$, the eigenspace~$\cH_\lambda$ corresponds to states of energy~$\lambda$. The topology on objects of $\qm$ allows for perturbations of the energy values, as well as splitting and merging of energy eigenspaces. The eigenspace with the smallest eigenvalue contains the \emph{ground states} of the theory. 
When the domain $\dom(H)\subset \cH$ of $H$ is not dense, the complement of the closure $\overline{\dom(H)}{}^\perp$ is the eigenspace ``at infinity," i.e., the infinite-energy states. From this perspective, it is perhaps better to view the true state space of the theory as the closure of the finite-energy eigenstates $\bigoplus_{\lambda<\infty} \cH_\lambda\subset \cH$. With this point of view, the state space of a topological theory is entirely comprised of ground states. 
\end{rmk}

\begin{prop}[{\cite[\S6]{HST} \cite[\S3.1]{ST04}}]\label{prop:KZ3}
When $\cH$ is an infinite-dimensional separable Hilbert space, the geometric realization of $\qm$ is the Eilenberg--MacLane space
\beq\label{eq:realization}
|\qm|\simeq K(\Z,3).
\eeq

\end{prop}

\begin{proof}[Proof sketch.]
Contractibility of the space of Hamiltonians is \cite[Theorem 6.35]{HST}. In brief, $\Ham$ has a basepoint given by the unbounded operator with domain $\{0\}$. There is a contracting homotopy to this basepoint that deforms Hamiltonians by uniformly pushing all eigenspaces to infinity. We conclude $|\qm|\simeq |\pt\sq \PU(\cH)|\simeq B\PU(\H)$.

Next, we claim $\U(\cH)$ is contractible; in the norm topology this is a classical result of Kuiper~\cite{Kuiper} and in the compact-open topology this is~\cite[Proposition A2.1]{AtiyahSegaltwistedK}. Using the fiber sequence $\U(1)\to \U(\cH)\to \PU(\cH)$ to compute $\pi_*\PU(\cH)$, we conclude that~$B\PU(\cH)\simeq BK(\Z,2)\simeq  K(\Z,3)$. 
\ep

The classifying space~\eqref{eq:realization} captures the homotopical content of quantum mechanics, but more refined information comes from a careful analysis of the stack~$\qm$. Below we briefly describe connections with projective representations of compact Lie groups. 

%
%
%

\begin{defn}\label{defn:sym1}
Let $G$ be a compact Lie group. The stack of \emph{quantum mechanical systems with $G$-symmetry} is the mapping stack $\Map(\pt\sq G, \qm)$. 
\end{defn}

To unpack this definition, an object in $\Map(\pt\sq G, \qm)$ (i.e., the data of a map of stacks $\pt\sq G\to \qm$) is a point in $H\in \Ham$ together with a homomorphism $\alpha\colon G\to \PU(\cH)$ satisfying a compatibility property. To state this compatibility, the pullback
\beq\label{eq:pullbacksequence}
\begin{array}{ccccccccc} 1 & \to & \U(1) & \to &  G^\alpha & \to &  G& \to & 1
\\
 &  & =\downarrow &  &  \downarrow & \lrcorner &  \phantom{\alpha}\downarrow \alpha &  & \\
1 & \to & \U(1) & \to &  \U(\cH) & \to &  \PU(\cH)& \to & 1\end{array}
\eeq
determines a central extension $G^\alpha$ of $G$ and a unitary $G^\alpha$-representation on $\cH$. The compatibility property between $\alpha$ and $H$ is that the $G^\alpha$-action on $\cH$ commutes with~$H$. 

\begin{rmk}
A quantum system with $G$-symmetry enhances the configurations in Figure~\ref{fig1}: each $\cH_\lambda\subset \cH$ is a finite-dimensional $G^
\alpha$-representation, and colliding points in the configuration corresponds to the direct sum of representations. 
\end{rmk}

\begin{ex}[Quantum mechanics on a $G$-manifold]\label{ex:GQM}
Let $M$ be an oriented, compact Riemannian manifold with isometric and orientation preserving $G$-action. Then the example discussed at the beginning of the section inherits additional structure: $\cH=L^2(M)$ acquires a $G$-action that commutes with the Hamiltonian $H=\Delta$. This gives a map 
$$
\pt\sq G\to\Ham\sq \U(\cH)\twoheadrightarrow \Ham\sq \PU(\cH)=\qm
$$
defining a quantum system with (non-projective) $G$-symmetry.
\end{ex}

\begin{ex}[Free fermions]\label{ex:fermion}
In dimension~1, the classical theory of \emph{$n$ free fermions} has fields and action functional
\beq\label{Eq:FFaction}
\psi\colon \R \to \R^{n},\qquad 
S(\psi)=\frac{1}{2}\int_\R \langle \psi,\partial_t\psi \rangle dt
\eeq
where $\langle-,-\rangle$ is the standard inner product on $\R^n$. The classical solutions are the constant maps to~$\R^{n}$. Viewing $\R^n$ as a purely odd vector space endows these classical solutions with an odd symplectic structure from the standard inner product on $\R^n$. Quantization is equivalent to constructing an irreducible module~$V$ over $\cCl_n$ \cite[pg.~623]{strings1}, which defines the state space of the (quantum) free fermions. Because classical solutions are constant maps (i.e., there are no dynamics) the Hamiltonian is $H\equiv 0$, and so the quantum theory is topological. The orthogonal group $\O_n$ acts on the Clifford algebra $\cCl_n$ by algebra isomorphisms. Twisting the $\cCl_n$-action on $V$ by this algebra isomorphism produces an irreducible $\cCl_n$-module for each element of $\O_n$. By Schur's lemma, this determines a projective $\O_n$-action on $V$. When~$n$ is even, $\Lambda^\bullet \C^{n/2}$ is the unique (graded) irreducible $\cCl_n$-module, constructed as a Fock space~\cite[Definition 2.2.4]{ST04}. After choosing an isometric embedding $V\subset \cH$ (the space of which is contractible), the free fermions determine a map
\beq\label{eq:FFsymmetry}
\pt\sq \O_n \to \qm
\eeq
whose associated central extension~\eqref{eq:pullbacksequence} is $1\to \U(1)\to \Pin^c_n\to \O_n\to 1$. Hence, the free fermions determine a topological quantum theory with (projective) $\O_n$-symmetry. 
\end{ex}

\begin{rmk}
The above description of the free fermions fails to account for the (important) $\Z/2$-gradings on Clifford algebras and Clifford modules. We remedy this in Example~\ref{ex:spintheories}. 
\end{rmk}

\subsection{Quantum systems as representations of cobordism categories}\label{sec:1bord}
Specializing Definition~\ref{defn:twistedQFT}, let $\EBord_1$ denote the \emph{1-dimensional Euclidean bordism category}. Objects are finite disjoint unions of $\bullet=\pt^+$ and $\circ=\pt^-$, the point with its two choices of orientation. Morphisms are compact, oriented 1-dimensional Riemannian bordisms with boundary. One finds (e.g., by Morse theory) that such bordisms are disjoint unions of those in Figure~\ref{fig:EBord}. 

Bordisms compose by gluing along matching source and target data. The disjoint union of 0- and 1-manifolds endows the category $\EBord_1$ with a symmetric monoidal structure. 
Composition of bordisms introduces relations, e.g.,  for $s,t\in \R_{\ge 0}$, $\sI_s^\pm\circ \sI_t^\pm \simeq \sI_{s+t}^\pm$
corresponds to the fact that lengths of intervals add, giving a geometric counterpart to the semigroup~$\R_{\ge 0}$ from~\eqref{eq:semigroup}. Other relations can similarly be deduced from the pictures. A somewhat different relation uses the symmetry isomorphism $\sigma \colon \pt^+\coprod \pt^-\xrightarrow{\sim} \pt^-\coprod \pt^+$, giving $\sigma \circ \sL^\pm_t\simeq \sL^\mp_t$ and $ \sR^\pm_t\circ \sigma \simeq \sR^\mp_t$. A complete set of generators and relations for the closely related unoriented Riemannian bordism category are given in~\cite[Theorem 6.11]{HST}, see also \cite[Proof of Proposition 3.1.1]{ST04}. 

A quantum mechanical system in the sense of Definition~\ref{defn:QMsystem} determines an (untwisted) field theory with values
\beq\label{eq:repof1EB0}
&&E\colon \EBord_1 \to \Vect,\qquad  \begin{array}{c} E(\pt^+)=V, \ E(\pt^-)=\overline{V}, \ E(\sI^\pm_t)=e^{-tH}, \\ 
E(\sL^+_t)=\langle-,e^{-tH}-\rangle,\  E(S_t)=\Tr(e^{-tH})
\\
E(\sR^+)\colon \C\to \overline{V}\otimes V\hookrightarrow \End(V), \ E(\sR^+)(1)=e^{-tH} \end{array}
\eeq
where $V:={\rm dom}(H)\subset \cH$ is the domain of $H$. 
Conversely, the values $E(\sI^+_t)$ of a symmetric monoidal functor $E$ determine a semigroup representation~\eqref{eq:semigroup}. To guarantee that $E(\sI^+_t)$ is valued in self-adjoint operators (and hence determines a quantum mechanical system)
one must impose reflection positivity. To explain this data and structure in the present setting, take the $\Z/2$-action on $\EBord_1$ generated by orientation reversal of 1-manifolds. On the generators in Figure~\ref{fig:EBord}, orientation reversal is
\beq
&&\pt^\pm\mapsto \pt^\mp,\qquad \sI_t^\pm\mapsto \sI_t^\mp,\qquad \sR_t^\pm\mapsto \sR_t^\mp,\qquad \sL_t^\pm\mapsto \sL_t^\mp. \label{eq:orientationreverse}
\eeq
Physically, this action reverses the arrow of time. We also have a $\Z/2$-action on vector spaces generated by complex conjugation. 

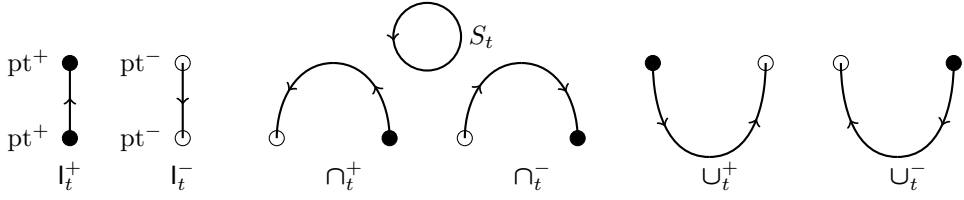
\begin{figure}
\beq
\begin{tikzpicture}[baseline=(basepoint)];
\draw [fill]  (0,1.5) circle [radius=0.1] node [black,left=4] {$\pt^+$};
\draw  [fill] (0,.5,0) circle [radius=0.1] node [black,left=4] {$\pt^+$};
\draw[decoration={markings,
        mark=at position \halfway with \arrow{<}},
        postaction=decorate,thick] (0,1.5) to (0,.5);
\node (Iplus) at (0,0) {$\sI_t^+$};

\draw (1.5,1.5) circle [radius=0.1] node [black,left=4] {$\pt^-$};
\draw (1.5,0.5) circle [radius=0.1] node [black,left=4] {$\pt^-$};
\draw[decoration={markings,
        mark=at position \halfway with \arrow{>}},
        postaction=decorate,thick] (1.5,1.5) to  (1.5,0.5);
\node (Iminus) at (1.5,0) {$\sI_t^-$};

\draw  (2.75,.5) circle [radius=0.1] node [black,left=4] {};
\draw[fill]  (4.25,.5) circle [radius=0.1] node [black,right=4] {};
\draw[decoration={markings,
        mark=at position \halfway with \arrow{<}},
        postaction=decorate,thick] (2.75,.5) to [out=90,in=180] (3.5,1.5);
        \draw[decoration={markings,
        mark=at position \halfway with \arrow{<}},
        postaction=decorate,thick] (3.5,1.5) to [out=0,in=90] (4.25,0.5);

\node (Rminus) at (3.65,0) {$\sL_t^+$};

\draw  (5.25,.5) circle [radius=0.1] {};
\draw [fill] (6.75,.5) circle [radius=0.1] {};
\draw[decoration={markings,
        mark=at position \halfway with \arrow{>}},
        postaction=decorate,thick] (5.25,.5) to [out=90,in=180] (6,1.5);
        \draw[decoration={markings,
        mark=at position \halfway with \arrow{>}},
        postaction=decorate,thick] (6,1.5) to [out=0,in=90] (6.75,0.5);

\node (Rplus) at (6.15,0) {$\sL_t^-$};

\draw [fill] (7.75,1.5) circle [radius=0.1] node [black,left=4] {};
\draw  (9.25,1.5) circle [radius=0.1] node [black,right=4] {};
\draw[decoration={markings,
        mark=at position \halfway with \arrow{>}},
        postaction=decorate,thick]  (7.75,1.5) to [in=180,out=270] (8.5,.25);
\draw[decoration={markings,
        mark=at position \halfway with \arrow{>}},
        postaction=decorate,thick]  (8.5,.25) to [in=270,out=0] (9.25,1.5);

\node (Lplus) at (8.65,0) {$\sR_t^+$};

\draw  (10.25,1.5) circle [radius=0.1]{};
\draw [fill](11.75,1.5) circle [radius=0.1] node [black,right=4] {};
\draw[decoration={markings,
        mark=at position \halfway with \arrow{<}},
        postaction=decorate,thick]  (10.25,1.5) to [in=180,out=270] (11,.25);
\draw[decoration={markings,
        mark=at position \halfway with \arrow{<}},
        postaction=decorate,thick]  (11,.25) to [in=270,out=0] (11.75,1.5);

\node (Lminus) at (11.15,0) {$\sR_t^-$};

\draw[decoration={markings,
        mark=at position \halfway with \arrow{>}},
        postaction=decorate,thick] (4.75,1.85) circle [radius=.45] node [right=12] {$S_t$};


\path (0,.75) coordinate (basepoint);
\end{tikzpicture}\nonumber
\eeq
\caption{Morphisms in $\EBord_1$ are finite disjoint unions of the above generators. The 1-manifold $\sS_t$ is the circle of circumference~$t$, viewed as a bordism from the empty set to itself. The 1-manifolds underlying the bordisms $\sI_t^\pm$, $\sR_t^\pm$ and~$\sL_t^\pm$ are all intervals of length~$t$; the distinctions between them come from their orientation (indicated by $+$ or $-$) together with source and target data. We follow the convention that reads bordisms in the upwards direction, e.g., the source of $\sL_t^+$ is $\pt^-\coprod \pt^+$. }
\label{fig:EBord}
\end{figure}
\begin{defn}[Reflection positivity, sketch] \label{defn:RP1}
A \emph{reflection structure} for the functor~\eqref{eq:repof1EB0} is equivariance data with respect to $\Z/2$-actions generated by orientation reversal and complex conjugation, which includes the data of an involutive isomorphism
\beq\label{eq:reflectionstructure}
E(\pt^+)\stackrel{\sim}{\to} \overline{E(\pt^-)}.
\eeq
A representation~\eqref{eq:repof1EB0} with a reflection structure is \emph{reflection positive} if the hermitian pairing
\beq\label{eq:positivity}
\langle-,-\rangle\colon \overline{E(\pt^+)}\otimes E(\pt^+)\simeq E(\pt^-)\otimes E(\pt^+)\xrightarrow{E(\sL_0^-)} \C
\eeq
is positive. Above, $\sL_0^-$ is the $t\to 0$ limit of $\sL_t^-$. 
\end{defn}

Following Definition~\ref{def:QFT},~\eqref{eq:repof1EB0} and~\eqref{eq:reflectionstructure} are the data of an object in the groupoid $\QFT_1:=\QFT_1^\one(\pt)$ of untwisted 1-dimensional Euclidean field theories over~$\X=\pt$.

\begin{hyp}[{Compare \cite[Definition 3.2.10]{ST04} and \cite[Corollary 6.28]{HST}}] \label{hyp:QMandBord}
There is an epimorphism of stacks $\QFT_1\to \qm$ sending a 1-dimensional Euclidean field theory in the sense of Definition~\ref{def:QFT} to a quantum mechanical system in the sense of Definition~\ref{eq:QM}. 

\end{hyp}
\begin{proof}[Sketch.] 
On objects, the map of stacks extracts a semigroup representation from an object $E\in \QFT_1$ via $E(\sI^+_t)$ and a (choice of) isometric embedding~$E(\pt^+)\subset \cH$. An isomorphism between trivially twisted field theories determines an invertible $\C$-$\C$ bimodule (i.e., complex line)~$L$ and a unitary map $L\otimes E(\spt^+)\to E'(\spt^+)$ of pre-Hilbert spaces~\cite[\S5.1]{ST11}. 
Such data is determined by a projective isomorphism of pre-Hilbert spaces, which extends to an element $u\in \PU(\cH)$. This completes the sketch of a functor~$\QFT_1\to \qm$ 

To show this functor is an epimorphism, the basic idea
is to promote the pictures in Figure~\ref{fig:EBord} to a generators and relations presentation of $\EBord_1$ (compare~\cite[Theorem~6.23]{HST}). 
Such a presentation reveals that a 1-dimensional Euclidean field theory is determined by the inner product space $(E(\pt^+),\langle-,-\rangle)$ and the self-adjoint semigroup representation determined by~$E(\sI_t^+)$. By \cite[Lemma 5.10]{HST}, we have $E(\sI_t^+)=e^{-tH}$  from a uniquely determined unbounded operator $H$ with domain $E(\pt^+)\subset \cH$. The construction~\eqref{eq:repof1EB0} then determines a 1-sided inverse to $\QFT_1\to \qm$. 
\ep


\begin{rmk}\label{rmk:weirdsheafdef}
One might hope for the map $\QFT_1\to \qm$ to be an equivalence, but this is not possible because distinct topological inner product spaces (underlying non-isomorphic functors~\eqref{eq:repof1EB0}) can have the same Hilbert completion. We discuss this further in Remark~\ref{rmk:weirdsheaf2}. 

\end{rmk}

\begin{cor}
Topological theories in the sense of Example~\ref{ex:TQFT} correspond to 1-dimensional hermitian topological theories in the sense of Atiyah~\cite{AtiyahTQFT}.
\end{cor}
\bp
Under the functor in Hypothesis~\ref{hyp:QMandBord}, quantum mechanical systems determined by finite-dimensional subspace $V\subset \cH$ with the zero Hamiltonian correspond to 1-dimensional field theories determined by a finite-dimensional inner product space $V$ and $E(\sI^+)=\id$. This recovers a hermitian topological field theory~\cite[Eq.~5]{AtiyahTQFT}.
\ep

\subsection{Additional insights from the bordism framework}\label{sec:1gens}

The bordism framework of Definition~\ref{def:QFT} illuminates various possible generalizations of quantum mechanical systems. Although most of the following can be understood using enhancements of Definition~\ref{defn:QMsystem} alone, when passing to 2-dimensional quantum theories the desired (e.g., equivariant) refinements become more difficult to glean without the bordism language as a guide. 




 \begin{ex}[Orientation reversal and Real structures]\label{ex:OrReal1}
As previously indicated in~\eqref{eq:orientationreverse}, the bordism category $\EBord_1(\X)$ supports an involution from orientation reversal of bordisms. This leads to a $\Z/2$-action on the groupoid of twisted field theories over $\X$. The (homotopy) fixed points for this action are \emph{unoriented 1-dimensional Euclidean field theories}. Being $\Z/2$-fixed is additional data, including an involutive isomorphism 
\beq\label{eq:realstructure}
E(\pt^+)\stackrel{\sim}{\to} \overline{E(\pt^+)},\qquad\qquad E\in (\QFT_1)^{h\Z/2}
\eeq
i.e., a real structure on $E(\pt^+)$. Hence, unoriented theories correspond to the category of quantum systems with real state spaces, compare~\cite[Corollary 6.28]{HST}. In terms of Definition~\ref{defn:QMsystem}, an unoriented theory is the data of a real structure $r\colon \cH\simeq \overline{\cH}$ with the property that~$H$ is a real operator, $r\circ H=H\circ r$. 
 \end{ex}


For a smooth stack $\X$, Definition~\ref{def:QFT} provides generalizations of Hypothesis~\ref{hyp:QMandBord} via the category $\EBord_1(\X)$ of 1-dimensional Euclidean bordisms over $\X$.
\begin{ex}[Background gauge fields and $G$-representations]\label{ex:rep1}
Consider the bordism category $\EBord_1(\pt\nsq G)$ over the stack $\pt\nsq G$ of $G$-bundles with connection. The objects of $\EBord_1(\pt\nsq G)$ are $G$-bundles over finite disjoint unions of oriented points (the connection datum is trivial on a 0-manifold). Up to isomorphism, these bundles can be identified with the trivial bundle with its $G$-action by automorphisms. Hence for $E\in \QFT_1(\pt\nsq G)$, the value on the trivial $G$-bundle over $\pt^+$ carries the data of a $G$-action, i.e., is a $G$-representation. After fixing a trivialization at the boundaries, a $G$-bundle over an interval is equivalent to an element $g\in G$ gotten from parallel transport; let $P_g$ denote such a $G$-bundle with connection over the interval. Then the semigroup from composition of intervals is enhanced by
\beq
(\sI_s^\pm,P_g)\circ (\sI_t^\pm,P_h) \simeq (\sI_{s+t}^\pm,P_{gh}),\qquad s,t\in \R_{>0}, \ g,h\in G\label{eq:Gsemigroup}
\eeq
and therefore includes the data of multiplication in $G$. Compatibility with source and target means that the resulting $G$-action is the same as the action by automorphisms of the trivial $G$-bundle. We further observe that orientation reversal changes the direction of parallel transport, sending $(\sI_t^\pm,P_g)$ to $(\sI^\mp_t,P_{g^{-1}})$. Hence, a reflection-positive 1-dimensional Euclidean field theory over $\pt\nsq G$ is the data of (i) a \emph{unitary} $G$-representation on an inner product space, and (ii) a self-adjoint, trace class semigroup representation commuting with the $G$-representation.  Such semigroup representations recover quantum systems with (non-projective) $G$-symmetry in the sense of Definition~\ref{defn:sym1}. If the theory is topological (and so does not depend on the geometric parameter $t\in \R_{>0}$) a representation of the bordism category is determined by a finite-dimensional unitary $G$-representation. 
\end{ex}

\begin{ex}[Projective $G$-representations]\label{ex:prep1}
Continuing with the twisted generalization of the previous example, an invertible twist $\EBord_1(\pt\nsq G)\to \TA$ will assign the invertible complex algebra $\C$ to the point together with a $G$-action by invertible $\C$-modules; this is equivalent data to a central extension of $G$ by $\U(1)$, e.g., see \cite[Example 1.76]{Freedalg}. Following similar analysis to the previous example, we conclude that a twisted 1-dimensional Euclidean field theory over $\pt\nsq G$ is (i) a unitary $G^\alpha$-representation for a central extension of $G$ in the notation of~\eqref{eq:pullbacksequence}, and (ii) a semigroup representation commuting with the $G^\alpha$-action. This recovers a quantum system with $G$-symmetry in the sense of Definition~\ref{defn:sym1}. 



\end{ex}

%

\begin{ex}[Vector bundles with connection] \label{ex:VB1}
Next consider the category $\EBord_1(X)$ of bordisms over a smooth manifold $X$. A functor $E\colon \EBord_1(X)\to \Vect$ assigns to every point $x\in X$ a vector space $E(x)=V_x$, and for every metrized path $\sI_t^+\colon x\to y$ in $X$ a linear map $E(\sI_t^+)\colon V_x\to V_y$. Applying a similar construction to~\eqref{eq:repof1EB0}, given a metrized vector bundle with compatible connection parallel transport determines a functor~\cite{LudewigStoffel,BEPTFT}
\beq\label{eq:dumbparallel}
\Vect^\nabla(X)\to \QFT^\one_1(X).
\eeq
Since parallel transport does not depend on the length of a path, the image of~\eqref{eq:dumbparallel} are the topological field theories over $X$, see~\cite{LudewigStoffel,BEPTFT}. More generally, a functor $E\colon \EBord_1(X)\to \Vect$ is determined by a length-dependent version of parallel transport on a vector bundle $V\to X$. In particular, the constant paths to $x\in X$ give an $X$-family of semigroup representations, determining an $X$-family of Hamiltonians. After choosing an embedding $V\subset X\times \cH$ as a summand of the trivial Hilbert bundle (the space of which is contractible) and forgetting the transport data, we obtain a map 
\beq\label{eq:famQFT1}
\QFT^\one_1(X)\to \Map(X,\Ham\sq \U(\cH))
\eeq
from theories over $X$ to the version of $\qm$ with unitary (non-projective) symmetries. 
\end{ex}

\begin{ex}[Gerbe-twisted vector bundles with connection] \label{ex:pVB1}
Continuing with the twisted generalization of the previous example, the data of a twist is an algebra bundle over~$X$ for a cocycle valued in the Morita bicategory. By \cite[Examples~1.71 and 1.74]{Freedalg}, this is equivalent to a gerbe on~$X$, see also \cite{Stoffeltwists}. An object of $\QFT^\alpha_1(X)$ can therefore be thought of as a gerbe-twisted vector bundle equipped with a metric-dependent (projective) parallel transport. Forgetting the transport data we obtain a map analogous to~\eqref{eq:famQFT1}
$$
\QFT^\alpha_1(X)\to \Map(X,\qm)
$$
where the composite $X\to \qm\to \pt\sq \PU(\cH)$ classifies the gerbe determined by $\alpha$. 

\end{ex}

\begin{ex}[Twisted theories valued in $A$-modules]\label{ex:Amod}
A $*$-algebra $A$ determines a reflection positive twist with values 
\beq\label{eq:twist1EB0}
&&\alpha_A\colon \EBord_1 \to \TA,\qquad  \begin{array}{ll} \alpha_A(\pt^+)=A, & \alpha_A(\pt^-)=\overline{A}, \\ \alpha_A(\sI^+_t)= \ _AA_A, & \alpha_A(\sI^-_t)= \ _{\overline{A}} {\overline{A}}_{\overline{A}} \\ 
\alpha_A(\sL^+_t)=A_{A\otimes \overline{A}}, & \alpha_A(\sL^-_t)=\overline{A}_{\overline{A}\otimes A}  \\  \alpha_A(\sR^+_t)= \ _{A\otimes \overline{A}} A & \alpha_A(\sR^+_t)= \  _{\overline{A}\otimes A} \overline{A} \\ 
\alpha_A(S_t)=A\otimes_{A\otimes \overline{A}} \overline{A} 
  \end{array}
\eeq
where the $*$-structure provides an involutive isomorphism $A^\op\simeq \bar A$ that is used in defining the module structures. We observe that the twist determined by~\eqref{eq:twist1EB0} is \emph{topological}: its values are independent of the length parameters~$t\in \R_{>0}$. A $\alpha_A$-twisted field theory $E$ is determined by
\beq\label{eq:twisted1EB0}
&&\begin{tikzpicture}[baseline=(basepoint)];
\node (A) at (0,0) {$\EBord_1$};
\node (B) at (2.5,0) {$\TA,$};
\node (C) at (1.25,0) {$E \Downarrow$};
\draw[->,bend left=15] (A) to node [above] {$\one$} (B);
\draw[->,bend right=15] (A) to node [below] {$\alpha_n$} (B);
\path (0,0) coordinate (basepoint);
\end{tikzpicture}\qquad  \begin{array}{ll} E(\pt^+)=\ _AV, & E(\pt^-)=\ _{\overline{A}} \overline{V}, \\ E(\sI^+_t)= \ _AV \to {_AV}, & E(\sI^-_t)= \ _{\overline{A}}\overline{V}\to _{\overline{A}}\overline{V} \\ 
E(\sL^+_t)=V\otimes_{\overline{A}} \overline{V}\to \C , & E(\sL^-_t)=\overline{V}\otimes_{A} V \to \C  \\  E(\sR^+_t)= \C\to V\otimes_{\overline{A}} \overline{V} & E(\sR^+_t)=  \C\to \overline{V}\otimes_A V \\ 
E(S_t)\in V\otimes_{A\otimes \overline{A}} \overline{V}
  \end{array}
\eeq
Translating back to Definition~\ref{defn:QMsystem}, this introduces an $A$-action on $\cH$ and requires the Hamiltonian $H$ be $A$-linear. 
\end{ex}

\begin{rmk}\label{rmk:symmetry}
Operators commuting with the Hamiltonian $H$ are \emph{symmetries} of a quantum theory. Hence, the examples of $A$-twisted field theories in the previous example can also be viewed as theories in which the algebra $A$ acts by symmetries. 
\end{rmk}

\begin{ex}[Spin theories]\label{ex:spintheories}
To handle fermions (as in Example~\ref{ex:fermion}), one enhances the bordism category $\EBord_1$ to $\EBord_1^{\rm Spin}$ where bordisms are now Euclidean 1-manifolds with spin structure. The appropriate target for the spin field theories enhances $\TA$ to a bicategory of $\Z/2$-graded algebras, $\Z/2$-graded bimodules, and grading-preserving bimodule maps equipped with the $\Z/2$-graded tensor product. Spin 1-manifolds have a $\Z/2$-automorphism from the spin flip, and similarly $\Z/2$-graded algebras and bimodules have a $\Z/2$-automorphism from their grading involution $(-1)^F$; see~\eqref{eq:susyalgebra} below. More formally, $\EBord^{\rm Spin}_1$ and $\TA$ have canonical $B\Z/2$-actions. A \emph{spin field theory} is then a reflection positive, $B\Z/2$-equivariant functor as in~\eqref{eq:repof1EB0}, and a \emph{twisted spin theory} is a $B\Z/2$-equivariant natural transformation relative to a $B\Z/2$-equivariant functor as in~\eqref{eq:twisted1EB0}. This $B\Z/2$-equivariance is a version of the spin statistics theorem~\cite[
Definition 6]{HST}. 

The free fermions of Example~\ref{ex:fermion} provide an important example of a twisted spin field theory theory with $A=\cCl_n$ in~\eqref{eq:twist1EB0} defining the twist, and $E(\pt^+)=V$ a $\cCl_n$-module defining the twisted spin theory. The semigroup is the identity, and hence the 1-dimensional free fermions determine a \emph{topological} twisted spin theory. Following the previous remark, the $\cCl_n$-action is the quantization of the classical translational symmetries that preserve the action functional~\eqref{Eq:FFaction}. Incorporating these symmetries with the correct signs in the bordism framework necessitates the formalism of spin field theories.
\end{ex}


%
%

 

\section{Supersymmetric quantum mechanics and K-theory}\label{sec:SQM}

The goal of this section is to explain the following analog of Conjecture~\ref{conj}.

\begin{conj}[Stolz--Teichner]\label{thm11}
For each $n\in \Z$ there exist dashed arrows 
\beq
\begin{tikzpicture}[baseline=(basepoint)]
    \node (A) at (0,0) {$  \left\{\begin{array}{c}  {\rm families\ of\ } \\ {n\hbox{-}{\rm dimensional\ spin }^c}\\ {\rm  \ manifolds } \ M\to X \end{array}\right\}$};  
  \node (B) at (8,-1) {$\K^{-n}(X)$};
  \node (C) at (0,-2.5) {$\left\{\begin{array}{c}  {\rm degree} \ -\!n \\ 1d \ {\rm supersymmetric }\\ {\rm field\ theories\ over\ } X\end{array}\right\}$};
  \draw[->] (A) to node [above=4pt] {${\rm topological\ index}$} (B);
  \draw[->,dashed] (A) to node [left] {${\rm quantize}$} (C);
  \draw[->>,dashed] (C) to node [below] {${\rm cocycle}_{-n}$} (B);
    \path (0,-1) coordinate (basepoint);
    \end{tikzpicture}\label{eq:Thm11}
\eeq
making the triangle commute, where ``topological index" is the spin$^c$ orientation of Atiyah--Bott--Shapiro~\cite{ABS}, ``quantize" is quantization of supersymmetric mechanics on~$M$, and ``${\rm cocycle}_n$" sends a quantum theory over $X$ to a class in~$\K^n(X)$. For negative $n$, the empty manifold is the only spin manifold and the content of~\eqref{eq:Thm11} is the existence of a cocycle map.
\end{conj}

Subject to some technical caveats, Conjecture~\ref{thm11} is a theorem. Indeed, it has been verified for the definition of 1-dimensional supersymmetric field theory from \cite[\S3.2]{ST04} and \cite{HST}, reviewed in \S\ref{sec:littlethm} and Hypothesis~\ref{hyp:SQMandBord} below. However, it is believed that these prior arguments will not generalize well to the $2|1$-dimensional case, and so variations on Conjecture~\ref{thm11} remain open; one such variation is described in Remark~\ref{rmk:Thm11}.

\subsection{Motivation from supersymmetric quantum mechanics}\label{sec:space11}
Given a quantum mechanical system $(\cH,H)$ in the sense of Definition~\ref{defn:QMsystem}, suppose that $\cH=\cH^+\oplus \cH^-$ is equipped with a $\Z/2$-grading and that $H$ respects the grading (i.e., is an even map). Equivalently, $\mathcal{H}$ is equipped with an involution $(-1)^F\colon \cH\to \cH$ with eigenspaces
\beq\label{eq:susyalgebra}
&&(-1)^F=\left\{\begin{array}{cl} +1 & {\rm on } \ \cH^+ \\ -1 & {\rm on } \ \cH^- \end{array}\right. \quad {\rm and} \quad  H(-1)^F=(-1)^FH.
\eeq
\emph{Supersymmetric quantum mechanics} enhances these data by $\mathcal{N}$ odd operators $Q_i$ satisfying
\beq\label{eq:Lie11}
&&Q_i^2=H, \quad [Q_i,Q_j]=Q_iQ_j+Q_jQ_i=0 \ i\ne j \quad [(-1)^F,Q_i]=0\quad i\in \{1,2,\dots, \mathcal{N}\}
\eeq
where $[-,-]$ is the graded commutator \cite[\S2]{Wittenbreaking}. The $Q_i$ are \emph{supersymmetry operators}. 

\begin{ex}[$\mathcal{N}=2$ supersymmetric quantum mechanics on an oriented manifold]\label{rmk:SUSYmorse}
The classic example of supersymmetric quantum mechanics takes as its space of states~$\cH=\Omega^\bullet(M)$, differential forms on an oriented Riemannian manifold~\cite[\S2]{susymorse}. The Hamiltonian is then defined to be the Hodge Laplacian $H=dd^*+d^*d$, which has a pair of commuting odd square roots, 
\beq
Q_1=d+d^*, \quad Q_2=\sqrt{-1}(d-d^*), \quad [Q_1,Q_2]=0, \quad Q_1^2=Q_2^2=H.\label{eq:Nis2}
\eeq
This is the data of $\mathcal{N}=2$ supersymmetry. Forgetting one of the $Q_i$ in~\eqref{eq:Nis2} extracts an $\mathcal{N}=1$ theory. Other $\mathcal{N}=1$ theories come from (twisted) Dirac operators, see Example~\ref{ex:spinmanifold}.
\end{ex}

The focus in Conjecture~\ref{thm11} (and below) is on $\mathcal{N}=1$ supersymmetry. In this case, the operators~\eqref{eq:Lie11} determine a representation of the super Lie algebra of the super Lie group~$\R^{1|1}$ with multiplication\footnote{This formula is to be understood in terms of the functor of points of~$\R^{1|1}$, e.g., see~\cite{DM}.}
\beq
(t,\theta)\cdot (s,\eta)=(t+s+\theta\eta,\theta+\eta)\qquad (t,\theta),(s,\eta)\in \R^{1|1}. \label{eq:superEuc}
\eeq
Including along $\theta=0$, one obtains the abelian subgroup $\R<\R^{1|1}$ for the usual addition on~$\R$. Extending along this inclusion provides a generalization of the usual time-evolution~\eqref{eq:semigroup} to a representation of a \emph{supersemigroup}~\cite[\S5]{HST} \cite[\S3.2]{ST04}
\beq
\R^{1|1}_{\ge 0}\to \End(\cH)\qquad (t,\theta)\mapsto e^{-tH+\theta Q},\label{eq:N1super}
\eeq
where $(t,\theta)$ are coordinates on the sub supermanifold $\R^{1|1}_{\ge 0}\subset \R^{1|1}$ with $t\ge 0$.

\begin{defn}\label{defn:SQMsystem}
A \emph{quantum mechanical system with ($\mathcal{N}=1$) supersymmetry} is a pair $(\cH,Q)$ where~$Q$ is an odd unbounded self-adjoint operator on the Hilbert space $\cH$ generating a trace class supersemigroup representation~\eqref{eq:N1super}. Systems~$(\cH,Q)$ and~$(\cH',Q')$ are \emph{equivalent} if there exists a unitary map $u\colon \cH\to \cH'$ with $Q'=uQu^{-1}$. 
\end{defn}

We expand on this below.

\subsection{The space of degree~$n$ quantum systems with $\mathcal{N}=1$ supersymmetry} 

Analogous to the space of Hamiltonians described in \S\ref{sec:QMdef}, we extend Definition~\ref{defn:SQMsystem} to a space---and ultimately a stack---of quantum systems with supersymmetry. 
With an eye towards the $\Z$-grading on K-theory, we also incorporate Clifford module structures; physically these arise from an algebra of (projective) symmetries of the free fermions, see Example~\ref{ex:fermion} and Remark~\ref{rmk:symmetry}. 

Fix an infinite-dimensional separable $\Z/2$-graded Hilbert space~$\cH$ whose even and odd subspaces are both infinite-dimensional. Set $\cH_n:=\cCl_n\otimes \cH$ where $\cCl_n$ is the $n$th Clifford algebra. Let $Q$ be a Clifford linear odd self-adjoint unbounded operator\footnote{An odd unbounded operator on a $\Z/2$-graded Hilbert space has domain $\dom(Q)\subset \cH$ a graded subspace and $Q\colon \dom(Q)\to \cH$ an odd linear map.} on $\cH_n$. We observe that $H=Q^2$ is a Clifford linear, nonnegative, even operator. Applying the spectral theorem to~$H=Q^2$ results in a collection of mutually orthogonal eigenspaces $\{\cH_\lambda\subset \cH_n\}_{\lambda_{\ge 0}\in \R}$ as in Figure~\ref{fig1}. These subspaces are finite-rank $\Z/2$-graded $\cCl_n$-modules with the action of an odd operator $Q_\lambda$ with eigenvalues $\pm\sqrt{\lambda}$. When $\lambda\ne 0$, the operator $\frac{1}{\sqrt\lambda} Q$ generates a $\cCl_{-1}$-action on $\cH_\lambda$ commuting with the $\cCl_n$-action. Hence, the $\cH_\lambda$ for $\lambda\ne 0$ acquire the structure of a $\cCl_{-1}\otimes \cCl_n$-module. Deformations of these configurations (e.g., see Figure~\ref{fig2}) therefore recover the Atiyah--Bott--Shapiro description of the coefficients of $\pi_*\K^n$ as isomorphism classes of $\cCl_n$-modules modulo those that extend to $\cCl_{-1}\otimes \cCl_n$-modules~\cite{ABS}. 

\begin{figure}
\beq\nonumber
\begin{tikzpicture}[baseline=(basepoint)];
\node (A) at (1,1.7) {$\bullet$};
\node (AA) at (1.5,1.7) {$\cH_{\lambda_3}$};
\node (B) at (1,1.1) {$\bullet$};
\node (BB) at (1.5,1.1) {$\cH_{\lambda_2}$};
\node (C) at (1,0) {$\bullet$};
\node (CC) at (1.5,0) {$\cH_{0}$};
\node (D) at (1,.7) {$\bullet$};
\node (DD) at (1.5,.7) {$\cH_{\lambda_1}$};
\draw[thick,<->] (1,-.3) to (1,2);
\draw[<->,thick] (2,1.7) to  node [below] {$_{\sdim\cH_{\lambda_3}=2|2}$} (3.25,1.7);
\draw[<->,thick] (2,1.15) to  node [below] {$_{\sdim\cH_{\lambda_2}=1|1}$} (3.25,1.15);
\draw[<->,thick] (2,0) to node [below] {$_{\sdim\cH_{\lambda_0}=2|0}$} (3.25,0);
\draw[<->,thick] (2.,.65) to  node [below] {$_{\sdim\cH_{\lambda_1}=2|2}$} (3.25,.65);
\draw[thick,<->] (4,-.3) to (4,2);
\node (E) at (3.5,0) {};
\node (1) at (4.5,0) {$\circ$};
\node (2) at (5,0) {$\circ$};
\node (EE) at (5.5,0) {};
\node (F) at (3.5,1.7) {};
\node (1) at (4.75,1.7) {$\circ$};
\node (1) at (4.5,1.7) {$\circ$};
\node (1) at (5,1.7) {$\times$};
\node (1) at (5.25,1.7) {$\times$};
\node (FF) at (5.5,1.7) {};
\node (G) at (3.5,1.1) {};
\node (1) at (4.5,1.1) {$\circ$};
\node (1) at (5,1.1) {$\times$};
\node (GG) at (5.5,1.1) {};
\node (H) at (3.5,.7) {};
\node (1) at (4.5,.7) {$\circ$};
\node (1) at (4.75,.7) {$\circ$};
\node (1) at (5,.7) {$\times$};
\node (1) at (5.25,.7) {$\times$};
\node (HH) at (5.5,.7) {};
\draw[-,thick] (E) to (EE);
\draw[-,thick] (F) to (FF);
\draw[-,thick] (G) to (GG);
\draw[-,thick] (H) to (HH);

\draw[->,thick,dashed,bend left=15] (5.7,.65) to (6.75,.35);
\draw[thick,<->] (7.5,-.3) to  (7.5,2);
\node (J) at (7,0) {};
\node (1) at (8,0) {$\circ$};
\node (2) at (8.5,0) {$\circ$};
\node (JJ) at (9,0) {};
\node (K) at (7,1.7) {};
\node (1) at (8,1.7) {$\circ$};
\node (1) at (8.25,1.7) {$\circ$};
\node (1) at (8.5,1.7) {$\times$};
\node (1) at (8.75,1.7) {$\times$};
\node (KK) at (9,1.7) {};
\node (L) at (7,1.1) {};
\node (1) at (8,1.1) {$\circ$};
\node (1) at (8.5,1.1) {$\times$};
\node (LL) at (9,1.1) {};
\node (M) at (7,.3) {};
\node (1) at (8,.3) {$\circ$};
\node (1) at (8.25,.3) {$\circ$};
\node (1) at (8.5,.3) {$\times$};
\node (1) at (8.75,.3) {$\times$};
\node (MM) at (9,.3) {};
\draw[-,thick] (J) to (JJ);
\draw[-,thick] (K) to (KK);
\draw[-,thick] (L) to (LL);
\draw[-,thick] (M) to (MM);

\draw[->,thick,dashed, bend left=10] (9.2,.3) to (10.25,.1);
\draw[thick,<->] (11,-.3) to  (11,2);
\node (J) at (10.5,0) {};
\node (1) at (11.2,0) {$\circ$};
\node (1) at (11.4,0) {$\circ$};
\node (2) at (12,0) {$\times$};
\node (JJ) at (12.5,0) {};
\node (K) at (10.5,1.7) {};
\node (1) at (11.5,1.7) {$\circ$};
\node (1) at (11.75,1.7) {$\circ$};
\node (1) at (12,1.7) {$\times$};
\node (1) at (12.25,1.7) {$\times$};
\node (KK) at (12.5,1.7) {};
\node (L) at (10.5,1.1) {};
\node (1) at (11.5,1.1) {$\circ$};
\node (1) at (12,1.1) {$\times$};
\node (LL) at (12.5,1.1) {};
\node (1) at (11.6,0) {$\circ$};
\node (1) at (11.8,0) {$\circ$};
\node (1) at (12,0) {$\times$};
\node (1) at (12.2,0) {$\times$};
\draw[-,thick] (J) to (JJ);
\draw[-,thick] (K) to (KK);
\draw[-,thick] (L) to (LL);
\path (0,.5) coordinate (basepoint);
\end{tikzpicture}\qquad 
\eeq

\caption{Generalizing Figure~\ref{fig1}, supersymmetry enhances the label $\cH_\lambda\subset \cH$ by a $\Z/2$-grading and odd operator $Q_\lambda$ satisfying $Q_\lambda^2=\lambda \cdot 1$ where $\lambda\ge 0$. Following \cite{Wittenbreaking} we record this by $\circ$'s counting the number of even dimensions and $\times$'s counting odd dimensions (representing bosonic and fermionic states, respectively). As $Q_\lambda$ is invertible for $\lambda\ne 0$, the number of even dimensions must equal the number of odd dimensions on $\cH_\lambda$ with $\lambda\ne 0$. Hence, the difference between even and odd dimensions in the zero eigenspace is the deformation invariant $[\cH_0]\in \K(\pt)$. 
}
\label{fig2}
\end{figure}
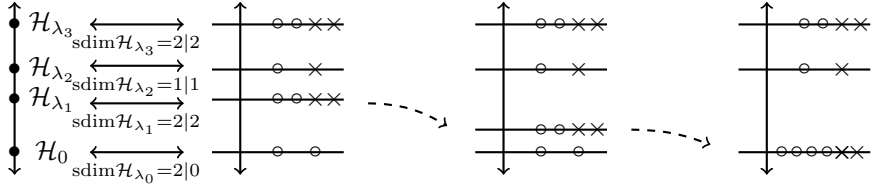

%

\begin{defn} \label{defn:sHam}
Let $\sHam^n$ denote the space of odd self-adjoint $\cCl_n$-linear unbounded operators $Q$ that generate trace class supersemigroup representations~\eqref{eq:N1super}.
\end{defn}
\begin{rmk}
The analytical conditions under which an unbounded operator $Q$ determines a trace class supersemigroup representation are given in~\cite[\S5]{HST}. 
\end{rmk} 

\begin{defn} For a $\Z/2$-graded Hilbert space, define the (super) \emph{unitary group}
$$
\U^\pm(\cH)\simeq \{u\in \U(\cH)\mid u=\left\{\left[\begin{array}{cc} * & 0 \\ 0 & *\end{array}\right] \ {\rm or} \ \left[\begin{array}{cc} 0 & * \\ * & 0\end{array}\right]\right\}
$$
for block diagonal form relative to the decomposition~$\cH=\cH^+\oplus \cH^-$ into even and odd subspaces. The group $\U^\pm(\cH)$ is graded (e.g., in the sense of~\cite{Freedalg}) via the homomorphism 
$$
\U^\pm(\cH)\to \pi_0(\U^\pm(\cH))\simeq \Z/2
$$
that reads off whether a unitary map preserves or reverses the grading on $\cH$. Let $\PU^\pm(\cH):=\U^\pm(\cH)/\U(1)$ denote the projectivization. Similarly for a graded unitary $\cCl_n$-module $\cH_n$, define $\U^\pm(\cH_n)$ as the unitary transformations commuting with the Clifford action and $\PU^\pm(\cH_n)=\U^\pm(\cH_n)/\U(1)$ as the projectivization. 
\end{defn}



\begin{defn}\label{defn:SQM}
 Define the topological groupoid $\sqm^n:=\sHam^n\sq \PU^\pm(\cH_n)$ as the quotient for the action induced by the conjugation action of the unitary group. 
\end{defn}

We have the following supersymmetric generalization of Definition~\ref{defn:sym1}. 

\begin{defn}\label{defn:sym11}
For $G$ a compact Lie group, the stack of \emph{degree~$n$ supersymmetric quantum mechanical systems with $G$-symmetry} is the mapping stack $\Map(\pt\sq G,\sqm^n)$. 
\end{defn}

Unpacking this, a point in $\Map(\pt\sq G,\sqm^n)$ is a homomorphism $\alpha\colon G\to \PU^\pm(\cH)$ defining a central extension $G^\alpha$ of $G$ with a unitary representation $G^\alpha\to \U^\pm(\cH)$ (analogous to~\eqref{eq:pullbacksequence}), together with a supersymmetry generator $Q$ on $\cH$ commuting with the $G^\alpha$-action. The composition $G^\alpha\to \U^\pm(\cH)\to \Z/2$ endows the group~$G^\alpha$ with a grading that reads off which elements act by grading preserving or grading reversing automorphisms of~$\cH$. 

\subsection{Supersymmetric quantum mechanics and K-theory}\label{sec:littlethm}

The Atiyah--Bott--Shapiro description of the graded ring $\pi_{-*}\K=\K^*(\pt)$ \cite{ABS} gives a map 
\beq\label{eq:Wittenindex}
&&\sHam^n\to \K^n(\pt), \qquad Q\mapsto [{\rm Index}(Q)]:=[{\rm Ker}(Q)]=[\cH_0]
\eeq
sending a supersymmetry generator $Q$ to its kernel, viewed as a finite rank $\cCl_n$-module. 

\begin{defn} 
The \emph{Witten index} of $Q\in \sHam^n$ is the image under~\eqref{eq:Wittenindex}. 
\end{defn}

In fact, we obtain the following.

\begin{thm}[\cite{HST},\cite{PokmanPhD}]\label{thm:HamK}
The space $\sHam^n$ represents the functor $\K^n(-)$. 
\end{thm}

\begin{proof}[Proof idea.]
The space $\sHam^n$ is homotopy equivalent to a space of Clifford linear Fredholm operators that is known to represent the $n$th space in the K-theory spectrum~\cite{AtiyahSingerskew}. Indeed, the growth conditions on the eigenvalues of operators in $Q\in \sHam^n$ ensure that the operator $Q/(1+Q^2)$ defined via functional calculus is Fredholm. 
\ep

\begin{ex}[Topological theories]\label{ex:TQFTQ}
Generalizing Example~\ref{ex:TQFT}, any finite-dimensional $\Z/2$-graded submodule $V\subset\mathcal{H}_n$ determines a point in $\sHam^n$ by taking $Q= 0$ with domain~$V$. This corresponds to the configuration supported at $0\in \R_{\ge 0}$ with label $\mathcal{H}_0=V$. We observe that the natural inclusion 
\beq\label{Eq:1dgroupcomplete}
\left\{\begin{array}{c} {\rm degree\ } n \\ {\rm topological\ theories}\end{array}\right\} \xhookrightarrow{{\rm group\ completion}} \sHam^n\to \sqm^n
\eeq
can be identified with the group completion of $\cCl_n$-modules. 
\end{ex}

\begin{prop} \label{hyp:SQM}
$\{\sqm^n\}_{n\in \Z}$ represents twisted K-theory of topological stacks, and in particular twisted equivariant $\K$-theory of $G$-spaces for $G$ a compact Lie group~\cite{FHTI}. \end{prop}

\begin{proof}[Proof sketch]
Given a map $\alpha$ of topological stacks, 
\beq\label{eq:twistedKtheory}
\X\xrightarrow{\alpha} \sqm^n=\sHam^n\sq \PU^\pm(\cH_n)
\eeq
postcomposing with the projection determines a map $\tau\colon \X\to \pt\sq \PU^\pm(\cH_n)$ classifying a projective Hilbert bundle $\mathbb{P}\mathcal{H}^\tau\to \X$ with Clifford action (also known as a \emph{twisted Hilbert bundle}~ \cite[\S3.2]{FHTI}). This bundle is furthermore classified by the pullback of the graded gerbe $\pt\sq \U^\pm(\cH_n)\to\pt \sq  \PU^\pm(\cH_n)$ along $\tau$. The remaining data of the map $\alpha$ is a section of the associated bundle over $\X$ with fiber $\s\Ham_n$. Since  $\s\Ham_n$ represents $\K^n(-)$, the homotopy class of this section is a twisted K-theory class. Specializing to the quotient stack $\X=X\sq G$ of a (nice) space $X$ by the action of compact Lie group $G$ recovers the $G$-equivariant twisted $\K$-theory of $X$~\cite[\S3.5.4]{FHTI}.
\ep


\subsection{Supersymmetric quantum mechanics and the index theorem}

\begin{ex}[Supersymmetric quantum mechanics on a spin$^c$ manifold]\label{ex:spinmanifold}
Let $M$ be a Riemannian manifold and consider the classical theory of paths $x\colon \R\to M$ together with a (fermion) $\psi\in \Gamma(\R;x^*TM)$ and action 
\beq\label{eq:SUSYMech}
S(x,\psi)=\frac{1}{2}\int_\R (\| \dot x \| +\langle \psi,\nabla_{\dot{x}}\psi\rangle)dt. 
\eeq
Suppose that $M$ has a chosen spin$^c$ structure with spinor bundle $\bS:=P\times_{\Spin(n)}\cCl_n$ supporting the $\cCl_{-n}$-linear\footnote{The bundle of right $\cCl_n$-modules $\bS=P\times_{\Spin(n)}\cCl_n$ determines a bundle of left $\cCl_{-n}=\cCl_n^\op$-modules.} Dirac operator $\slashed{D}$. Then the classical theory~\eqref{eq:SUSYMech} quantizes with state space $\cH=\Gamma(\bS)$ and Hamiltonian $H=\slashed{D}{}^2$~\cite[pages~43-44]{Freed5}. Taking $Q=\slashed{D}$ defines a degree~$-n$ supersymmetric quantum system in the sense of Definition~\ref{defn:SQM}, giving 
\beq\label{Eq:11quantize}
\left\{\begin{array}{c}  {n\hbox{-}{\rm dimensional\ spin }^c}\\ {\rm  \ manifolds } \ M \end{array}\right\}\to \sqm_{-n},\qquad M\mapsto \slashed{D}.
\eeq
If $M$ has a $G$-action preserving the spin$^c$ structure, then $\Gamma(\bS)$ is naturally a $G$-representation. If instead $W_3(TM)\in \H_G^3(M;\Z)$ pulls back from the cohomology of $BG$ (i.e., the equivariant cohomology of the point), then $\Gamma(\bS)$ only has a projective $G$-action. Both cases lead to a supersymmetric quantum mechanical system with $G$-symmetry in the sense of Definition~\ref{defn:sym11}. 
\end{ex}

\begin{proof}[Proof sketch of Conjecture~\ref{thm11}]
If one defines a degree~$-n$, 1-dimensional supersymmetric Euclidean field theory over $X$ to be a map $X\to \Ham(\cH_{-n})$ (compare  \cite[Definition~3.2.10]{ST04}), the cocycle map in~\eqref{eq:Thm11} follows from Proposition~\ref{hyp:SQM}. Given a family of spin$^c$ manifolds $M\to X$, the fiberwise $\cCl_{-n}$-linear Dirac operator gives a families version of~\eqref{Eq:11quantize}, producing a map $X\to \Ham(\cH_{-n})$ corresponding to the families analytic index~\cite[Remark 3.2.21]{ST04}. Commutativity of the triangle~\eqref{eq:Thm11} is the families index theorem~\cite{AtiyahSinger4}. 
\end{proof}

\begin{rmk}
The physics argument for the index theorem considers the limits of the time-evolution semigroup of supersymmetric quantum mechanics on~$M$,
\beq\label{eq:partition function2}
&&\sTr_{\Cl_{-n}}(e^{-t\slashed{D}^2})=\left\{ \begin{array}{ll} \displaystyle{\rm Index}(\slashed{D})  & t\to \infty \\ \displaystyle \int_X\hat{A}(X) & t\to 0.\end{array}\right.
\eeq
The McKean--Singer theorem~\cite{McKeanSinger} shows that this super trace is constant in $t$ and identifies the $t\to\infty$ limit with the (Witten) index of $Q=\slashed{D}$. On the other hand, the $t\to 0$ limit can be computed using path integral techniques. Specifically, supersymmetric localization identifies the path integral with the integral of the $\hat{A}$-form on~$M$~\cite{susymorse,Alvarez}. For mathematical treatments, see~\cite{AtiyahCircular,JonesPetrack,BGV}.
\end{rmk}

The topological pushforward in K-theory comes from the Thom isomorphism, which also has a description in terms of supersymmetric quantum mechanics. 



\begin{ex}[Supersymmetric deformations of free fermions]\label{ex:freefermion11}
By Example~\ref{ex:TQFTQ}, the free fermions from Examples~\ref{ex:fermion} and~\ref{ex:spintheories} have a canonical extension to a degree~$n$ supersymmetric theory with state space a graded $\cCl_n$-module $V$, and supersymmetry generator~$Q=0$. Consider the modified free fermion classical field theory~\eqref{Eq:FFaction}
\beq\label{eq:FFsource}
\psi\colon \R\to \R^n,\quad S(\psi)=\frac{1}{2}\int_\R (\langle \psi,\partial_t\psi\rangle +\langle x,\psi\rangle) dt
\eeq
that adds the source term $\langle x,\psi\rangle$ to the action for a fixed $x\in \R^n$. Quantization of this system leads to state space $V=\cCl_n$ as a left $\cCl_n$-module together with time-evolution operator determined by the supersymmetry generator $Q={\rm cl}_x\colon V\to V$ where ${\rm cl}_x$ is right Clifford multiplication by $x\in \Pi \R^n\subset \cCl_n$. Varying $x\in \R^n$, we obtain a map 
\beq\label{eq:suspension}
\R^n\to \sqm^n,\qquad x\mapsto Q={\rm cl}_x
\eeq
that at $x=0$ is a topological theory. One can identify~\eqref{eq:suspension} with the $n$th suspension class in K-theory~\cite{ABS}. There is an $\O_n$-action on~\eqref{eq:FFsource} that rotates $\psi$ and the source term~$x$, and the map~\eqref{eq:suspension} is (projectively) $\O_n$-equivariant. This gives the generalization of~\eqref{eq:FFsymmetry}
\beq\label{eq:Thom}
&&\Th_n\colon \R^n\sq \O_n \to \sqm^n=\sHam^n\sq \PU^\pm(\cH), \qquad \begin{array}{c} \O_n\to \PU^\pm(\cH),\\ x\mapsto \cl_x=Q(x), \ x\in \R^n\end{array}
\eeq
where the map $\O_n\to \PU^\pm(\cH_n)$ classifies the Pin$^c$-extension. Under Proposition~\ref{hyp:SQM}, the twisted K-theory class underlying~\eqref{eq:Thom} is the $n$th universal (twisted equivariant) Thom class $[\Th_n]\in \K^{w_1+W_3}(\R^n\sq \O_n)_c\simeq \K^{w_1+W_3}_{\O_n}(\R^n)_c$ \cite[\S3.6]{FHTI} \cite[Remark 3.2.22]{ST04}.
\end{ex}

\begin{rmk}
The Mathai--Quillen formalism~\cite{MathaiQuillen} computes the partition function of~\eqref{eq:FFsource}, recovering a universal Thom class determined by the (inverse) $\hat{A}$-form on $\mathfrak{o}_n$.
\end{rmk}

\subsection{Supersymmetric quantum mechanics in the bordism framework}\label{sec:11bord}

One can repackage supersymmetric quantum mechanics as a representation of a $1|1$-Euclidean bordism category $\EBord_{1|1}$ that generalizes $\EBord_1$. The objects of $\EBord_{1|1}$ are finite disjoint unions of $\bullet=\spt^+$ and $\circ=\spt^-$, the \emph{super point} with its two choices of orientation. Morphisms are pictured in Figure~\ref{fig:11EBord}, where $(t,\theta)\in \R^{1|1}_{\ge 0}$ encodes a ``super length." Analogously to before, we find relations associated to compositions of bordisms, e.g.,
\beq\label{eq:SSrep}
&&\sI_{s,\eta}^\pm\circ \sI_{t,\theta}^\pm \simeq \sI_{s+t\pm \eta\theta,\eta+\theta}^\pm,
\eeq
giving the geometric counterpart to multiplication in the supersemigroups $\R_{\ge 0}^{1|1}$ from~\eqref{eq:N1super}. 

\begin{rmk}
Super Euclidean geometry in dimension~$1|1$ is determined by the super Lie group~\eqref{eq:superEuc} in the following sense: one declares that the group $\R^{1|1}$ acts on the supermanifold~$\R^{1|1}$ by Euclidean isometries. This is analogous to how Euclidean geometry in dimension~1 is determined by declaring the group $\R$ to be the oriented symmetry group of~$\R$, acting through translations. There is additional $1|1$-Euclidean isometry determined by $\theta\mapsto -\theta$, playing the role of the spin flip. Consequently, $1|1$-Euclidean field theories generalize the spin Euclidean field theories from Example~\ref{ex:spintheories}. We impose the corresponding $B\Z/2$-equivariance property implicitly below, referring to \cite[\S4.3]{ST11} for details. 
\end{rmk}

Following~\eqref{eq:twist1EB0}, we construct a twist for each Clifford algebra $\cCl_n$.

\begin{defn}[Degree twist]\label{defn:11degreetw}
The \emph{degree twist} is determined by the assignments 
\beq\label{eq:degreetwist11}
\alpha_n\colon \EBord_{1|1} \to \TA,\qquad  \begin{array}{c} \alpha_n(\spt^+)=\cCl_n, \ \alpha_n(\spt^-)=\cCl_n^\op=\cCl_{-n}, \\ \alpha_n(\sI^\pm_{t,\theta})=_{\cCl_n}\!(\cCl_n)_{\cCl_n},  
\alpha_n(\sL^+_{t,\theta})=(\cCl_n)_{\cCl_n\otimes \cCl_n^\op},  \\  \alpha_n(S_{t,\theta}^\pm)=\cCl_n\otimes_{\cCl_n\otimes \cCl_n^\op} \cCl_n,
\\
\alpha_n(\sR^+_{t,\theta})= _{\cCl_n\otimes \cCl_n^\op} \cCl_n,  \end{array}
\eeq
with reflection structure from the $*$-structure on the Clifford algebra $*\colon \cCl_n\xrightarrow{\sim} \overline{\cCl}{}_n^\op=\overline{\cCl}_{-n}$. These values are independent of $(t,\theta)$, and hence the degree twist is topological. The values of a degree~$n$ twisted field theory can be read off analogously to~\eqref{eq:twisted1EB0}. The Morita equivalences $\cCl_n\otimes \cCl_m\simeq \cCl_{n+m}$ yield isomorphisms of twists $\alpha_n\otimes \alpha_m\simeq\alpha_{n+m}$ for $n,m\in \Z$.
\end{defn}

\begin{figure}
\beq
\begin{tikzpicture}[baseline=(basepoint)];
\draw [fill]  (0,1.5) circle [radius=0.1] node [black,left=4] {$\pt^+$};
\draw  [fill] (0,.5,0) circle [radius=0.1] node [black,left=4] {$\pt^+$};
\draw[decoration={markings,
        mark=at position \halfway with \arrow{<}},
        postaction=decorate,thick] (0,1.5) to (0,.5);
\node (Iplus) at (0,0) {$\sI_{t,\theta}^+$};

\draw (1.5,1.5) circle [radius=0.1] node [black,left=4] {$\pt^-$};
\draw (1.5,0.5) circle [radius=0.1] node [black,left=4] {$\pt^-$};
\draw[decoration={markings,
        mark=at position \halfway with \arrow{>}},
        postaction=decorate,thick] (1.5,1.5) to  (1.5,0.5);
\node (Iminus) at (1.5,0) {$\sI_{t,\theta}^-$};

\draw  (2.75,.5) circle [radius=0.1] node [black,left=4] {};
\draw[fill]  (4.25,.5) circle [radius=0.1] node [black,right=4] {};
\draw[decoration={markings,
        mark=at position \halfway with \arrow{<}},
        postaction=decorate,thick] (2.75,.5) to [out=90,in=180] (3.5,1.5);
        \draw[decoration={markings,
        mark=at position \halfway with \arrow{<}},
        postaction=decorate,thick] (3.5,1.5) to [out=0,in=90] (4.25,0.5);

\node (Rminus) at (3.65,0) {$\sL_{t,\theta}^+$};

\draw  (5.25,.5) circle [radius=0.1] {};
\draw [fill] (6.75,.5) circle [radius=0.1] {};
\draw[decoration={markings,
        mark=at position \halfway with \arrow{>}},
        postaction=decorate,thick] (5.25,.5) to [out=90,in=180] (6,1.5);
        \draw[decoration={markings,
        mark=at position \halfway with \arrow{>}},
        postaction=decorate,thick] (6,1.5) to [out=0,in=90] (6.75,0.5);

\node (Rplus) at (6.15,0) {$\sL_{t,\theta}^-$};

\draw [fill] (7.75,1.5) circle [radius=0.1] node [black,left=4] {};
\draw  (9.25,1.5) circle [radius=0.1] node [black,right=4] {};
\draw[decoration={markings,
        mark=at position \halfway with \arrow{>}},
        postaction=decorate,thick]  (7.75,1.5) to [in=180,out=270] (8.5,.25);
\draw[decoration={markings,
        mark=at position \halfway with \arrow{>}},
        postaction=decorate,thick]  (8.5,.25) to [in=270,out=0] (9.25,1.5);

\node (Lplus) at (8.65,0) {$\sR_{t,\theta}^+$};

\draw  (10.25,1.5) circle [radius=0.1]{};
\draw [fill](11.75,1.5) circle [radius=0.1] node [black,right=4] {};
\draw[decoration={markings,
        mark=at position \halfway with \arrow{<}},
        postaction=decorate,thick]  (10.25,1.5) to [in=180,out=270] (11,.25);
\draw[decoration={markings,
        mark=at position \halfway with \arrow{<}},
        postaction=decorate,thick]  (11,.25) to [in=270,out=0] (11.75,1.5);

\node (Lminus) at (11.15,0) {$\sR_{t,\theta}^-$};

\draw[decoration={markings,
        mark=at position \halfway with \arrow{>}},
        postaction=decorate,thick] (4.75,1.85) circle [radius=.45] node [right=12] {$S_{t,\theta}^\pm$};


\path (0,.75) coordinate (basepoint);
\end{tikzpicture}\nonumber
\eeq
\caption{Morphisms in $\EBord_{1|1}$ are finite disjoint unions of the above generators. These enhance the generators in Figure~\ref{fig:EBord} by adding spin structures and replacing the length parameter $t\in \R_{\ge 0}$ by a ``super length" $(t,\theta)$ valued in $\R^{1|1}_{\ge 0}$. }
\label{fig:11EBord}
\end{figure}
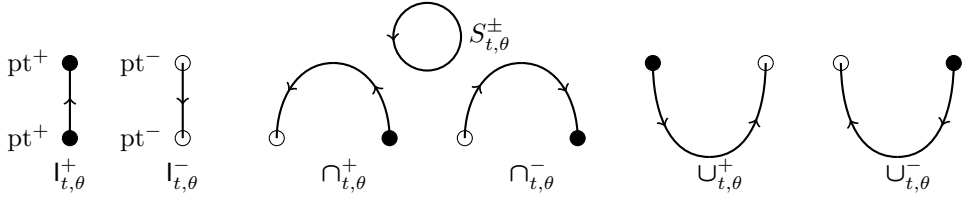

\begin{rmk}
The Clifford algebras are invertible objects in the 2-category $\TA$, and hence are fully dualizable. By the cobordism hypothesis~\cite{Lurie_cob}, $\cCl_n$ determines a fully extended 2-dimensional topological field theory. The values of this theory on 0- and 1-manifolds are given by~\eqref{eq:degreetwist11}~\cite{Sam}. Hence, the degree twist is a truncation of this fully-extended topological field theory with reflection structure inherited the $*$-structure on~$\cCl_n$. 
\end{rmk}

Let $\QFT_{1|1}^\alpha(\X)$ denote the groupoid of $\alpha$-twisted reflection positive $1|1$-dimensional Euclidean field theories over $\X$ (see Definition~\ref{def:QFT}). Reflection positivity is essentially the same as in Definition~\ref{defn:RP1}. 
We start by analyzing the groupoid $\QFT^n_{1|1}:=\QFT_{1|1}^{\alpha_n}(\pt)$ for the degree twist~\eqref{eq:degreetwist11}.

\begin{hyp}[{Compare \cite[Proposition 3.2.6]{ST04} and \cite[Corollary 6.38]{HST}}] \label{hyp:SQMandBord}
There is an epimorphism of stacks $\QFT_{1|1}^n\to \sqm^n$.
\end{hyp}

\begin{proof}[Sketch.]
The argument is analogous to the one for Hypothesis~\ref{hyp:QMandBord}, and so we will be (even more) brief. The functor $\QFT_{1|1}^n\to \sqm^n$ sends a $1|1$-Euclidean field theory $E$ to the generator~$Q$ of its supersemigroup representation $E(\sI^+_{t,\theta})$. Indeed, by \cite[Lemma 5.10]{HST}, we have $E(\sI_{t,\theta}^+)=e^{-tH+\theta Q}$ from a uniquely determined unbounded odd Clifford linear operator $Q$ with domain $E(\spt^+)\subset \cH$. 

To show this assignment determines an epimorphism, a generators and relations presentation of $\EBord_{1|1}$ (compare \cite[Theorem 6.42]{HST}) shows that a $1|1$-Euclidean field theory of degree~$n$ is determined by a self-adjoint Clifford module $(E(\spt^+),\langle-,-\rangle)$ and a self-adjoint $\cCl_n$-linear supersemigroup representation from $E(\sI_{t,\theta}^+)$. A construction analogous to~\eqref{eq:repof1EB0} then determines the desired 1-sided inverse functor. 
\ep

\begin{rmk}\label{rmk:weirdsheaf2}
A hope for an equivalence of stacks $\QFT_{1|1}^n\to \sqm^n$ runs into a problem: there are many non-isomorphic topological vector spaces~$E(\pt^+)$ leading to non-isomorphic field theories that become isomorphic as objects in $\sqm^n$ (after completing a topological inner product space to a Hilbert space). One possible fix is to impose a weaker equivalence relation. \emph{Concordance classes} of field theories~\cite[Definition 1.12]{ST11} are the most promising option, and there is an expectation that concordance classes of field theories over $X$ are in bijection with concordance classes of maps from $X$ to the stack $\QFT_{1|1}^n$. However, a concordance equivalence $\QFT_{1|1}^n\to \sqm^n$ still requires further technical modifications to the definition of field theory, see \cite[Remark 3.16]{ST11} and~\cite{Ulrickson}. Roughly, one must allow families of field theories whose state spaces are vector bundles with nonconstant rank, where changes in rank mimic the jumps in index bundle for a family of Dirac operators. 
\end{rmk}



\subsection{Additional insights from the bordism framework} Extending the ideas in~\S\ref{sec:1gens}, the bordism framework offers the following embellishments of quantum mechanical systems with supersymmetry.

 \begin{ex}[Orientation reversal, KO, and KR]\label{ex:OrReal11}
Following Example~\ref{ex:OrReal1}, orientation reversal on $\EBord_{1|1}(\X)$ affords a $\Z/2$-action on the groupoid of (twisted) field theories over~$\X$ by precomposition. Taking (homotopy) fixed points yields \emph{unoriented $1|1$-dimensional Euclidean field theories} (see \cite[Example 6.16]{HST}). The $\Z/2$-fixed data is again a real structure~\eqref{eq:realstructure} on the space of states, and time-evolution operators are required to be real. Hence, unoriented theories correspond to KO-theory, see \cite[\S6.8]{HST}. A further variation (related to Atiyah's KR-theory) takes as input a stack $\X$ with $\Z/2$-action; this determines a $\Z/2$-action on $\EBord_{1|1}(\X)$ through both orientation reversal of bordisms and the $\Z/2$-action on $\X$, which in turn gives a $\Z/2$-action on field theories. 
 \end{ex}

\begin{ex}[Super vector bundles with super connection]
Dumitrescu~\cite{Florin} and Stoffel~\cite{Stoffeltwists} described generalizations of Examples~\ref{ex:VB1} and~\ref{ex:pVB1} to $1|1$-dimensional Euclidean field theories. Specifically, given a (finite-rank) super vector bundle $V\to X$ over a smooth manifold $X$ with Quillen superconnection~$\A$, Dumitrescu constructs \emph{super parallel transport} along $1|1$-dimensional super paths in $X$. This gives a functor
\beq\label{eq:sPar}
{\rm sPar}\colon \Vect^\A(X)\to \QFT_{1|1}^\one (X).
\eeq
For a gerbe $\alpha$ on $X$, Stoffel provides an $\alpha$-twisted generalization of~\eqref{eq:sPar}.
\end{ex}

\begin{rmk}\label{rmk:Thm11}
A one-sided inverse to~\eqref{eq:sPar} takes values in super connections on vector bundles that are not necessarily finite-rank~\cite{DBEEFT}. This proves part of Conjecture~\ref{thm11} in the framework from~\cite{ST11}: the index bundle of a Clifford linear super connection provides the desired cocycle map from field theories to K-theory. Bismut's super connection~\cite{Bismutindex} for a family of spin$^c$ manifolds furthermore gives a partial construction of the quantization map in~\eqref{eq:Thm11}. To complete the proof of Conjecture~\ref{thm11} in this framework, it would suffice to to construct an infinite-rank generalization of Dumitrescu's super parallel transport~\eqref{eq:sPar} for the Bismut super connection, yielding the ``quantize" arrow in~\eqref{eq:Thm11}. 
\end{rmk}

\begin{rmk}\label{rmk:Chernchar}
As the bordisms in $\EBord_{1|1}(\X)$ are supermanifolds, the moduli spaces involved can be subtle, e.g., see~\cite{WittenRiemann}. For field theories over $\X=\pt$, all the super geometry is essentially funneled into the ``super length" parameter $(t,\theta)\in \R^{1|1}_{\ge 0}$. However, for bordisms over a stack~$\X$ the super aspect becomes more involved. For example, for a manifold~$X$ the value of $E\in \QFT_{1|1}^\one(X)$ on super circles determines a differential form $Z_E\in C^\infty(\Map(\R^{0|1},X))\simeq \Omega^\bullet(X)$, that (in conjunction with~\eqref{eq:sPar}) can be interpreted as the Chern character of the family of field theories~\cite{Han}. Generalizing this, functions on the superstacks $\Map(\R^{0|1},X\nsq G)$ provide a model for equivariant de~Rham cohomology~\cite{Wu,Stolzclass}, leading to equivariant Chern characters of $1|1$-Euclidean field theories~\cite{DBEHan}. We do not know a way of accessing this refined information without the ``super bordism" framework from~\cite{ST11}. 
\end{rmk}

\begin{ex}[Projective $G$-representations] 
Adding supersymmetry to Examples~\ref{ex:rep1} and~\ref{ex:prep1}, the category $\EBord_{1|1}(\pt\nsq G)$ adds $G$-bundles with connection to the bordisms depicted in Figure~\ref{fig:11EBord}. Restricting the supersemigroup of super intervals $\sI_{t,\theta}^+$ to $\theta=0$ recovers a subcategory of $\EBord_{1|1}(\pt\nsq G)$ equivalent to $G$-bundles over ordinary intervals pictured in~\eqref{eq:Gsemigroup}. Hence, one recovers the same information as before, but with extra enhancements. First, the supersymmetry determines a square root $Q$ of the previous Hamiltonian $H$. Second, connections on $G$-bundles over supermanifolds are more intricate, involving the geometry of equivariant differential forms. We refer to~\cite{DBEHan} for further discussion. 
%
\end{ex}

%

\section{Elliptic cohomology and 2-dimensional supersymmetric field theories}\label{sec:4}

Identifying the complete mathematical content of a 2-dimensional supersymmetric quantum field theory remains an open problem. The rough contours come from loop space generalizations of supersymmetric quantum mechanics as described by Witten~\cite{Witten_Dirac} and reviewed in~\S\ref{sec:41}-\ref{sec:42}. However, since the beginning it was known that this simplistic generalization fails to capture additional data and property of a 2-dimensional supersymmetric field theory, e.g., real structures on spaces of states~\cite{Witten_Elliptic} and anomaly cancelation data~\cite[\S3]{Witten_Dirac}. Such data anticipated from physical considerations emerge naturally---and for purely geometric reasons---within Stolz and Teichner's bordism framework~\cite{ST11}. We explain aspects in~\S\ref{sec:21bord}-\ref{sec:21add}. In~\S\ref{sec:final}, we outline expected features of fully-extended supersymmetric field theories and speculate on some inroads to Conjecture~\ref{conj}.

\subsection{The basic data of a 2-dimensional supersymmetric quantum field theory}\label{sec:41}
To first approximation, 2-dimensional quantum field theory is a loop space version of quantum mechanics~\cite{Witten_Dirac}: in addition to a state space and Hamiltonian, we have an $S_r^1$-action on states for $S^1_r=\R/2\pi r \Z$ a circle of radius~$r$. The infinitesimal~$S^1_r$-action defines the \emph{momentum operator} $P^{(r)}\in \End(\cH)$ with spectrum in~$\frac{1}{r}+\frac{n}{24r}\cdot \Z\subset \R$.\footnote{The shift by $n/24$ for $n\in \Z$ encodes the \emph{degree} of the theory related to an anomaly described below.}
For each $r\in \R_{>0}$, the Hamiltonian $H^{(r)}$ and momentum operator $P^{(r)}$ define a semigroup representation 
\vspace{-.1in}
\beq\label{eq:QFTsemigroup}
&&\begin{array}{c}\displaystyle e^{-t H^{(r)}+i xP^{(r)}}=q^{ H^{(r)}_+}\bar q^{ H^{(r)}_-}\colon \mathbb{H}\to \End(\cH),\\ \null \\ \mathbb{H}=\{\tau\in \C \mid {\rm im}(\tau)>0\}  \\ \tau=\frac{x+it}{2\pi r}, \ q=e^{2\pi i \tau}\end{array}\qquad \qquad \qquad 
\begin{tikzpicture}[baseline=(basepoint)];
\node (D) at (-.5,0) {$\bullet$};
\node (G) at (.3,.8) {$\bullet$};
\draw[thick] (.5,0) to (.5,1);
\draw[thick] (-.5,0) to (-.5,1);
\draw[thick,<->] (-1,0) to node [left] {$t$} (-1,1);
\draw[thick,->, bend left=30] (-.5,1.35) to node [above] {$x$} (.4,1.36);
\node [draw, thick, ellipse, minimum width=1cm,minimum height=.45cm] (b) at (0,0){};
\node [draw, thick, ellipse, minimum width=1cm,minimum height=.45cm] (c) at (0,1){};
\node (G) at (1.5,.5) {$\rC_{r,\tau}$};
\path (0,.5) coordinate (basepoint);
\end{tikzpicture}
\eeq
generalizing the time-evolution semigroup~\eqref{eq:semigroup} in quantum mechanics, where (for later convenience) we introduce operators~\cite[Eq.~106]{WittenICM} 
\beq\label{eq:L0ops}
&&H_+^{(r)}=\frac{r}{2}(H^{(r)}+P^{(r)}),\qquad  H_-^{(r)}=\frac{r}{2}(H^{(r)}-P^{(r)}).
\eeq
As indicated on the right of~\eqref{eq:QFTsemigroup}, one can reformulate this semigroup representation geometrically: the parameters $(r,\tau)$ determine a cylinder $\rC_{r,\tau}$ viewed as a Euclidean (i.e., flat Riemannian) bordism between circles of radius~$r$, where the complex parameter $\tau$ encodes the height~$t$ of the cylinder and the rotation $x \ {\rm mod} \ 2\pi r$ that measures the difference in parameterizations of the boundary circles~\cite[Figure~2]{WittenICM}. The \emph{partition function} is the character (i.e., trace) of the semigroup representation~\eqref{eq:QFTsemigroup},
\beq\nonumber
&&Z(r,\tau):=\Tr(q^{ H_+^{(r)}}\bar q^{ H_-^{(r)}})\in C^\infty(\R_{>0}\times \mathbb{H}),\quad r\in \R_{>0}, \ \tau\in \HH. 
\eeq
 Path integral arguments~\cite[\S3]{WittenICM} identify $Z(r,\tau)$ with the value of the theory on the torus $\rT_{r,\tau}=\R^2/(2\pi r\Z\oplus 2\pi\tau\Z)$ gotten from gluing the ends of a cylinder together. General invariance properties lead one to expect that this value only depends on the underlying torus, not the based lattice generated by~$r$ and~$\tau$. The group $\SL_2(\Z)$ acts by changing the basis of this lattice, 
\beq\nonumber
\left[\begin{array}{cc} a & b \\ c& d\end{array}\right]\cdot (r,\tau)=\left(r|c\tau+d|,\frac{a\tau+b}{c\tau+d}\right),\quad \left[\begin{array}{cc} a & b \\ c& d\end{array}\right]\in \SL_2(\Z)
\eeq 
and so naively one might expect $Z(r,\tau)$ to be an $\SL_2(\Z)$-invariant function. However, tracking the shift by $n/24r$ in the spectrum of~$P^{(r)}$, one finds $Z(r,\tau)$ is only invariant up to a root of unity  encoded by a character $\chi^n\colon \MP_2(\Z)\to \C^\times$ \cite{Witten_Elliptic}
\beq\label{eq:SL2invariance}
&&Z(r,\tau)\in C^\infty(\R_{>0}\times \mathbb{H};\chi^n)^{\MP_2(\Z)}
\eeq
for the metaplectic double cover $\MP_2(\Z)\to \SL_2(\Z)$. This failure of~$Z(r,\tau)$ to be $\SL_2(\Z)$-invariant is a hallmark of an anomalous theory~\cite[\S4]{Freed_Det}. Indeed, $Z(r,\tau)$ determines a section of the $n$th tensor power of the Pfaffian line bundle over Euclidean tori whose sections are closely related to modular forms of weight $n/2$; see Remark~\ref{rmk:etaMF}.

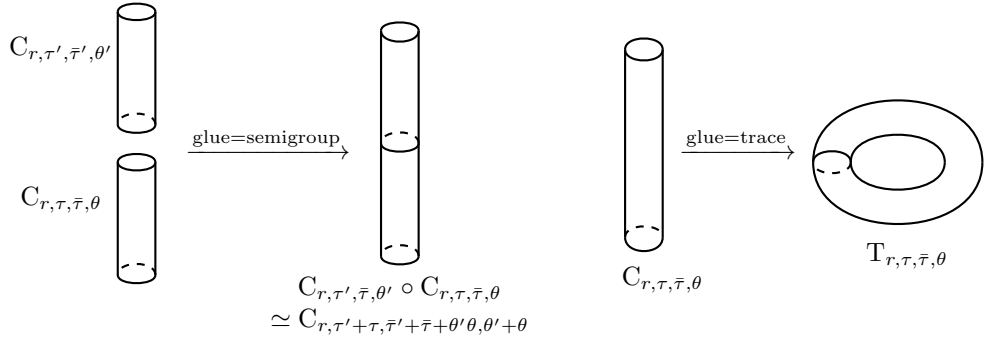
\begin{figure}
\beq\label{eq:Euclideanpicture}\label{eq:gluecylinderpic}\nonumber
&&\begin{tikzpicture}[baseline=(basepoint)];
	\begin{pgfonlayer}{nodelayer}
		\node [style=none] (0) at (1.5, 1.75) {};
		\node [style=none] (1) at (1, 1.75) {};
		\node [style=none] (2) at (1.5, 0.25) {};
		\node [style=none] (3) at (1, 0.25) {};
		\node [style=none] (4) at (1.5, -1.25) {};
		\node [style=none] (5) at (1, -1.25) {};
		\node [style=none] (6) at (-2, 2) {};
		\node [style=none] (7) at (-2.5, 2) {};
		\node [style=none] (8) at (-2, 0.5) {};
		\node [style=none] (9) at (-2.5, 0.5) {};
		\node [style=none] (10) at (-2, 0) {};
		\node [style=none] (11) at (-2.5, 0) {};
		\node [style=none] (12) at (-2, -1.5) {};
		\node [style=none] (13) at (-2.5, -1.5) {};
		\node [style=none] (14) at (-2.5, 0.5) {};
		\node [style=none] (15) at (-2.5, 0.5) {};
				\node [style=none] (21) at (-.5, .25) {$\xrightarrow{\rm glue=semigroup}$};
		\node [style=none] (20) at (1.25,-1.9) {$\begin{array}{c} \rC_{r,\tau',\bar\tau,\theta'}\circ \rC_{r,\tau,\bar\tau,\theta}\\ \simeq \rC_{r,\tau'+\tau,\bar\tau'+\bar\tau+\theta'\theta,\theta'+\theta}\end{array}$};
\node [style=none] (19) at (-3.25,-.5) {$\rC_{r,\tau,\bar\tau,\theta}$};
\node [style=none] (18) at (-3.25,1.5) {$\rC_{r,\tau',\bar\tau',\theta'}$};
	\end{pgfonlayer}
	\begin{pgfonlayer}{edgelayer}
		\draw [bend left=270, looseness=0.75,thick] (0.center) to (1.center);
		\draw [bend left=270,thick,dashed] (2.center) to (3.center);
		\draw [bend left=270,thick,dashed] (4.center) to (5.center);
		\draw [thick] (0.center) to (4.center);
		\draw [thick] (1.center) to (5.center);
		\draw [bend left=90,thick] (0.center) to (1.center);
		\draw [bend left=90, looseness=0.75,thick] (2.center) to (3.center);
		\draw [bend left=90, looseness=0.75,thick] (4.center) to (5.center);
		\draw [bend left=270, looseness=0.75,thick] (6.center) to (7.center);
		\draw [bend left=270,thick,dashed] (8.center) to (9.center);
		\draw [bend left=90, looseness=0.75,thick] (6.center) to (7.center);
		\draw [bend left=90, looseness=0.75,thick] (8.center) to (9.center);
		\draw [thick] (8.center) to (6.center);
		\draw [thick] (9.center) to (7.center);
		\draw [bend left=270, looseness=0.75,thick] (10.center) to (11.center);
		\draw [bend left=270,thick,dashed] (12.center) to (13.center);
		\draw [bend left=90, looseness=0.75,thick] (10.center) to (11.center);
		\draw [bend left=90, looseness=0.75,thick] (12.center) to (13.center);
		\draw [thick] (12.center) to (10.center);
		\draw [thick] (13.center) to (11.center);
	\end{pgfonlayer}
	\path (0,0) coordinate (basepoint);
\end{tikzpicture}
\qquad\quad\begin{tikzpicture}[baseline=(basepoint)];
\node [style=none] (19) at (-2,-1.57) {$\rC_{r,\tau,\bar\tau,\theta}$};
\node [style=none] (19) at (1.25,-1.25) {$\rT_{r,\tau,\bar\tau,\theta}$};
	\begin{pgfonlayer}{nodelayer}
		\node [style=none] (1) at (2.25, 0) {};
		\node [style=none] (2) at (1.75, 0) {};
		\node [style=none] (4) at (0.5, 0) {};
		\node [style=none] (5) at (0, 0) {};
		\node [style=none] (8) at (-2.5, 1.5) {};
		\node [style=none] (9) at (-2, 1.5) {};
		\node [style=none] (12) at (-2.5, -1) {};
		\node [style=none] (13) at (-2, -1) {};
						\node [style=none] (21) at (-1, .25) {$\xrightarrow{\rm glue=trace}$};
	\end{pgfonlayer}
	\begin{pgfonlayer}{edgelayer}
		\draw [bend left=270, looseness=1.25,thick] (1.center) to (5.center);
		\draw [bend right=90,thick] (2.center) to (4.center);
		\draw [bend left=90,thick,dashed] (4.center) to (5.center);
		\draw [bend right=90,thick] (4.center) to (5.center);
		\draw [bend right=90, looseness=1.25,thick] (5.center) to (1.center);
		\draw [bend right=90,thick] (4.center) to (2.center);
		\draw [bend left=90,thick] (8.center) to (9.center);
		\draw [bend right=90,thick] (8.center) to (9.center);
		\draw [bend right=90, looseness=1.25,thick] (12.center) to (13.center);
		\draw [bend left=90,thick,dashed] (12.center) to (13.center);
		\draw [thick] (9.center) to (13.center);
		\draw [thick] (8.center) to (12.center);
	\end{pgfonlayer}
	\path (0,0) coordinate (basepoint);
\end{tikzpicture}
\eeq
\caption{Composition in the semigroup corresponds to the gluing of cylinders, and the value on the torus corresponds to a trace.}
\label{fig:2d}
\end{figure}


Next we add supersymmetry:\footnote{This discussion only considers the Ramond sector; in~\S\ref{sec:21bord} we incorporate the Neveu--Schwarz sector.} the data of $\mathcal{N}_+$ odd commuting square roots of $ H_+^{(r)}$ and~$\mathcal{N}_-$ odd commuting square roots of $ H_-^{(r)}$ is called \emph{$\mathcal{N}=(\mathcal{N}_+,\mathcal{N}_-)$ supersymmetry}. We will be interested in $\mathcal{N}=(0,1)$ supersymmetry, i.e., a single odd operator $Q^{(r)}$ with
\beq\label{Eq:N01susy}
(Q^{(r)})^2=\frac{1}{2}[Q^{(r)},Q^{(r)}]= H_-^{(r)},\qquad [Q^{(r)},P^{(r)}]=[Q^{(r)},H^{(r)}]=0. 
\eeq
Generalizing time evolution~\eqref{eq:N1super} in supersymmetric quantum mechanics, the operator $Q^{(r)}$ promotes the semigroup representation~\eqref{eq:QFTsemigroup} to a \emph{supersemigroup} representation 
\beq\label{Eq:ogsupersemi}
q^{ H_+^{(r)}}\bar{q}^{ H_-^{(r)}+\sqrt{2\pi i}\theta Q^{(r)}},\qquad \mathbb{H}^{2|1}\to \End(\mathcal{H})
\eeq
where $\mathbb{H}^{2|1}\subset \R^{2|1}$ is the upper half-space of the super Lie group $\R^{2|1}$ with multiplication
\beq\label{Eq;21group}
&&(\tau',\bar\tau',\theta')\cdot (\tau,\bar\tau,\theta)=(\tau'+\tau,\bar \tau'+\bar \tau+\theta'\theta,\theta'+\theta), \qquad (\tau',\bar\tau',\theta'),(\tau,\bar\tau,\theta)\in \R^{2|1}.
\eeq
This supersemigroup representation can again be interpreted geometrically, but now as a representation of a category of \emph{super cylinders} whose geometric parameters are $(r,\tau,\bar\tau,\theta)\in \R_{>0}\times \HH^{2|1}$. Similarly, the trace of this semigroup representation corresponds to the value of the theory on a \emph{super torus}, which imposes an expected invariance property analogous to~\eqref{eq:SL2invariance}. This geometry is pictured schematically in Figure~\ref{fig:2d}. 
\begin{rmk}\label{rmk:etaMF}
Consider the quotient stack $\HH\sq \MP_2(\Z)$ classifying spin elliptic curves with odd spin structure, where $\MP_2(\Z)$ acts by fractional linear transformations on $\HH$ through the double cover $\MP_2(\Z)\to \SL_2(\Z)$. This stack supports a line bundle determined by the character
\beq\label{eq:abelianization}
\chi\colon \MP_2(\Z)\to \MP_2(\Z)/[\MP_2(\Z),\MP_2(\Z)]\simeq \Z/24\hookrightarrow \C^\times
\eeq
which is the \emph{Pfaffian line bundle} $\Pf\to \HH\sq \MP_2(\Z)$~\cite[\S4]{Freed_Det}. A second line bundle has sections modular forms of weight 1/2, i.e., holomorphic functions $f(\tau)$ satisfying 
\beq\label{eq:mfhalf}
f\left(\frac{a\tau+b}{c\tau+d}\right)=\sqrt{c\tau+d}f(\tau). 
\eeq
The line determined by~\eqref{eq:mfhalf} and the Pfaffian line are generators of the Picard group $\Pic(\HH\sq \MP_2(\Z))\simeq \Z/24$, and therefore are isomorphic. The Dedekind $\eta$-function furnishes an explicit isomorphism, and in particular a bijection on global sections
$$
\Gamma(\HH\sq \MP_2(\Z);\Pf^{\otimes n})\xrightarrow{\sim} \MF_{n/2},\qquad Z(q)\mapsto \eta(q)^n Z(q),\qquad \eta(q)=q^{1/24}\prod_{k>0} (1-q^k)
$$ 
as can be checked from the standard transformation properties of $\eta(q)$. This bijection between modular forms and functions on $\HH$ transforming with a phase under $\MP_2(\Z)$ is used extensively in the physics literature, e.g., \cite[Eq.~18]{Witten_Dirac}. 
\end{rmk}

\begin{ex}[Chiral theories]\label{ex:chiral}
The data~\eqref{eq:L0ops} determine a \emph{chiral} field theory\footnote{The more common usage of \emph{chiral} assumes that $H_+^{(r)}=L_0$ and $H_-^{(r)}=\bar L_0$ are part of the data of a conformal field theory with Virasoro generators $L_k$ and $\overline{L}_k=0$ for $k\in \Z$.} if $H_-^{(r)}=0$. A chiral theory has a canonical supersymmetric extension by taking $Q^{(r)}=0$, and so (following Example~\ref{ex:TQFTQ}) chiral theories are the relevant 2-dimensional generalization of topological theories; see~\eqref{eq:chiralinclusion} for further discussion. For now we observe that any chiral conformal field theory (e.g., the data of a vertex operator algebra) gives a 2-dimensional chiral theory and hence a supersymmetric theory. 
\end{ex}


\begin{ex}[The supersymmetric sigma model and the Witten genus]\label{ex:Wittengenus}
Let $M$ be a Riemannian spin manifold. Witten's spinor bundle on $LM$ restricted along $M\hookrightarrow LM$ is (by definition) the formal power series of vector bundles~\cite[Eq.~24]{Witten_Dirac}
\beq\label{eq:Wittenpower}
&&\bS \otimes \left(\bigotimes_{m\ge 1} \Sym_{q^m}(TM)\right)=\bS+q\bS\otimes TM+q^2\bS\otimes (TM+\Sym^2(TM))+\dots
\eeq
where $\bS$ is the spinor bundle on $M$ and $\Sym_q(TM)=1+qTM+q^2\Sym^2(TM)+\dots $. Let $W_k(TM)$ denote the coefficient of $q^k$ in~\eqref{eq:Wittenpower}, and $\cH$ denote the Hilbert completion of the sections of~\eqref{eq:Wittenpower}. For each $r\in \R_{>0}$ define operators $P^{(r)}$ and $Q^{(r)}$ on $\cH$ via
$$
Q^{(r)}|_{W_k(TM)}=\slashed{D}\otimes W_k(TM)+\sqrt{k}/\sqrt{r},\qquad P^{(r)}|_{W_k(TM)}=k/r-n/24r,\qquad n=\dim(M)
$$
where $\slashed{D}$ is the $\cCl_{-n}$-linear Dirac operator on $M$. This gives the input data~\eqref{Eq:N01susy} for a supersymmetric field theory. The Atiyah--Singer index theorem computes the (Clifford) super trace of the associated semigroup representation (see Remark~\ref{rmk:Cliffsuper})
\beq\label{eq:sigmamodelZ}
Z(r,\tau)=\Tr_{\cCl_{-n}}(q^{ H^{(r)}_+}\bar q^{ H^{(r)}_-})=q^{-n/24}\sum_{k=0}^\infty q^k \hat{A}(M,W_n(TM))
\eeq
where $\hat{A}(M,W_n(TM))$ is the twisted $\hat{A}$-class. Zagier~\cite{Zagiermodular} verified that the normalized partition function $\eta(q)^nq^{-n/24}\sum_{k=0}^\infty q^k \hat{A}(M,W_n(TM))$ is a modular form of weight~$n/2$ when the 1st rational Pontryagin class $p_1(TM)\in \H^4(M;\Q)$ vanishes. Remark~\ref{rmk:etaMF} then confirms that the action of $\MP_2(\Z)$ on $Z(r,\tau)$ is through a 24th root of unity. 
\end{ex}

\begin{rmk}\label{rmk:Cliffsuper}
The \emph{Clifford super trace} of a $\cCl_n$-linear endomorphism $T$ on a Clifford module~$V$ is defined as~\cite[Definition~3.2.16]{ST04}
$$
\sTr_{\cCl_n}(T)=\sTr(\Gamma_n\circ T), \qquad T\colon V\to V
$$
where $\Gamma_n=(2i)^{-n/2}e_1\cdots e_n\in \cCl_n$ for $\{e_1,\dots, e_n\}$ the standard basis of $\R^n$, identified with (odd) elements of the Clifford algebra under the inclusion $\R^n \subset \cCl_n$. 
\end{rmk}

\begin{rmk}
There are several ongoing efforts using analytical tools to rigorously construct a spinor bundle on loop space and its attendant Dirac operator in a way that recovers Example~\ref{ex:Wittengenus} upon restriction to the constant loops, e.g., see \cite{ST_spinors,KottkeMelrose1,KottkeMelrose2,KLW,LudewigWaldorf2grp,Ludewigspinor}. More algebraic approaches have constructed closely related objects in a formal neighborhood of the constant loops, e.g.,~\cite{CDO1,CDO2,costello_WG1,costello_WG2}.

\end{rmk}

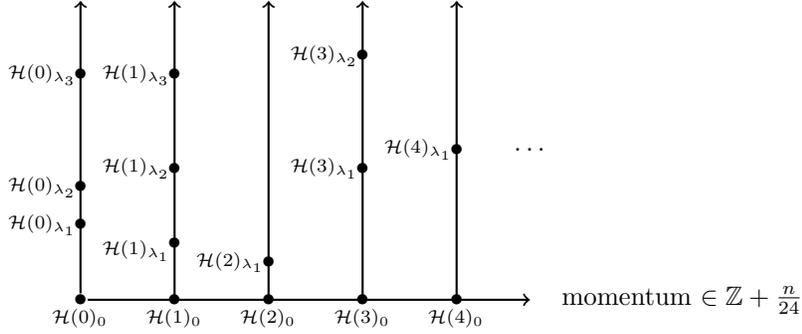
\begin{figure}
\beq\nonumber
\begin{tikzpicture}

		\node [style=none] (Q) at (6, 2) {$\dots$};
		
		\node [style=none] (A) at (0, 0) {$\bullet$};
		\node [style=none] (B) at (6, 0) {};
		\node [style=none] (C) at (0, 4) {};
				\node [style=none] (1) at (0, 1) {$\bullet$};
				\node [style=none] (1) at (-.5, 1) {$_{\cH(0)_{\lambda_1}}$};
				\node [style=none] (2) at (0, 1.5) {$\bullet$};
				\node [style=none] (2) at (-.5, 1.5) {$_{\cH(0)_{\lambda_2}}$};
				\node [style=none] (3) at (0, 3) {$\bullet$};
				\node [style=none] (3) at (-.5, 3) {$_{\cH(0)_{\lambda_3}}$};

		\node [style=none] (AA) at (1.25, 0) {$\bullet$};										
		\node [style=none] (E) at (1.25, 0) {};
		\node [style=none] (F) at (1.25, 4) {};

			\node [style=none] (1) at (1.25, 0.75) {$\bullet$};
						\node [style=none] (1) at (.75, 0.65) {$_{\cH(1)_{\lambda_1}}$};
			\node [style=none] (2) at (1.25, 1.75) {$\bullet$};
						\node [style=none] (2) at (.75, 1.75) {$_{\cH(1)_{\lambda_2}}$};
			\node [style=none] (3) at (1.25, 3) {$\bullet$};
				\node [style=none] (3) at (.75, 3) {$_{\cH(1)_{\lambda_3}}$};

		\node [style=none] (AAA) at (2.5, 0) {$\bullet$};										
		\node [style=none] (H) at (2.5, 4) {};
		\node [style=none] (G) at (2.5, 0) {};
				\node [style=none] (1) at (2.5, 0.5) {$\bullet$};
								\node [style=none] (1) at (2, 0.5)  {$_{\cH(2)_{\lambda_1}}$};

		\node [style=none] (AAA) at (3.75, 0) {$\bullet$};												
		\node [style=none] (I) at (3.75, 0) {};		
		\node [style=none] (J) at (3.75, 4) {};
			\node [style=none] (1) at (3.75, 1.75) {$\bullet$};
						\node [style=none] (1) at (3.25, 1.75) {$_{\cH(3)_{\lambda_1}}$};
			\node [style=none] (2) at (3.75, 3.25) {$\bullet$};
						\node [style=none] (2) at (3.25, 3.25) {$_{\cH(3)_{\lambda_2}}$};
		
		\node [style=none] (AAAA) at (5, 0) {$\bullet$};														
		\node [style=none] (K) at (5, 0) {};
		\node [style=none] (L) at (5, 4) {};
		
		\node [style=none] (22) at (5, 2) {$\bullet$};		
		\node [style=none] (22) at (4.5, 2) {$_{\cH(4)_{\lambda_1}}$};

		\node [style=none] (31) at (8, 0) {momentum $\in \Z+\frac{n}{24}$};
		\node [style=none] (26) at (5, -0.25) {$_{\mathcal{H}(4)_{0}}$};
		\node [style=none] (27) at (3.75, -0.25) {$_{\mathcal{H}(3)_{0}}$};
		\node [style=none] (28) at (2.5, -0.25) {$_{\mathcal{H}(2)_{0}}$};
		\node [style=none] (29) at (1.25, -0.25) {$_{\mathcal{H}(1)_{0}}$};
		\node [style=none] (30) at (0, -0.25) {$_{\mathcal{H}(0)_{0}}$};
		\draw[thick,->] (A) to (B);
		\draw[thick,->] (A) to (C);
		\draw[thick,->] (E) to (F);
		\draw[thick,->] (G) to (H);
		\draw[thick,->] (I) to (J);
		\draw[thick,->] (K) to (L);
\end{tikzpicture}
\eeq
\caption{An object in $\qft^n$ determines a continuous map from $\R_{>0}$ to configurations in $\R^2$ where labels are subspaces~$\cH(k)_\lambda\subset \cH_n$ indexed by $k\in \Z$ and $\lambda\in \R$, where $P^{(r)}$ and $H^{(r)}$ share an eigenspace with eigenvalues $k/r+n/24$ and $\lambda$, respectively. It is convenient to arrange these eigenspaces in the plane with horizontal axis eigenvalues of $rP^{(r)}$ and the vertical axis by eigenvalues of $H_-^{(r)}$ (i.e., the difference between energy and momentum). 
Supersymmetry then requires the eigenvalues in the vertical direction be nonnegative and decorates each label with an odd operator gotten from restricting $Q^{(r)}$ to $\cH(k)_\lambda$. This gives a copy of Figure~\ref{fig2} for each $k\in \Z$.
}
\label{fig3}
\end{figure}

\subsection{Approximating categories of 2-dimensional field theories}\label{sec:42}

Next we formulate preliminary (and incomplete) definitions of 2-dimensional quantum field theories that collect the data and properties of the operators~\eqref{eq:L0ops} and~\eqref{Eq:N01susy} described above. 
With Example~\ref{ex:Wittengenus} in mind, we add a Clifford degree to the definitions below: for each $n\in \Z$ fix a graded $\cCl_n$-module $\cH_n$, e.g., for an infinite dimensional separable Hilbert space~$\cH$, set~$\cH_n=\cCl_n\otimes \cH$.

\begin{defn}\label{defn:QFTsystem}
Define a topological groupoid $\qft^n$ whose objects are continuous $\R_{>0}$-families of $\cCl_n$-linear unbounded even operators $\{H_+^{(r)},H_-^{(r)}\}_{r\in \R_{>0}}$ on $\cH_n$ satisfying:
\begin{itemize}
\item[(i)] $H_+^{(r)}$ and $H_-^{(r)}$ determine self-adjoint trace class semigroup representations~\eqref{eq:QFTsemigroup};
\item[(ii)] the Clifford super trace
\beq\label{eq:anomalous}
Z(r,\tau)&:=&\sTr_{\cCl_n}(q^{ H_+^{(r)}}\bar q^{ H_-^{(r)}}),\quad Z\left(r|c\tau+d|,\frac{a\tau+b}{c\tau+d}\right)=\chi^n(\gamma)Z(r,\tau),\\
&&\gamma=\left(\left[\begin{array}{cc} a & b \\ c& d\end{array}\right],(c\tau+d)^{1/2}\right)\in \MP_2(\Z)\nonumber
\eeq
 is $\MP_2(\Z)$-invariant up to a phase specified by the character~$\chi^n$ from~\eqref{eq:abelianization}. 
 \end{itemize}
 Morphisms in $\qft^n$ are continuous maps $\R_{>0}\to \PU^\pm (\cH_n)$, acting on $H_\pm^{(r)}$ by conjugation. 
\end{defn}

The operator $rP=H_+^{(r)}+H_-^{(r)}$ has spectrum in the discrete set $\Z+\frac{n}{24}$, and therefore this spectrum does not depend on the continuous parameter~$r$. Up to isomorphism, this gives a decomposition of the state space
\beq\label{eq:weighspaces}
\cH_n\simeq \bigoplus_{k\in \Z} \cH(k),\qquad H_-^{(r)}(k):=H_-^{(r)}|_{\cH(k)}
\eeq 
for which $rP^{(r)}$ acts on $\cH(k)$ by $k+\frac{n}{24}\in \Z$, and~$H_-^{(r)}(k)$ in~\eqref{eq:weighspaces} is defined by restriction. In this way, an object in $\qft^n$ can be viewed as an infinite collection of quantum mechanical systems for the Hamiltonians $H_-^{(r)}(k)$, one for each momentum number $k\in \Z$ and radius $r\in \R_{>0}$, where the dependence on the radius parameter is continuous.\footnote{More precisely, $H_-^{(r)}(k)$ is a combination of a Hamiltonian and a generator of a projective $S^1$-symmetry for the associated quantum mechanical system, see~\cite[\S1.3]{Witteninstrings}.} This intuition can be made precise through a configuration space description of the operators $H_+^{(r)}$ and $H_-^{(r)}$ pictured in Figure~\ref{fig3}: each vertical line is a configuration defining a quantum mechanical system, see also \cite[Figure~4]{PokmanPhD}.

\begin{prop}
For each $n\in \Z$, there is a homotopy equivalence $|\qft^n|\simeq B\PU^\pm(\cH_n)$.
\end{prop}
\begin{proof}[Sketch.]
As described above, for each weight space $\cH(k)\subset \cH_n$ we obtain a continuous~$\R_{>0}$-family of Hamiltonians~$H(k)^{(r)}$. For each such Hamiltonian $H^{(r)}_k$ there is a deformation sending all eigenspaces to infinity. This determines a contracting homotopy to the basepoint with domain~$\{0\}\subset \cH_n$. Applying these homotopies for each $k\in \Z$ simultaneously, the space of objects of $\qft^n$ has a deformation retract to the zero theory corresponding to the empty configuration in Figure~\ref{fig3}. By the same argument as in Proposition~\ref{prop:KZ3}, we obtain the claimed homotopy type. 
\ep

\begin{rmk}
The space $B\PU^\pm(\cH_n)$ has a pair of nontrivial homotopy groups, $\pi_1B\PU^\pm(\cH_n)\simeq \Z/2$ and $\pi_3B\PU^\pm(\cH_n)\simeq \Z$. These are connected by a nontrivial $k$-invariant. 
\end{rmk}

As in the case of quantum mechanics, the homotopy type of $\qft^n$ becomes more interesting if we add supersymmetry. 

\begin{defn}\label{defn:SQFTsystem}

Define a topological groupoid $\sqft^n$ whose objects are continuous $\R_{>0}$-families of $\cCl_n$-linear unbounded operators $\{H_+^{(r)},Q^{(r)}\}_{r\in \R_{>0}}$, where $H_+^{(r)}$ is even and $Q^{(r)}$ is odd. These data satisfy two properties.
\begin{itemize}
\item[(i)] $H_+^{(r)}$ and $Q^{(r)}$ determine self-adjoint trace class supersemigroup representations~\eqref{Eq:ogsupersemi}.
\item[(ii)] Setting $H_-^{(r)}=(Q^{(r)})^2$, the Clifford super trace~\eqref{eq:anomalous} is $\MP_2(\Z)$-invariant up to a phase specified by the character~$\chi^n$ from~\eqref{eq:abelianization}. 
 \end{itemize}
 Morphisms in $\sqft^n$ are continuous maps $\R_{>0}\to \PU^\pm (\cH_n)$, acting on $H_+^{(r)}$ and $Q^{(r)}$ by conjugation. 
\end{defn}

Following the same discussion as the one after Definition~\ref{defn:QFTsystem}, an object in $\sqft^n$ can be viewed as an infinite collection of supersymmetric quantum mechanical systems for the supersymmetry generators $Q^{(r)}(k):=Q^{(r)}|_{\cH(k)}$, one for each momentum number $k\in \Z$ and radius $r\in \R_{>0}$. We again refer to Figure~\ref{fig3} for a configuration space perspective.

\begin{thm} [{\cite[Theorem~1.0.2]{ST04} and \cite[Theorem 2.2.2]{PokmanPhD}}]\label{thm:Ktate}
The space of objects in $\sqft^n$ maps to $\K^n(\!(q)\!)$, the $n$th space in the spectrum representing elliptic cohomology at the Tate curve. 
\end{thm}
\begin{proof}[Sketch.]
Decomposing $\cH\simeq \bigoplus \cH(k)$ into weight spaces for $rP^{(r)}$, we get an $\R_{>0}$-family of odd, $\cCl_n$-linear operators $\{Q^{(r)}(k)\}_{k\in \Z}$. Choosing a fixed radius (e.g., $r=1$) for each $k$ this determines a map into the space representing $\K^n$ from Theorem~\ref{thm:HamK}. The trace class condition on the semigroup demands that $Q^{(r)}(k)$ represents the zero class for $k$ sufficiently negative, and hence all together we get a map to the space~$\K^n(\!(q)\!)$.
\ep

\begin{thm}[{\cite[page~46]{ST04}}]\label{thm:modularity}
Up to a normalization by $\eta(q)^n$, the partition function defines a continuous map from the space of objects of $\sqft^n$ to (the discrete set of) weight $-n/2$ integral modular forms. 
\end{thm}

\bp
The McKean--Singer argument that results in $t$-independence in~\eqref{eq:partition function2} results in $\bar q$-independence for the $\eta$-function normalized partition function
\beq
\frac{Z(r,\tau)}{\eta(q)^n}&=&\eta(q)^{-n}\Tr_A((-1)^Fq^{ H_+^{(r)}}\bar q^{ H_-^{(r)}})=\eta(q)^{-n}\sTr_A(q^{ H_+^{(r)}}\bar q^{ H_-^{(r)}})\nonumber \\
&=&\eta(q)^{-n}\sTr_A(q^{ H_+^{(r)}}|_{{\rm ker}( H_-^{(r)})})=\eta(q)^{-n}\sum_{k\in \Z} \sTr_A(q^{ H_+^{(r)}}|_{{\rm ker}( Q^{(r)})\bigcap \mathcal{H}(k)})\nonumber\\
&=&\eta(q)^{-n}q^{n/24}\sum_{k\in \Z} q^k{\rm sdim}_{\Fer_n}({\rm ker}(Q^{(r)})\bigcap \mathcal{H}(k))\nonumber\\
&=&\left(\prod_{k>0}(1-q^k)\right)^n\sum_{k\in \Z} q^k{\rm Index}(Q^{(r)}|_{\cH(k)})\in \Z(\!(q)\!),\nonumber 
\eeq
where we have used the decomposition~\eqref{eq:weighspaces}. In the last line above, all powers of $q$ and coefficients of $q^k$ are integers, and the trace class condition on the semigroup demands that the coefficients of $q^k$ vanish for sufficiently negative powers of~$q$. By continuity, the integral power series $\eta(q)^{-n} Z(r,\tau)$ must also be independent of the continuous parameter~$r\in \R_{>0}$. Hence we conclude that $\frac{Z(r,\tau)}{\eta(q)^n}\in \Z(\!(q)\!)$. 

Next we consider the $\MP_2(\Z)$-equivariant property. The equivariance property in Definition~\ref{defn:SQFTsystem} determines the modular transformation property (see Remark~\ref{rmk:etaMF})
$$
 \eta\left(\frac{a\tau+b}{c \tau +d}\right)^{-n}Z\left(r|c\tau+d|,\frac{a\tau+b}{c\tau+d}\right)=(c\tau+d)^{-n/2} \eta(\tau)^{-n}Z(r,\tau)
$$
yielding the claimed weight $-n/2$ integral modular form. 
\ep

\begin{rmk}
Above we take a modular form to be a \emph{weakly holomorphic modular form}, i.e., a holomorphic function $f\in \mathcal{O}(\HH)$ that transforms with weight under the $\SL_2(\Z)$ action and with the property that there is at worst a finite order pole as $\tau\to i\infty$. 
\end{rmk}

Let $\chi\qft^n\subset \qft^n$ denote the subgroupoid of chiral field theories from Example~\ref{ex:chiral}, i.e., $(H_+^{(r)},H_-^{(r)})=(H_+^{(r)},0)$.

\begin{cor}[{\cite[Theorem 1.15]{ST11}}]\label{cor:chiral}
Up to a normalization by $\eta(q)^n$, the partition function defines a continuous map from $\chi\qft^n$ to the set of integral modular forms. Furthermore, every integral modular form arises as the partition function of an object in $\chi\qft^n$ and (in particular) as the partition function of an object in~$\sqft^n$. 
\end{cor}

\bp
The first statement follows from applying Example~\ref{ex:chiral} to Proposition~\ref{thm:modularity}. The last statement follows from taking the canonical supersymmetric theory associated to a chiral one (see Example~\ref{ex:chiral}), provided that every integral modular form arises as the partition function of a chiral theory. To verify this remaining claim, let $f(q)\in \MF^\Z$ be a weight $n/2$ modular form. If $n$ is odd then $f(q)=0$ and we may take the zero theory as our chiral theory. If $n$ is even, then consider $f(q)\eta^n=q^{n/24}\sum a_n q^n\in q^{n/24}\Z(\!(q)\!)$. Define 
\beq\label{eq:chiralconstruction}
\cH(k)= (_{\cCl_2}\C^{1|1})^{\otimes n/2}\otimes \C^{\times a_n} ,\qquad H_+^{(r)}=\bigoplus_{k=-1}^\infty \frac{1}{r}\left(k +\frac{n}{24}\right)\cdot \id _{\cH(k)}
\eeq
where $_{\cCl_2}(\C^{1|1})_\C$ is the bimodule implementing the Morita equivalence between $\cCl_2$ and $\C$. This bimodule has Clifford super trace equal to $1$, so (using multiplicativity for traces of tensor product) the theory~\eqref{eq:chiralconstruction} has partition function
\beq
Z(r,\tau)&:=&\Tr_{\cCl_n}(q^{ H_+^{(r)}}\bar q^{ H_-^{(r)}})=\Tr_{\cCl_n}(q^{ H_+^{(r)}})\nonumber\\
&=&\sum_{k\in \Z}^\infty \Tr_{\cCl_n}(_{\cCl_2}\C^{1|1})^{\otimes n/2})\cdot \Tr(q^{k+n/24}\cdot \id_{\C^{\times a_n}})=q^{n/24}\sum_{k\in \Z}^\infty a_k q^k,\nonumber
\eeq
which transforms by an $n$th root of unity under the $\MP_2(\Z)$-action using the modularity properties of $f(q)$ and Remark~\ref{rmk:etaMF}. 
\ep

\begin{rmk}
Corollary~\ref{cor:chiral} raises a red flag: the edge homomorphism 
\beq\label{eq:edgehomo}
\pi_*\TMF\to \MF^\Z
\eeq
from the coefficients of $\TMF$ to integral modular forms is \emph{not} surjective~\cite[Proposition 4.6]{HopkinsICM2002}. Therefore, if Conjecture~\ref{conj} is to hold, not all integral modular forms should be realized as partition functions of supersymmetric field theories. In particular, Conjecture~\ref{conj} is false if one uses Definition~\ref{defn:SQFTsystem}. It bears noting that the most well-established examples of supersymmetric field theories in physics agree with the image of~\eqref{eq:edgehomo}, e.g., see \cite[Table 2]{GPPV}. However, the question of realizing a given integral modular form as a field theory story is not completely settled in the physics literature, see~\cite{LinPei}. The mathematical expectation is that some (chiral) theories fail to extend down in the sense of~\S\ref{sec:final} below, and extending down need not be unique. This prevents certain modular forms from being realized as partition functions, and also provides multiple possible realizations of a given modular form as the partition function of a field theory. This expectation fits with~\eqref{eq:edgehomo} being neither surjective nor injective. 
\end{rmk}



\subsection{Repackaging as representations of a bordism category}\label{sec:21bord}
One can repackage the analytic data in Definitions~\ref{defn:QFTsystem} and~\ref{defn:SQFTsystem} as part of a representation of the bordism categories~$\EBord_{2}$ and~$\EBord_{2|1}$, respectively. We briefly review these bordism categories, referring to~\cite{ST11} for details. 

The objects of $\EBord_{2}$ are finite disjoint unions of metrized circles $S^{\Ra}_r$ and $S^{\NS}_r$ with spin structure, where $r\in \R_{>0}$ records the radius, $\Ra$ stands for the \emph{Ramond} (alias: periodic, odd, or nonbounding) spin structure and $\NS$ stands for \emph{Neveu--Schwarz} (alias: antiperiodic, even, or bounding) spin structure. The morphisms of $\EBord_{2}$ are flat 2-manifolds with spin structure and geodesic boundary, as pictured in Figure~\ref{fig:2EBord}. In addition to the cylinders and tori in Figure~\ref{fig:2d}, changing source and target data of a cylinder provides the ``macaronis" in Figure~\ref{fig:2EBord}. From the discussion at the beginning of \S\ref{sec:41}, the moduli space of Euclidean macaronis for a fixed source and target is $\HH/r\Z$, where $r\in \R_{>0}$ records the radius, and $\tau=\frac{x+it}{2\pi r}\in \HH$ gives the height $t$ of a cylinder and $x$ is the difference in boundary parameterizations. Composition corresponds to addition in $\HH$, e.g., $\rC^\Ra_{r,\tau}\circ \rC^\Ra_{r,\tau'}\simeq \rC^\Ra_{r,\tau+\tau'}.$
There are 4 spin structures on the torus $\C/2\pi r \Z\oplus 2 \pi \tau\Z$. The action by the metaplectic group~$\MP_2(\Z)$ permutes the three even spin structures, while the odd spin structure is preserved by this action. We emphasize that the flat metric prohibits closed surfaces of higher genus. Flat metrics together with geodesic boundaries prohibits bordisms between circles of different sizes or different numbers of circles, e.g., pairs of pants are not morphisms in $\EBord_2$. 

Adding supersymmetry, the objects of $\EBord_{2|1}$ are $1|1$-dimensional super circles; the moduli of such are again parameterized by a spin structure and radius \cite[\S4.5]{DBEtorsion}, so we use the same notation, $S^{\Ra}_r,S^{\NS}_r\in \EBord_{2|1}$. Bordisms between these are $2|1$-dimensional Euclidean manifolds with geodesic boundary. Analogous to the $1|1$-dimensional case, this introduces odd directions in the moduli space of morphisms in $\EBord_{2|1}$. For example, super Euclidean cylinders in the Ramond sector are parameterized by $\R_{>0}\times \HH^{2|1}$, with composition 
\beq\label{Eq:Ramondcylinders}
&&\rC^\Ra_{r,\tau,\bar\tau,\theta}\circ \rC^\Ra_{r,\tau',\bar\tau',\theta'}\simeq \rC^\Ra_{r,\tau+\tau',\bar\tau+\bar\tau'+\theta\theta',\theta+\theta'},\qquad r\in \R_{>0}, \ (\tau,\bar\tau,\theta)\in \HH^{2|1}.
\eeq
This gives the geometric counterpart to the supersemigroup representation~\eqref{Eq:ogsupersemi}.

To set up the degree~$n$ super Euclidean field theories relevant to Conjecture~\ref{conj}, we first must describe the degree twist. Let $\cCl(L_r\R^n)$ denote the Clifford algebra of the vector space $L_r\R^n=C^\infty(S_r^1;\R^n)$ for bilinear pairing $\int_0^{2\pi r} \langle \psi(t),\phi(t)\rangle dt$. The rotation action on~$L_r\R^n$ by~$S^1_r$ preserves this bilinear pairing, leading to an $S^1_r$-action on $\cCl(L_r\R^n)$ by algebra automorphisms. For each $\tau\in \HH/r\Z$, there is an algebra map $\varphi_r(\tau) \colon \cCl(L_r\R^n)\to \cCl(L_r\R^n)$ gotten from analytic continuation of this $S^1_r$-action. 
We refer to~\cite[page~59]{ST11} for details. 

\begin{defn}[{Degree twist, sketch of Ramond sector only}]\label{defn:21degreetw}
The \emph{degree twist} has values 
\beq\label{eq:degreetwist21}
\alpha_n\colon \EBord_{2|1} \to \TA,\qquad  \begin{array}{c} \alpha_n(S^{\Ra}_r)=\cCl(L_r\R^n), \\ \alpha_n(\rC^\Ra_{r,\tau,\bar\tau,\theta})=_{^{\varphi_r(\tau)}\cCl(L_r\R^n)}\! \cCl(L_r\R^n)_{\cCl(L_r\R^n)}, \\ 
\alpha_n(\upmac{}_{r,\tau,\bar\tau,\theta})= \cCl(L_r\R^n)_{^{\varphi_r(\tau)}\cCl(L_r\R^n)\otimes \cCl(L_r\R^n)} \\
\alpha_n(\downmac{}_{r,\tau,\bar\tau,\theta})=_{^{\varphi_r(\tau)}\cCl(L_r\R^n)\otimes \cCl(L_r\R^n)}\! \cCl(L_r\R^n),  \\
\alpha_n(\torus{}_{r,\tau,\bar\tau,\theta})=\Pf^{\otimes -n} \end{array}
\eeq
where the modules all come from $\cCl(L_r\R^n)$ as a bimodule over itself with left action twisted by the algebra automorphism $\varphi_r(\tau)$, and $\Pf$ is the Pfaffian line bundle, see Remark~\ref{rmk:etaMF}. These values are independent of $(\bar\tau,\theta)$, and hence the degree twist is chiral. 
\end{defn}

\begin{rmk}
The ``superfication" functor $\EBord_2\to \EBord_{2|1}$ \cite[Eq.~4.14]{ST11} allows one to restrict a twist $\alpha\colon \EBord_{2|1}\to \TA$ to a twist $\EBord_2\to \TA$. In particular, the twists~\eqref{eq:degreetwist21} determine twists for the 2-dimensional Euclidean bordism category. 
\end{rmk}

\begin{rmk}
The Neveu--Schwarz values are omitted from Definition~\ref{defn:21degreetw}, but come from a similar Clifford algebra of the vector space $C^\infty(S^1;\mu \otimes \R^n)$ where $\mu$ is the M\"obius bundle, i.e., the spinor bundle for the Neveu--Schwarz spin structure. A detailed sketch of the degree twist as a functor out of the Euclidean bordism category is given in~\cite[\S4.3]{ST11}, see also~\cite{Matthias1}. The theory is chiral, and so expected to admit a canonical supersymmetric extension, generalizing Example~\ref{ex:chiral}. However, the details of this extension in the cobordism framework of~\cite{ST11} have not been worked out. 
\end{rmk}
%



Let $\QFT_{2|1}^\alpha(\X)$ denote the groupoid of $\alpha$-twisted reflection positive $2|1$-dimensional Euclidean field theories over $\X$ (see Definition~\ref{def:QFT}). Reflection positivity is more subtle in the presence of chiral supersymmetry, but ultimately the data of a reflection structure follows a similar pattern as in Definition~\ref{defn:RP1}. 
We start by analyzing the groupoids $\QFT_2^n:=\QFT_2^{\alpha_n}(\pt)$ and $\QFT^n_{2|1}:=\QFT_{2|1}^{\alpha_n}(\pt)$ for field theories over the point twisted by the degree twist~\eqref{eq:degreetwist21}.

\begin{hyp}[{\cite[Theorem 1.0.2]{ST04} \cite[Theorem 2.2.2]{PokmanPhD}}] \label{hyp:QFTandBord}
There are map of stacks 
\beq\label{Eq:QFTmaps}
\QFT_{2}^n\to \qft^n \qquad{\rm and} \qquad \QFT_{2|1}^n\to \sqft^n
\eeq 
gotten by restriction to the Ramond sector. Together with Theorems~\ref{thm:Ktate} and~\ref{thm:modularity}, the stack $\QFT_{2|1}^n$ admits a map to the $n$th space in the spectrum representing $\K(\!(q)\!)$ and a map to the set of integral modular forms. 
\end{hyp}

\begin{proof}[Sketch.]
We focus on $\QFT_{2|1}^n\to \sqft^n$ below; the case without supersymmetry is easier.

First one shows that twisted field theories for~\eqref{eq:degreetwist21} restricted to the Ramond sector determine modules over the ordinary Clifford algebras~$\cCl_n$, following~\cite[\S6]{ST11}. The idea is to decompose the algebra~$\cCl(L_r\R^n)$ according to its Fourier modes, i.e., weight spaces for the $S^1_r$-action. These weights provide a dense inclusion
$$
\cCl_n\otimes \bigotimes_{k>0} \cCl(\C^n_k\oplus \C^n_{-k})) \hookrightarrow \cCl(L_r\R^n)
$$
where $\cCl_n$ is the Clifford algebra of the $S^1_r$-fixed points (i.e., constant maps) and $\cCl(\C^n_k\oplus \C^n_{-k})$ is the Clifford algebra of loops with of Fourier numbers $k$ and $-k$. For the pairing $\int_{S^1}\langle \phi(t),\psi(t)\rangle dt$, such loops only pair nontrivially with each other, and the pairing is the hyperbolic pairing. Hence, we obtain the Clifford algebras $\cCl(\C^n_k\oplus \C^n_{-k})\subset \cCl(L_r\R^n)$ as subalgebras. Furthermore, for each $k$, $\cCl(\C^n_k\oplus \C^n_{-k})\simeq \cCl_{n,n}$ is Morita equivalent to~$\C$. The completion of the tensor product of Morita bimodules gives the desired functor from $\cCl(L_r\R^n)$-modules to $\cCl_n$-modules, see~\cite[Eq.~6.2]{ST11}. 

Applying this to the value of $E\in \QFT_{2|1}^n$ on the subcategory of Ramond sector super annuli~\eqref{Eq:Ramondcylinders} produces a $\cCl_n$-linear supersemigroup representation of the form~\eqref{Eq:ogsupersemi}. This gives the data of an object in $\sqft^n$. The value of $E$ on Ramond--Ramond super tori imposes the modularity condition~\eqref{eq:anomalous} \cite[page~60]{ST11}. Applying an argument analogous to the one in Hypothesis~\ref{hyp:QMandBord}, an isomorphism between twisted field theories comes from a projective isomorphism between state spaces. Hence we obtain a functor $\QFT^n_{2|1}\to \sqft^n$. 
\ep




\begin{figure}
\beq
&&\begin{tikzpicture}[baseline=(basepoint)];
	\begin{pgfonlayer}{nodelayer}
		\node [style=none] (0) at (-1.5, -0.25) {};
		\node [style=none] (2) at (-3, 0.25) {};
		\node [style=none] (3) at (-3.5, 0.25) {};
		\node [style=none] (4) at (-2, -0.25) {};
		\node [style=none] (5) at (-4.5, 0.25) {};
		\node [style=none] (6) at (-5, 0.25) {};
		\node [style=none] (13) at (-3, -0.25) {};
		\node [style=none] (14) at (-3.5, -0.25) {};
		\node [style=none] (15) at (0.5, 1.5) {};
		\node [style=none] (16) at (0, 1.5) {};
		\node [style=none] (17) at (0.5, 0) {};
		\node [style=none] (18) at (0, 0) {};
		\node [style=none] (19) at (0.5, -1.5) {};
		\node [style=none] (20) at (0, -1.5) {};
		\node [style=none] (24) at (-.75, 0) {$\xrightarrow{\rm compose}$};
	\end{pgfonlayer}
	\begin{pgfonlayer}{edgelayer}
		\draw [bend right=90, looseness=1.25,thick] (0.center) to (4.center);
		\draw [bend right=105, looseness=1.25,thick,dashed] (2.center) to (3.center);
		\draw [bend left=90, looseness=0.75,thick] (2.center) to (3.center);
		\draw [bend left=90,thick] (0.center) to (4.center);
		\draw [bend right=90, looseness=1.50,thick] (2.center) to (6.center);
		\draw [bend right=90, looseness=1.50,thick] (3.center) to (5.center);
		\draw [bend left=90,thick] (5.center) to (6.center);
		\draw [bend right=90,thick,dashed] (5.center) to (6.center);
		\draw [bend right=90, looseness=1.25,thick] (13.center) to (14.center);
		\draw [bend left=90,thick] (13.center) to (14.center);
		\draw [bend right=90, looseness=1.75,thick] (13.center) to (4.center);
		\draw [bend right=90, looseness=1.50,thick] (14.center) to (0.center);
		\draw [bend left=270, looseness=0.75,thick] (15.center) to (16.center);
		\draw [bend left=270,thick,dashed] (17.center) to (18.center);
		\draw [bend left=270,thick,dashed] (19.center) to (20.center);
		\draw[thick] (15.center) to (19.center);
		\draw[thick] (16.center) to (20.center);
		\draw [bend left=90,thick] (15.center) to (16.center);
		\draw [bend left=90, looseness=0.75,thick] (17.center) to (18.center);
		\draw [bend left=90, looseness=0.75,thick] (19.center) to (20.center);
	\end{pgfonlayer}
		\path (0,0) coordinate (basepoint);
\end{tikzpicture}\qquad \qquad 
\begin{tikzpicture}[baseline=(basepoint)];
	\begin{pgfonlayer}{nodelayer}
		\node [style=none] (0) at (2.5, -0.5) {};
		\node [style=none] (1) at (2, -0.5) {};
		\node [style=none] (2) at (0.5, -0.5) {};
		\node [style=none] (3) at (0, -0.5) {};
		\node [style=none] (4) at (0.5, 0) {};
		\node [style=none] (5) at (0, 0) {};
		\node [style=none] (6) at (0.5, 1.5) {};
		\node [style=none] (7) at (0, 1.5) {};
		\node [style=none] (8) at (2.5, 0) {};
		\node [style=none] (9) at (2, 0) {};
		\node [style=none] (10) at (2.5, 1.5) {};
		\node [style=none] (11) at (2, 1.5) {};
		\node [style=none] (12) at (6.5, 0) {};
		\node [style=none] (13) at (6, 0) {};
		\node [style=none] (14) at (4.5, 0) {};
		\node [style=none] (15) at (4, 0) {};
		\node [style=none] (18) at (4.5, 1.25) {};
		\node [style=none] (19) at (4, 1.25) {};
		\node [style=none] (20) at (6.5, 1.25) {};
		\node [style=none] (21) at (6, 1.25) {};
		\node [style=none] (22) at (3.25, 0) {$\xrightarrow{\rm compose}$};
	\end{pgfonlayer}
	\begin{pgfonlayer}{edgelayer}
		\draw [bend left=105, looseness=1.25,thick] (0.center) to (1.center);
		\draw [bend right=90, looseness=0.75,thick] (0.center) to (1.center);
		\draw [bend left=90, looseness=1.50,thick] (0.center) to (3.center);
		\draw [bend left=90, looseness=1.50,thick] (1.center) to (2.center);
		\draw [bend right=90,thick] (2.center) to (3.center);
		\draw [bend left=90,thick] (2.center) to (3.center);
		\draw [bend right=270,thick] (4.center) to (5.center);
		\draw [bend right=285,thick,dashed] (6.center) to (7.center);
		\draw [bend right=90, looseness=0.75,thick,dashed] (4.center) to (5.center);
		\draw [bend right=90,thick] (6.center) to (7.center);
		\draw [thick] (6.center) to (4.center);
		\draw [thick] (7.center) to (5.center);
		\draw [bend right=270,thick] (8.center) to (9.center);
		\draw [bend right=270,thick,dashed] (10.center) to (11.center);
		\draw [bend right=90, looseness=0.75,thick,dashed] (8.center) to (9.center);
		\draw [bend right=75, looseness=0.75,thick] (10.center) to (11.center);
		\draw [thick] (10.center) to (8.center);
		\draw [thick] (11.center) to (9.center);
		\draw [bend left=105, looseness=1.25,thick] (12.center) to (13.center);
		\draw [bend right=90, looseness=0.75,thick,dashed] (12.center) to (13.center);
		\draw [bend left=90, looseness=1.50,thick] (12.center) to (15.center);
		\draw [bend left=90, looseness=1.50,thick] (13.center) to (14.center);
		\draw [bend right=90,thick,dashed] (14.center) to (15.center);
		\draw [bend left=90,thick] (14.center) to (15.center);
		\draw [bend right=270,thick,dashed] (18.center) to (19.center);
		\draw [bend right=75, looseness=0.75,thick] (18.center) to (19.center);
		\draw [bend right=270,thick,dashed] (20.center) to (21.center);
		\draw [bend right=75, looseness=0.75,thick] (20.center) to (21.center);
		\draw [thick] (15.center) to (19.center);
		\draw [thick] (14.center) to (18.center);
		\draw [thick] (13.center) to (21.center);
		\draw [thick] (12.center) to (20.center);
	\end{pgfonlayer}
	\path (0,0) coordinate (basepoint);
\end{tikzpicture}\nonumber\\
&&\qquad\begin{array}{cc} (\upmac{}_{r,\tau,\bar\tau,\theta})\circ (\downmac{}_{r,\tau',\bar\tau',\theta'})\\
\simeq \rC^\Ra_{r,\tau'+\tau,\bar\tau'+\bar\tau+\theta'\theta,\theta'+\theta}
\end{array} \qquad\qquad \begin{array}{c} (\rC^\Ra_{r,\tau,\bar\tau,\theta}\coprod \rC^\Ra_{r,\tau'',\bar\tau'',\theta''})\circ (\downmac{}_{r,\tau',\bar\tau',\theta'})\\ \simeq \downmac{}_{r,\tau+\tau'+\tau'',\bar\tau'+\bar\tau'+\bar\tau''+\theta\theta'+(\theta+\theta')\theta'',\theta+\theta'+\theta''} \end{array}\nonumber
\eeq
\caption{
Generators and relations for $2|1$-dimensional bordism categories are analogous to Figures~\ref{fig:EBord}  and~\ref{fig:11EBord}, but the geometric parameters are given by $(r,\tau,\bar\tau,\theta)\in \R_{>0}\times \HH^{2|1}$. After keeping track of source and targets, composition corresponds to multiplication in the semigroup $\HH^{2|1}$. As usual, bordisms are read in the positive (upward) direction. On the left, gluing a sad macaroni to a happy macaroni gives a cylinder. On the right, gluing cylinders onto a happy macaroni gives a happier macaroni. 
}
\label{fig:2EBord}
\end{figure}
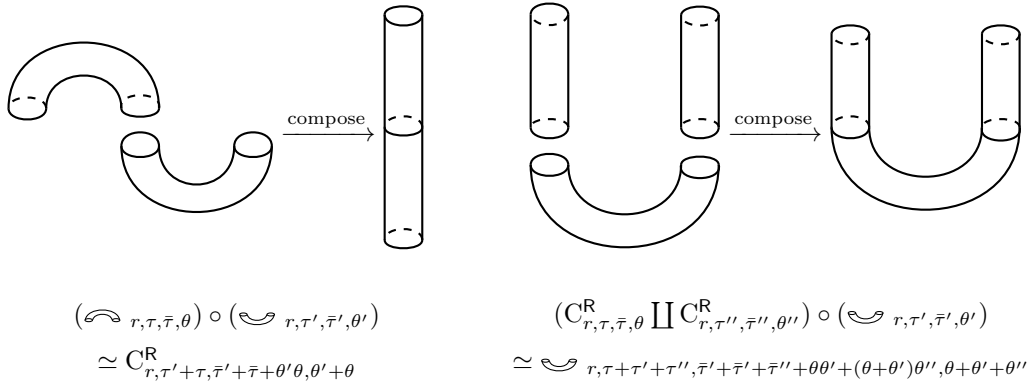

\subsection{Additional insights from the bordism framework}\label{sec:21add}
In contrast to Hypothesis~\ref{hyp:SQMandBord}, the map $\QFT_{2|1}^n\to \sqft^n$ in Hypothesis~\ref{hyp:QFTandBord} is no longer an epimorphism of stacks: the bordism definition contains strictly more information than Definition~\ref{defn:SQFTsystem}, e.g., values~$E(S^{\NS}_r)$ on the Neveu--Schwarz sector. 
We highlight a few additional enhancements below. 

\begin{ex}[Real structures]
In the bordism framework of~\cite{ST11}, both the objects and morphisms in $\EBord_{2|1}(\X)$ form stacks, meaning the isometry group of a bordism is encoded in the structure. This leads to an additional $\Z/2$-automorphism on objects that is not present in the 1-dimensional case, boiling down to the fact that $-\id\in \O_2$ is an orientation preserving isometry of $\R^2$ (in contrast to $-1\in \O(1)$ that reverses the orientation on~$\R$). Physically, this automorphism is a parity-time reversal. For reflection positive theories, this symmetry determines an isomorphism between the space of states and its complex conjugate, i.e., is the data of a real structure~\cite[\S8.4]{DBEtorsion}, thereby refining the $\K(\!(q)\!)$-valued cocycle map in \ref{hyp:QFTandBord} to a $\KO(\!(q)\!)$-valued map. This argument can also be understood in physical language, see~\cite{Witten_Elliptic} and \cite[3.2.2]{GPPV}. 
\end{ex}

\begin{ex}[Torsion invariants of field theories]
Following Remark~\ref{rmk:Chernchar}, the supermoduli spaces encoded by $\EBord_{2|1}(\X)$ can be quite subtle. In particular, the moduli of super Euclidean tori with a map to $\X$ is surprisingly rich. For $\X=X$ a smooth manifold, in~\cite{DBEChern} it was shown that functions on this moduli space determine a cocycle representative of a class in $\H(M;\MF)$, i.e., ordinary cohomology with values in modular forms. Subtle features of the moduli space of super tori encode holomorphic anomaly cancelation data in the cocycle representative, recovering the information anticipated by physical arguments from~\cite{GJFW,GJF2}. Together with Witten's partial construction of the supersymmetric sigma model (see Example~\ref{ex:Wittengenus}), this anomaly cancelation data can be used to provide a quantum field theory interpretation of Bunke and Naumann's torsion invariant of $4k-1$-dimensional string manifolds~\cite{BunkeNaumann}. These ideas are contextualized in the framework of~\cite{ST11} in \cite{DBEtorsion}, providing the first examples of $\TMF$-torsion captured by 2-dimensional supersymmetric field theories. 
\end{ex}


\begin{ex}[Symmetries of the degree twist and the elliptic Thom class]\label{ex:freefermion2} 
In analogy to Example~\ref{ex:freefermion11}, the degree twist from Definition~\ref{defn:21degreetw} can be constructed as the quantization of a theory of $n$ chiral free fermions~\cite[\S2.6]{ST04}, defined by the fields and action functional
$$
\psi\in \Gamma(\Sigma;\bS\otimes \R^n),\qquad S(\psi)=\frac{1}{2}\int_\Sigma \langle \psi,\bar\partial \psi\rangle
$$
where $\Sigma$ is a spin Riemann surface with spinor bundle $\bS\to \Sigma$ and chiral Dirac operator~$\bar\partial$. This theory has classical (gauge) symmetries given by maps $\Sigma\to \O(n)$, acting on~$\psi$ by rotating $\R^n$. In the quantum theory, this leads to a projective action of the loop group~$L\O(n)$ on the Hilbert space of states that intertwines with the action of the Clifford algebra $\cCl(L\R^n)$. This loop group representation is precisely the  fermionic Fock space construction for the level one representation of $\widehat{L\Spin(n)}$, the central extension of the loop group $L\Spin(n)$ classified by $\frac{p_1}{2}\in \H^4(B\Spin(n);\Z)$~\cite[\S12]{PressleySegal}. This chiral conformal field theory determines a Euclidean field theory with $\mathcal{N}=(0,1)$ supersymmetry. Under Conjecture~\ref{conj}, it is expected to recover the universal elliptic Thom class, see~\cite[Theorem 1.0.3]{ST04} and \cite{DouglasHenriques}. 

%
%

\end{ex}


\subsection{The final frontier: Extending down}\label{sec:final}

\begin{figure}
\beq\nonumber
\begin{tikzpicture}
		\node [style=none] (0) at (-5.5, -2) {};
		\node [style=none] (A) at (-4.75, -2.15) {$S_{\ell,\tau,\bar\tau,\theta}$};
		\node [style=none] (1) at (-4, -2) {};
		\node [style=none] (2) at (-3, -2) {};
		\node [style=none] (3) at (-1.5, -2) {};
		\node [style=none] (B) at (-2.25, -2.35) {$S_{\ell',\tau,\bar\tau,\theta}$};
				\node [style=none] (C) at (-.5, -.6) {$\stackrel{{\rm glue}}{\implies}$};
		\node [style=none] (4) at (-0.5, -2) {};
		\node [style=none] (5) at (2.5, -2) {};
		\node [style=none] (6) at (-4.5, 1) {};
		\node [style=none] (7) at (-3, 1) {};
		\node [style=none] (8) at (-2, 1) {};
		\node [style=none] (9) at (-0.5, 1) {};
		\node [style=none] (10) at (0.5, 1) {};
		\node [style=none] (11) at (3.5, 1) {};
		\node [style=none] (12) at (4.75, -2) {};
		\node [style=none] (13) at (6.25, -2) {};
		\node [style=none] (14) at (5.75, 1) {};
		\node [style=none] (15) at (7.25, 1) {};
		\node [style=none] (E) at (5.5, -2.75) {$C_{\ell+\ell',\tau,\bar\tau,\theta}$};
		\node [style=none] (F) at (4, -.5) {$\stackrel{{\rm glue}}{\implies}$};
		\node [style=none] (16) at (1, -2) {};
		\node [style=none] (D) at (1, -2.25) {$S_{\ell+\ell',\tau,\bar\tau,\theta}$};
		\node [style=none] (17) at (2, 1) {};
		\draw[thick] (6.center) to (0.center);
		\draw[thick] (7.center) to (1.center);
		\draw[thick] (8.center) to (2.center);
		\draw[thick] (9.center) to (3.center);
		\draw[thick] (10.center) to (4.center);
		\draw[thick] (11.center) to (5.center);
		\draw [bend left,thick] (0.center) to (1.center);
		\draw [bend left=15,thick] (6.center) to (7.center);
		\draw [bend left=15,thick] (8.center) to (9.center);
		\draw [bend right=15,thick] (2.center) to (3.center);
		\draw[thick] (12.center) to (14.center);
		\draw[thick] (15.center) to (13.center);
		\draw[thick] [in=-90, out=-90,thick] (12.center) to (13.center);

		\draw [bend left=80, looseness=1.25,thick] (12.center) to (13.center);
		\draw [bend left=80, looseness=1.25,thick] (14.center) to (15.center);
		\draw [bend right=90,thick,dashed] (14.center) to (15.center);
		\draw [bend left=15,thick] (4.center) to (16.center);
		\draw [bend right=15,thick] (16.center) to (5.center);
		\draw [bend left=15,thick] (10.center) to (17.center);
		\draw [bend right=15, looseness=1.25,thick] (17.center) to (11.center);
\end{tikzpicture}
\eeq
\caption{
Horizontal composition of 2-morphisms in the fully-extended $2|1$-dimensional Euclidean bordism category. The strips $S_{\ell,\tau,\bar\tau,\theta}$ are indexed by a length of the incoming and outgoing string as well as a parameter in $\HH^{2|1}\subset \R^{2|1}$ measuring an angle and a height of the strip. We view these as open strings that propagate in the upward (time) direction, and have a gluing operation in the horizontal (space) direction. Gluing the boundaries of a strip together recovers the previously considered category of cylinders.}
\label{fig:open}
\end{figure}
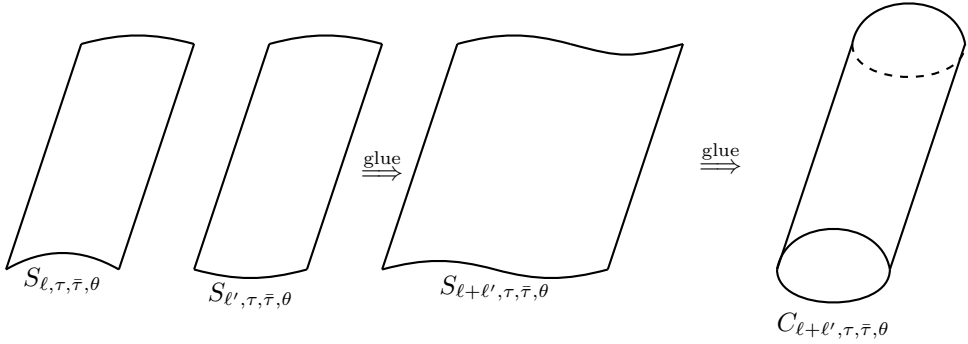

A fully-extended $2|1$-dimensional Euclidean bordism 2-category encodes the cutting and gluing of $2|1$-dimensional super Euclidean manifolds with corners: the theory of closed strings encoded by the operations in Figure~\ref{fig:2EBord} is refined to a theory of open strings pictured in Figure~\ref{fig:open}. The expectation is that extending down will refine the $\KO(\!(q)\!)$ and $\MF^\Z$-valued invariants described in Hypothesis~\ref{hyp:QFTandBord} to invariants of field theories valued in $\TMF$, i.e., the cocycle map in~\eqref{eq:conj}. 

Some rough features of such fully-extended field theories are known e.g.,~\cite[3.5]{MooreKtheory} and \cite[\S4]{SegalDbrane}. However, a lot of mystery remains: extending down is not uniquely defined, and no candidate definition has been spelled out in detail. As usual, the work can be divided into two interrelated pieces: (i) the geometry of the bordism 2-category and (ii) the analysis of the target 3-category. The language of $n$-fold Segal spaces provides a promising framework in which to study fully-extended bordism categories in general (e.g., see~\cite{theoclaudia,GradyPavlov}), and fully-extended $2|1$-dimensional Euclidean bordisms in particular. One candidate for the target is the symmetric monoidal 3-category of conformal nets~\cite{DouglasHenriques,CN1,CN2,CN3,CN4}, though there are some anticipated technical modifications as this 3-category does not deloop the 2-category $\TA$ for 1-extended field theories from Definition~\ref{defn:twistedQFT} and~\cite{ST11}. 


One way to guide the choices in extending down is to get a better handle on the examples. Perhaps the most fundamental example is the degree twist~\eqref{eq:degreetwist21} and its expected connection with the elliptic Thom class, as sketched in Example~\ref{ex:freefermion2}. This generalizes known structures in K-theory, see Example~\ref{ex:freefermion11}. Indeed, a fully-extended degree~$n$ field theory implementing a universal twisted Thom isomorphism of the form
$$
\QFT_{2|1}^k(\pt\nsq \O(n))\to \QFT_{2|1}^{n+k+\lambda}(\R^n\nsq \O(n))_c,\qquad E\mapsto E\boxtimes \Th_n
$$ 
would automatically construct a string oriented multiplicative cohomology theory from 2-dimensional supersymmetric field theories. Following Lurie's suggestion from \cite[\S5.5]{Lurie_Elliptic}, higher equivariant structures in field theories could then be used to construct a comparison with $\TMF$ as in~\eqref{eq:conj}. Evidence for higher equivariance comes from the expectation that symmetries of the free fermion field theory are the string 2-group~\cite{ST04,ChrisAndrestring,SchommerPries,Waldorfstring}. 

We offer one further speculation. Let $\chi\QFT_2^n\subset \QFT_2^n$ denote the subcategory of degree~$n$ theories whose image under~\eqref{Eq:QFTmaps} is chiral in the sense of Example~\ref{ex:chiral}. As in Corollary~\ref{cor:chiral}, it is expected that such chiral conformal field theories determine supersymmetric ones, giving an inclusion 
\beq\label{eq:chiralinclusion}
\chi\QFT_2^\bullet\xhookrightarrow{\rm group\ completion?} \QFT_{2|1}^\bullet.
\eeq
Analogizing with~\eqref{Eq:1dgroupcomplete}, one might hope that~\eqref{eq:chiralinclusion} implements a group completion of chiral field theories. Indeed, in dimension~1 the topological field theories are precisely the ones invariant under the renormalization group flow, and in dimension~2 these are the chiral conformal theories. Together with Conjecture~\ref{conj},~\eqref{eq:chiralinclusion} would realize $\TMF$ as a group completion of chiral conformal field theories. At present it is unclear how to parse this statement, but it is likely that the rigid objects in the image of~\eqref{eq:chiralinclusion} play an important role in Conjecture~\ref{conj} just as vector bundles are the basic geometric objects in K-theory.

\bibliographystyle{amsalpha}
\bibliography{references}

\newcommand{\etalchar}[1]{$^{#1}$}
\providecommand{\bysame}{\leavevmode\hbox to3em{\hrulefill}\thinspace}
\providecommand{\MR}{\relax\ifhmode\unskip\space\fi MR }
\providecommand{\MRhref}[2]{%
  \href{http://www.ams.org/mathscinet-getitem?mr=#1}{#2}
}
\providecommand{\href}[2]{#2}
\begin{thebibliography}{DHVW85}

\bibitem[AB02]{AndoBasterra}
M.~Ando and M.~Basterra, \emph{The {W}itten genus and equivariant elliptic
  cohomology}, Math. Z. \textbf{240} (2002), no.~4, 787--822.

\bibitem[ABS64]{ABS}
M.~Atiyah, R.~Bott, and A.~Shapiro, \emph{Clifford modules}, Topology
  \textbf{3} (1964), no.~suppl. 1, 3--38.

\bibitem[AG83]{Alvarez}
L.~Alvarez-Gaum{\'e}, \emph{Supersymmetry and the {Atiyah-Singer} index
  theorem}, Communications in Mathematical Physics \textbf{90} (1983),
  161--173.

\bibitem[AHR10]{AHR}
M.~Ando, M.~Hopkins, and C.~Rezk, \emph{Multiplicative orientations of
  {KO}-theory and the spectrum of topological modular forms}, Draft (2010).

\bibitem[AKL23]{TMFmoon2}
J.~Albert, J.~Kaidi, and Y.-H. Lin, \emph{Topological modularity of
  supermoonshine}, Progress of Theoretical and Experimental Physics
  \textbf{2023} (2023).

\bibitem[And00]{Ando}
M.~Ando, \emph{Power operations in elliptic cohomology and representations of
  loop groups}, Transactions of the American Mathematical Society \textbf{352}
  (2000).

\bibitem[AO21]{AganagicOkounkov}
M.~Aganagic and A.~Okounkov, \emph{Elliptic stable envelope}, J. Amer. Math.
  Soc. \textbf{34} (2021).

\bibitem[AS69]{AtiyahSingerskew}
M.~Atiyah and I.~Singer, \emph{Index theory for skew-adjoint {Fredholm}
  operators}, Inst. Hautes \'etudes Sci. Publ. Math. \textbf{37} (1969).

\bibitem[AS71]{AtiyahSinger4}
\bysame, \emph{The index of elliptic operators {IV}}, Annals of Mathematics
  \textbf{93} (1971).

\bibitem[AS04]{AtiyahSegaltwistedK}
M.~Atiyah and G.~Segal, \emph{Twisted {K}-theory}, Ukrainskyj Matematychnyj
  Visnyk \textbf{1} (2004).

\bibitem[Ati85]{AtiyahCircular}
M.~Atiyah, \emph{Circular symmetry and stationary phase approximation},
  Proceedings of the conference in honor of {L. Schwartz, Ast\'erisque}
  \textbf{131} (1985).

\bibitem[Ati89]{AtiyahTQFT}
\bysame, \emph{Topological quantum field theories}, {Publications
  Math\'ematiques de l'IH\'ES} \textbf{68} (1989).

\bibitem[Ati90]{AtiyahJonesWitten}
\bysame, \emph{The {Jones-Witten} invariants of knots}, S\'eminaire Bourbaki
  \textbf{715} (1989/90).

\bibitem[BBES22]{Powerops}
T.~Barthel, D.~Berwick-Evans, and N.~Stapleton, \emph{Power operations in the
  {Stolz--Teichner} program}, Geom. Topol. \textbf{26} (2022).

\bibitem[BDH15]{CN1}
A.~Bartels, C.~Douglas, and A.~Henriques, \emph{Conformal nets {I}:
  coordinate-free nets}, Int. Math. Res. Not. \textbf{13} (2015).

\bibitem[BDH17]{CN2}
\bysame, \emph{Conformal nets {II}: conformal blocks}, Comm. Math. Phys.
  \textbf{354} (2017).

\bibitem[BDH18]{CN4}
\bysame, \emph{Conformal nets {IV}: the 3-category}, Algebr. Geom. Topol.
  \textbf{18} (2018).

\bibitem[BDH19]{CN3}
\bysame, \emph{Fusion of defects}, Mem. Amer. Math. Soc. \textbf{1237} (2019).

\bibitem[BE23a]{DBEChern}
D.~Berwick-Evans, \emph{{Chern characters for supersymmetric field theories}},
  Geom. Topol. \textbf{27} (2023), no.~5, 1947--1986.

\bibitem[BE23b]{DBEEFT}
\bysame, \emph{The families index for $1|1$-dimensional {Euclidean} field
  theories}, preprint (2023).

\bibitem[BE23c]{DBEtorsion}
\bysame, \emph{How do field theories detect the torsion in topological modular
  forms?}, preprint (2023).

\bibitem[BEH16]{DBEHan}
D.~Berwick-Evans and F.~Han, \emph{The equivariant {Chern} character as super
  holonomy on loop stacks}, preprint (2016).

\bibitem[BEP23]{BEPTFT}
D.~Berwick-Evans and D.~Pavlov, \emph{{Smooth one-dimensional topological field
  theories are vector bundles with connection}}, Algebr. Geom. Topol.
  \textbf{23} (2023), no.~8, 3707--3743.

\bibitem[BET19]{BET1}
D.~Berwick-Evans and A.~Tripathy, \emph{A model for complex analytic
  equivariant elliptic cohomology from quantum field theory}, preprint (2019).

\bibitem[BET21]{BET0}
\bysame, \emph{A de rham model for complex analytic equivariant elliptic
  cohomology}, Advances in Mathematics \textbf{380} (2021).

\bibitem[BGV92]{BGV}
N.~Berline, E.~Getzler, and M.~Vergne, \emph{Heat kernels and {Dirac}
  operators}, Springer, 1992.

\bibitem[Bis86]{Bismutindex}
J.-M. Bismut, \emph{The {Atiyah-Singer} index theorem for families of {Dirac}
  operators: Two heat equation proofs}, Inventiones mathematicae \textbf{83}
  (1986).

\bibitem[BN14]{BunkeNaumann}
U.~Bunke and N~Naumann, \emph{Secondary invariants for string bordism and
  topological modular forms}, Bull. Sci. Math. \textbf{138} (2014).

\bibitem[Bor92]{Borcherds}
R.~Borcherds, \emph{Monstrous moonshine and monstrous {Lie} superalgebras},
  Invent. Math. \textbf{109} (1992).

\bibitem[BT89]{BottTaubes}
R.~Bott and C.~Taubes, \emph{On the rigidity theorems of {W}itten}, J. Amer.
  Math. Soc. \textbf{2} (1989), no.~1, 137--186.

\bibitem[BW97]{Weinsteinquantize}
S.~Bates and A.~Weinstein, \emph{Lectures on the geometry of quantization},
  Berkeley mathematics lecture notes, American Mathematical Society, 1997.

\bibitem[Che08]{PokmanPhD}
P.~Cheung, \emph{Supersymmetric field theories and cohomology}, preprint
  (2008).

\bibitem[Cos10]{costello_WG1}
K.~Costello, \emph{A geometric construction of the {Witten} genus {I}},
  {Proceedings of the International Congress of Mathematicians} (2010).

\bibitem[Cos11]{costello_WG2}
\bysame, \emph{A geometric construction of the {Witten} genus {II}}, preprint
  (2011).

\bibitem[DEF{\etalchar{+}}99]{strings1}
P.~Deligne, P.~Etingof, D.~Freed, L.~Jeffrey, D.~Kazhdan, J.~Morgan,
  D.~Morrison, and E.~Witten, \emph{{Quantum Fields and Strings: {A} Course for
  Mathematicians, Volume 1}}, American Mathematical Society, 1999.

\bibitem[Dev96]{DevotoI}
J.~Devoto, \emph{Equivariant elliptic homology and finite groups}, Michigan
  Math. J. \textbf{43} (1996).

\bibitem[Dev98]{DevotoII}
\bysame, \emph{An algebraic description of the elliptic cohomology of
  classifying spaces}, Journal of Pure and Applied Algebra \textbf{130} (1998).

\bibitem[DFHH14]{DFHH}
C.~Douglas, J.~Francis, A.~Henriques, and M.~Hill, \emph{Topological modular
  forms}, American Mathematical Society, 2014.

\bibitem[DH]{ChrisAndrestring}
C.~Douglas and A.~Henriques, \emph{Geometric string structures}, preprint.

\bibitem[DH11]{DouglasHenriques}
\bysame, \emph{Topological modular forms and conformal nets}, Mathematical
  foundations of quantum field theory and perturbative string theory, Proc.
  Sympos. Pure Math., vol.~83, Amer. Math. Soc., Providence, RI, 2011,
  pp.~341--354.

\bibitem[DHVW85]{DHVW}
L.~Dixon, J.~Harvey, C.~Vafa, and E.~Witten, \emph{Strings on orbifolds},
  Nuclear Phys. B \textbf{261} (1985).

\bibitem[DM99]{DM}
P.~Deligne and J.~Morgan, \emph{Notes on supermanifolds}, Quantum Fields and
  Strings: {A} Course for Mathematicians, Volume 1 (P.~Deligne, P.~Etingof,
  D.~Freed, L.~Jeffrey, D.~Kazhdan, J.~Morgan, D.~Morrison, and E.~Witten,
  eds.), American Mathematical Society, 1999.

\bibitem[DMVV97]{DMVV}
R.~Dijkgraaf, G.~Moore, E.~Verlinde, and H.~Verlinde, \emph{Elliptic genera of
  symmetric products and second quantized strings}, Comm. Math. Phys.
  \textbf{185} (1997).

\bibitem[Dou10]{Douglas_space}
M.~Douglas, \emph{Spaces of quantum field theories}, preprint (2010).

\bibitem[Dum12]{Florin}
F.~Dumitrescu, \emph{$1|1$ parallel transport and connections}, Differential
  Geometry and its Applications \textbf{30} (2012).

\bibitem[EU14]{EspinozaUribe}
J.~Espinoza and B.~Uribe, \emph{Topological properties of the unitary group},
  JP Journal of Geometry and Topology \textbf{16} (2014).

\bibitem[FH21]{FreedHopkins}
D.~Freed and M.~Hopkins, \emph{Reflection positivity and invertible topological
  phases}, Geom. Topol. \textbf{25} (2021).

\bibitem[FHT11a]{FHTI}
D.~Freed, J.~Hopkins, and C.~Teleman, \emph{Loop groups and twisted
  {$K$}-theory {I}}, J. Topol. \textbf{4} (2011), no.~4, 737--798.

\bibitem[FHT11b]{FHT3}
D.~Freed, M.~Hopkins, and C.~Teleman, \emph{{Twisted K-theory and loop group
  representations III}}, Ann. Math. \textbf{174} (2011).

\bibitem[FM13]{FreedMoore}
D.~Freed and G.~Moore, \emph{Twisted equivariant matter}, Ann. Henri Poincare
  \textbf{14} (2013).

\bibitem[Fre87]{Freed_Det}
D.~Freed, \emph{On determinant line bundles}, Math. aspects of string theory
  (1987).

\bibitem[Fre99]{Freed5}
\bysame, \emph{{Five Lectures on Supersymmetry}}, American Mathematical
  Society, 1999.

\bibitem[Fre12]{Freedalg}
\bysame, \emph{Lectures on twisted {K}-theory and orientifolds}, lectures at
  {ESI} {Vienna} (2012).

\bibitem[Fre14]{Freedanomaly}
\bysame, \emph{{Anomalies and Invertible Field Theories}}, Proc. Symp. Pure
  Math. \textbf{88} (2014), 25--46.

\bibitem[FRV17]{FRV}
G.~Felder, R.~Rim\'anyi, and A.~Varchenko, \emph{Elliptic dynamical quantum
  groups and equivariant elliptic cohomology}, Arxiv preprint (2017).

\bibitem[FT14]{FreedTelemanRelative}
D.~Freed and C.~Teleman, \emph{Relative quantum field theory}, Commun. Math.
  Phys. \textbf{326} (2014).

\bibitem[Gan06]{GanterOrb}
N.~Ganter, \emph{Orbifold genera, product formulas and power operations}, Adv.
  Math. \textbf{205} (2006).

\bibitem[Gan07]{Ganterstringy}
\bysame, \emph{Stringy power operations in {Tate K-theory}}, preprint (2007).

\bibitem[Gan09]{GanterHecke}
\bysame, \emph{Hecke operators in equivariant elliptic cohomology and
  generalized {M}oonshine}, Groups and symmetries, CRM Proc. Lecture Notes,
  vol.~47, Amer. Math. Soc., Providence, RI, 2009, pp.~173--209.

\bibitem[Gan14]{GanterEllipticWCF}
\bysame, \emph{The elliptic {W}eyl character formula}, Compos. Math.
  \textbf{150} (2014), no.~7, 1196--1234.

\bibitem[GG14]{gradygwilliam}
R.~Grady and O.~Gwilliam, \emph{{One-dimensional Chern\textendash{}Simons
  theory and the \^A genus}}, Algebr. Geom. Topol. \textbf{14} (2014), no.~4,
  2299--2377.

\bibitem[GJF23]{GJF2}
D.~Gaiotto and T.~Johnson-Freyd, \emph{Mock modularity and a secondary elliptic
  genus}, Journal of High Energy Physics \textbf{2023} (2023).

\bibitem[GJFW21]{GJFW}
D.~Gaiotto, T.~Johnson-Freyd, and E.~Witten, \emph{A note on some minimally
  supersymmetric models in two dimensions}, Integrability, Quantization, and
  Geometry: II. Quantum Theories and Algebraic Geometry (S.~Novikov et~al.,
  ed.), vol. 103, Proc. Symposia Pure Math., 2021.

\bibitem[GKV95]{GKV}
V.~Ginzburg, M.~Kapranov, and E.~Vasserot, \emph{Elliptic algebras and
  equivariant elliptic cohomology}, Arxiv preprint (1995).

\bibitem[GL22]{GanterLaures}
N.~Ganter and G.~Laures, \emph{Codes, vertex operators and topological modular
  forms}, Bulletin of the London Mathematical Society \textbf{54} (2022),
  no.~4, 1167--1196.

\bibitem[GM23]{LenartDavid}
D.~Gepner and L.~Meier, \emph{On equivariant topological modular forms},
  Compositio Mathematica \textbf{159} (2023), 2638--2693.

\bibitem[GMS00]{CDO1}
V.~Gorbounov, F.~Malikov, and V.~Schechtman, \emph{Gerbes of chiral
  differential operators}, Math. Res. Lett. \textbf{7} (2000).

\bibitem[GMS04]{CDO2}
\bysame, \emph{Gerbes of chiral differential operators {II}. vertex
  algebroids}, Invent. Math. \textbf{155} (2004).

\bibitem[GP20]{GradyPavlov}
D.~Grady and D.~Pavlov, \emph{Extended field theories are local and have
  classifying spaces}, arxiv preprint (2020).

\bibitem[GPPV21]{GPPV}
S.~Gukov, D.~Pei, P.~Putrov, and C.~Vafa, \emph{4-manifolds and topological
  modular forms}, J. High Energ. Phys. \textbf{84} (2021).

\bibitem[GR13]{GanterRam}
N.~Ganter and A.~Ram, \emph{Generalized {S}chubert calculus}, Journal of the
  Ramanujan Mathematical Society \textbf{28A} (2013).

\bibitem[Gro07]{Grojnowski}
I.~Grojnowski, \emph{Delocalised equivariant elliptic cohomology}, Elliptic
  cohomology: Geometry, applications, and higher chromatic analogues ({H.
  Miller and D. Ravenel}, ed.), London Mathematical Society, 2007.

\bibitem[Gun16]{Sam}
S.~Gunningham, \emph{Spin {Hurwitz} numbers and topological quantum field
  theory}, Geometry and Topology \textbf{20} (2016).

\bibitem[Han08]{Han}
F.~Han, \emph{Supersymmetric {QFT}s, super loop spaces and {Bismut-Chern}
  character}, {PhD} Thesis (2008).

\bibitem[HK04]{HuKriz}
P.~Hu and I.~Kriz, \emph{Conformal field theory and elliptic cohomology}, Adv.
  Math. \textbf{189} (2004), no.~2, 325--412.

\bibitem[HKR00]{HKR}
M.~Hopkins, N.~Kuhn, and D.~Ravenel, \emph{Generalized group characters and
  complex oriented cohomology theories}, J. Amer. Math. Soc. \textbf{13}
  (2000).

\bibitem[Hop02]{HopkinsICM2002}
M.~Hopkins, \emph{Algebraic topology and modular forms}, Proceedings of the ICM
  \textbf{1} (2002).

\bibitem[HS20]{ZhenMatt}
Z.~Huan and M.~Spong, \emph{Twisted quasi-elliptic cohomology and twisted
  equivariant elliptic cohomology}, preprint (2020).

\bibitem[HST10]{HST}
H.~Hohnhold, S.~Stolz, and P.~Teichner, \emph{From minimal geodesics to super
  symmetric field theories}, CRM Proceedings and Lecture Notes \textbf{50}
  (2010).

\bibitem[Hua18]{Huan}
Z.~Huan, \emph{Quasi-elliptic cohomology {I}}, Advances in Mathematics
  \textbf{337} (2018).

\bibitem[JF20]{TheoTMM}
T.~Johnson-Freyd, \emph{Topological {Mathieu} moonshine}, preprint (2020).

\bibitem[JFS17]{theoclaudia}
T.~Johnson-Freyd and C.~Scheimbauer, \emph{(op)lax natural transformations,
  twisted quantum field theories, and even higher€ morita categories},
  Advances in Mathematics \textbf{307} (2017), 147--223.

\bibitem[Jon87]{Jones}
V.~Jones, \emph{Hecke algebra representations of braid groups and link
  polynomials}, Annals of Mathematics \textbf{126} (1987), no.~2.

\bibitem[JP90]{JonesPetrack}
J.~Jones and S.~Petrack, \emph{The fixed point theorem in equivariant
  cohomology}, Trans. Amer. Math. Soc. \textbf{322} (1990).

\bibitem[Kit09a]{Kitaev}
A.~Kitaev, \emph{Periodic table for topological insulators and
  superconductors}, AIP Conf. Proc. \textbf{1134} (2009).

\bibitem[Kit09b]{Kitchloo}
N.~Kitchloo, \emph{Dominant {$K$}-theory and integrable highest weight
  representations of {K}ac-{M}oody groups}, Adv. Math. \textbf{221} (2009),
  no.~4, 1191--1226.

\bibitem[Kit14]{KitchlooII}
\bysame, \emph{Quantization of the modular functor and equivariant elliptic
  cohomology}, Arxiv preprint (2014).

\bibitem[KLW22]{KLW}
P.~Kristel, M.~Ludewig, and K.~Waldorf, \emph{A representation of the string
  2-group}, arxiv preprint (2022).

\bibitem[KM13]{KottkeMelrose1}
C.~Kottke and R.~Melrose, \emph{Equivalence of string and fusion loop-spin
  structures}, preprint (2013).

\bibitem[KM15]{KottkeMelrose2}
\bysame, \emph{Loop-fusion cohomology and transgression}, Math. Res. Lett.
  \textbf{22} (2015).

\bibitem[KRW19]{KRWell}
S.~Kumar, R.~Rim\'anyi, and A.~Weber, \emph{Elliptic classes of schubert
  varieties}, Mathematische Annalen (2019), 1--26.

\bibitem[KS21]{KontsevichSegal}
M.~Kontsevich and G.~B. Segal, \emph{Wick rotation and the positivity of energy
  in quantum field theory}, preprint (2021).

\bibitem[Kui65]{Kuiper}
N.~Kuiper, \emph{The homotopy type of the unitary group of {Hilbert} space},
  Topology \textbf{3} (1965).

\bibitem[Lic63]{Lichnerowicz}
A.~Lichnerowicz, \emph{Spineurs harmoniques}, C. R. Acad. Sci. Paris
  \textbf{S\'erie A} (1963), no.~257.

\bibitem[Lin22]{TMFmoon}
Y.-H. Lin, \emph{Topological modularity of monstrous moonshine}, preprint
  (2022).

\bibitem[Liu96]{Liu2}
Kefeng Liu, \emph{On elliptic genera and theta-functions}, Topology \textbf{35}
  (1996), no.~3, 617--640.

\bibitem[LP23]{LinPei}
Y.-H. Lin and D.~Pei, \emph{Holomorphic {CFTs} and topological modular forms},
  Communications in Mathematical Physics \textbf{401} (2023).

\bibitem[LS21]{LudewigStoffel}
M.~Ludewig and A.~Stoffel, \emph{A framework for geometric field theories and
  their classification in dimension one}, Symmetry, Integrability and Geometry:
  Methods and Applications (SIGMA) \textbf{17} (2021).

\bibitem[Lud22]{Matthias1}
M.~Ludewig, \emph{Categories of {Lagrangian} correspondences and fermionic
  functorial field theory}, preprint (2022).

\bibitem[Lud23]{Ludewigspinor}
\bysame, \emph{The spinor bundle on loop space}, arxiv preprint (2023).

\bibitem[Lue22]{Kiran}
K.~Luecke, \emph{Completed {K}-theory and equivariant elliptic cohomology},
  Advances in Mathematics \textbf{410} (2022).

\bibitem[Lur09a]{Lurie_cob}
J.~Lurie, \emph{On the classification of topological field theories}, Current
  Developments in Mathematics (2009).

\bibitem[Lur09b]{Lurie_Elliptic}
J.~Lurie, \emph{A survey of elliptic cohomology}, Algebraic Topology (N.~Baas,
  E.~Friedlander, J.~Bj\"orn, and P.~{O}st\ae{r}, eds.), vol.~4, Springer
  Berlin Heidelberg, 2009.

\bibitem[Lur19]{LurieIII}
J.~Lurie, \emph{Elliptic cohomology {III}: {T}empered cohomology}.

\bibitem[LW23]{LudewigWaldorf2grp}
M.~Ludewig and K.~Waldorf, \emph{Lie 2-groups from loop group extensions},
  arxiv preprint (2023).

\bibitem[Moo03]{MooreKtheory}
G.~Moore, \emph{{K theory from a physical perspective}}, {Symposium on
  Topology, Geometry and Quantum Field Theory (Segalfest)}, 4 2003,
  pp.~194--234.

\bibitem[Mor06]{MoravaHKR}
J.~Morava, \emph{H{KR} characters and higher twisted sectors}, Gromov-{W}itten
  theory of spin curves and orbifolds, Contemp. Math., vol. 403, Amer. Math.
  Soc., Providence, RI, 2006, pp.~143--152.

\bibitem[Mor09]{Moravamoon}
\bysame, \emph{Moonshine elements in elliptic cohomology}, Groups and
  symmetries, CRM Proc. Lecture Notes, vol.~47, Amer. Math. Soc., Providence,
  RI, 2009, pp.~247--257.

\bibitem[MQ86]{MathaiQuillen}
V.~Mathai and D.~Quillen, \emph{Superconnections, {Thom} classes and
  equivariant differential forms}, Topology \textbf{25} (1986).

\bibitem[MS67]{McKeanSinger}
H.P. McKean and I.M. Singer, \emph{Curvature and the eigenvalues of the
  {Laplacian}}, J. Differential Geometry \textbf{1} (1967).

\bibitem[PS86]{PressleySegal}
A.~Pressley and G.~Segal, \emph{Loop groups}, Oxford Mathematical Monographs,
  The Clarendon Press, Oxford University Press, New York, 1986, Oxford Science
  Publications.

\bibitem[Ros01]{RosuEquivariant}
I.~Rosu, \emph{Equivariant elliptic cohomology and rigidity}, Amer. J. Math.
  \textbf{123} (2001), no.~4, 647--677.

\bibitem[Ros03]{RosuDelocalized}
\bysame, \emph{Equivariant {$K$}-theory and equivariant cohomology}, Math. Z.
  \textbf{243} (2003), no.~3, 423--448, With an appendix by Allen Knutson and
  Rosu.

\bibitem[RSVZ22]{RSVZ}
R.~Rim\'anyi, A.~Smirnov, A.~Varchenko, and Z.~Zhou, \emph{Three-dimensional
  mirror symmetry and elliptic stable envelopes}, International Mathematics
  Research Notices \textbf{2022} (2022).

\bibitem[RTV17]{RTV}
R.~Rim\'anyi, V.~Tarasov, and A.~Varchenko, \emph{Elliptic and {K}-theoretic
  stable envelopes and {N}ewton polytopes}, Arxiv preprint (2017).

\bibitem[RW20]{RWschubert}
R.~Rim\'anyi and A.~Weber, \emph{Elliptic classes of {Schubert} varieties via
  {Bott--Samelson} resolution}, Journal of Topology \textbf{13} (2020).

\bibitem[RW22]{RWell}
\bysame, \emph{Elliptic classes on langlands dual flag varieties},
  Communications in Contemporary Mathematics \textbf{24} (2022), no.~01,
  2150014.

\bibitem[Seg88]{Segal_Elliptic}
G.~Segal, \emph{Elliptic cohomology}, S\'eminaire N. Bourbaki \textbf{695}
  (1988).

\bibitem[Seg01]{SegalDbrane}
\bysame, \emph{Topological structures in string theory}, Philosophical
  Transactions of the Royal Society of London. Series A: Mathematical, Physical
  and Engineering Sciences \textbf{359} (2001).

\bibitem[Seg04]{SegalCFT}
\bysame, \emph{The definition of conformal field theory}, Topology, geometry
  and quantum field theory, London Math. Soc. Lecture Note Ser., vol. 308,
  Cambridge Univ. Press, Cambridge, 2004, pp.~421--577.

\bibitem[Seg07]{SegalElliptic}
\bysame, \emph{What is an elliptic object?}, Elliptic cohomology, London Math.
  Soc. Lecture Note Ser., vol. 342, Cambridge Univ. Press, Cambridge, 2007,
  pp.~306--317.

\bibitem[SP11]{SchommerPries}
C.~Schommer-Pries, \emph{Central extensions of smooth 2-groups and a
  finite-dimensional string 2-group}, Geometry and Topology \textbf{15} (2011).

\bibitem[ST04]{ST04}
S.~Stolz and P.~Teichner, \emph{What is an elliptic object?}, Topology,
  geometry and quantum field theory, London Math. Soc. LNS 308, Cambridge Univ.
  Press (2004), 247--343.

\bibitem[ST05]{ST_spinors}
\bysame, \emph{The spinor bundle on loop space}, {MPIM} preprint (2005).

\bibitem[ST11]{ST11}
\bysame, \emph{Supersymmetric field theories and generalized cohomology},
  Mathematical Foundations of Quantum Field and Perturbative String Theory ({B.
  Jur{\v c}o, H. Sati, U. Schreiber}, ed.), Proceedings of Symposia in Pure
  Mathematics, 2011.

\bibitem[Sto96]{stolz_conj}
S.~Stolz, \emph{A conjecture concerning positive ricci curvature and the witten
  genus}, Math. Ann. \textbf{304} (1996), 785--800.

\bibitem[Sto13]{Stolzclass}
\bysame, \emph{Equivariant de {Rham} cohomology and gauged field theories},
  Course notes (2013).

\bibitem[Sto19]{Stoffeltwists}
A.~Stoffel, \emph{{Supersymmetric field theories from twisted vector bundles}},
  Commun. Math. Phys. \textbf{367} (2019), 417--453.

\bibitem[Sus07]{Susskind}
L.~Susskind, \emph{The anthropic landscape of string theory}, Universe or
  Multiverse? (B.~Carr, ed.), Cambridge University Press, 2007.

\bibitem[Tac21]{Tachikawa}
Y.~Tachikawa, \emph{{Topological modular forms and the absence of a heterotic
  global anomaly}}, Progress of Theoretical and Experimental Physics
  \textbf{2022} (2021), no.~4.

\bibitem[TY23a]{TachikawaYamashita2}
Y.~Tachikawa and M.~Yamashita, \emph{Anderson self-duality of topological
  modular forms, its differential-geometric manifestations, and vertex operator
  algebras}, preprint (2023).

\bibitem[TY23b]{TachikawaYamashita}
\bysame, \emph{Topological modular forms and the absence of all heterotic
  global anomalies}, Communications in Mathematical Physics \textbf{402}
  (2023), 1--36.

\bibitem[TYY23]{TYY}
Y.~Tachikawa, M.~Yamashita, and K.~Yonekura, \emph{Remarks on mod-2 elliptic
  genus}, preprint (2023).

\bibitem[Ulr21]{Ulrickson}
P.~Ulrickson, \emph{{Supersymmetric Euclidean field theories and K-theory}}, J.
  Geom. Phys. \textbf{161} (2021).

\bibitem[Wal13]{Waldorfstring}
K.~Waldorf, \emph{String connections and {Chern--Simons} theory}, Transactions
  of the American Mathematical Society \textbf{365} (2013).

\bibitem[Was98]{Wassermann}
A.~Wassermann, \emph{Operator algebras and conformal field theory {III}},
  Invent. math. \textbf{133} (1998).

\bibitem[Wit82a]{Wittenbreaking}
E.~Witten, \emph{Constraints on supersymmetry breaking}, Nuclear Physics B
  \textbf{202} (1982), no.~2, 253--316.

\bibitem[Wit82b]{susymorse}
\bysame, \emph{Supersymmetry and {Morse} theory}, Journal of Differential
  Geometry \textbf{17} (1982), 661--692.

\bibitem[Wit86]{WittenICM}
\bysame, \emph{Physics and geometry}, Proceedings of the International Congress
  of Mathematicians (1986).

\bibitem[Wit87]{Witten_Elliptic}
\bysame, \emph{Elliptic genera and quantum field theory}, Commun. Math. Phys.
  \textbf{109} (1987), 525--536.

\bibitem[Wit88a]{Witten_Dirac}
\bysame, \emph{The index of the {Dirac} operator in loop space}, Elliptic
  Curves and Modular Forms in Algebraic Topology (P.~Landweber, ed.), Springer
  Berlin Heidelberg, 1988.

\bibitem[Wit88b]{WittenTQFT}
\bysame, \emph{{Topological quantum field theory}}, Communications in
  Mathematical Physics \textbf{117} (1988), no.~3.

\bibitem[Wit89]{WittenJones}
\bysame, \emph{Quantum field theory and the {Jones} polynomial}, Commun. Math.
  Phys. \textbf{121} (1989).

\bibitem[Wit91]{Wittentoptwist}
\bysame, \emph{Introduction to cohomological field theory}, International
  Journal of Modern Physics A \textbf{6} (1991), no.~6.

\bibitem[Wit99]{Witteninstrings}
\bysame, \emph{Index of {Dirac} operators}, {Quantum Fields and Strings: {A}
  Course for Mathematicians, Volume 1} (P.~Deligne, P.~Etingof, D.~Freed,
  L.~Jeffrey, D.~Kazhdan, J.~Morgan, D.~Morrison, and E.~Witten, eds.),
  American Mathematical Society, 1999.

\bibitem[Wit19]{WittenRiemann}
\bysame, \emph{{Notes On Super Riemann Surfaces And Their Moduli}}, Pure Appl.
  Math. Quart. \textbf{15} (2019), no.~1.

\bibitem[Wu05]{Wu}
S.~Wu, \emph{{Mathai--Quillen} formalism}, Encyclopedia of Mathematical Physics
  (2005).

\bibitem[YZ17]{YangZhao}
Y.~Yang and G.~Zhao, \emph{Quiver varieties and elliptic quantum groups}, Arxiv
  preprint (2017).

\bibitem[Zag86]{Zagiermodular}
D.~Zagier, \emph{Note on the {Landweber-Stong} elliptic genus}, Elliptic curves
  and modular forms in Alg. Top. (1986).

\bibitem[ZZ15]{ZhaoZhong}
G.~Zhao and C.~Zhong, \emph{Elliptic affine hecke algebras and their
  representations}, Advances in Mathematics \textbf{395} (2015).

\end{thebibliography}

\end{document}